\pdfoutput=1
\RequirePackage{ifpdf}
\ifpdf 
\documentclass[pdftex]{sigma}
\else
\documentclass{sigma}
\fi

\numberwithin{equation}{section}

\newtheorem{Theorem}{Theorem}[section]
\newtheorem{Corollary}[Theorem]{Corollary}
\newtheorem{Lemma}[Theorem]{Lemma}
\newtheorem{Proposition}[Theorem]{Proposition}
\newtheorem{Fact}[Theorem]{Fact}

\theoremstyle{definition}
\newtheorem{Definition}[Theorem]{Definition}

\newtheorem{Example}[Theorem]{Example}
\newtheorem{Remark}[Theorem]{Remark}

\usepackage{mathrsfs}
\usepackage[all]{xy}
\usepackage{enumitem}
\usepackage{mathtools}

\newcommand{\spn}[1]{\mathrm{span}_{\CC}\{#1\}}
\newcommand{\set}[1]{\left\{#1 \right\}}

\newcommand{\Coker}{\mathrm{Coker}}
\newcommand{\Ker}{\mathrm{Ker}}

\newcommand{\Hom}{\mathrm{Hom}}
\DeclareMathOperator{\End}{End}

\DeclareMathOperator{\id}{id}

\newcommand{\BasicRing}[1]{\mathbb{#1}}
\newcommand{\ZZ}{\BasicRing{Z}}
\newcommand{\NN}{\BasicRing{N}}

\newcommand{\RR}{\BasicRing{R}}
\newcommand{\CC}{\BasicRing{C}}

\newcommand{\Mod}{\mathrm{Mod}}
\newcommand{\Ann}{\mathrm{Ann}}
\newcommand{\PIdeg}{\mathrm{PI.deg}}
\newcommand{\Ext}{\mathrm{Ext}}

\newcommand{\alg}[1]{\mathcal{#1}}

\DeclareMathOperator{\gr}{gr}

\DeclareMathOperator{\rank}{rank}
\newcommand{\Ad}{\mathrm{Ad}}
\newcommand{\ad}{\mathrm{ad}}

\newcommand{\classicalg}[1]{\mathfrak{#1}}

\newcommand{\lieC}[1]{\mathfrak{#1}}
\newcommand{\lie}[1]{\mathfrak{#1}}

\newcommand{\univ}[1]{\mathcal{U}(\lieC{#1})}

\newcommand{\univcent}[1]{\mathcal{Z}(\lieC{#1})}
\newcommand{\ind}[2]{\mathrm{ind}^{\lieC{#1}}_{\lieC{#2}}}

\newcommand{\BGGcat}[2]{\mathcal{O}_{\lie{#2}}^{\lie{#1}}}

\newcommand{\Dsheaf}[2]{\mathscr{#1}_{#2}}
\newcommand{\ntDsheaf}[1]{\Dsheaf{D}{#1}}

\newcommand{\zuck}[2]{\Gamma^{#1}_{#2}}
\newcommand{\Dzuck}[3]{R^{#3}\zuck{#1}{#2}}

\newcommand{\regular}{\mathcal{O}}
\newcommand{\rring}[1]{\regular(#1)}

\let\Re\relax
\let\Im\relax
\DeclareMathOperator{\Re}{\mathrm{Re}}
\DeclareMathOperator{\Im}{\mathrm{Im}}

\newcommand{\rhog}[1]{\rho(#1)}
\DeclareMathOperator{\lattice}{\mathcal{L}}
\DeclareMathOperator{\Irr}{\mathrm{Irr}}
\DeclareMathOperator{\Co}{\mathrm{Co}}
\newcommand{\fdual}[1]{{#1}^{d}}

\newcommand{\lieSL}{\classicalg{sl}}
\newcommand{\lieGL}{\classicalg{gl}}
\newcommand{\lieU}{\classicalg{u}}
\newcommand{\lieSp}{\classicalg{sp}}
\newcommand{\lieSO}{\classicalg{so}}
\newcommand{\lieG}{\classicalg{g}_2}

\newcommand{\Sn}{\mathfrak{S}}
\DeclareMathOperator{\PI}{PI}
\DeclareMathOperator{\sgn}{sgn}

\newcommand{\lmod}[3]{\mathcal{L}^{\lie{#1}}_{\lie{#2}, #3}}

\begin{document}
\allowdisplaybreaks

\newcommand{\arXivNumber}{2410.17125}

\renewcommand{\PaperNumber}{095}

\FirstPageHeading

\ShortArticleName{Construction of Irreducible $\mathcal{U}(\mathfrak{g})^{G'}$-Modules and Discretely Decomposable Restrictions}

\ArticleName{Construction of Irreducible $\boldsymbol{\mathcal{U}(\mathfrak{g})^{G'}}$-Modules\\ and Discretely Decomposable Restrictions}

\Author{Masatoshi KITAGAWA}

\AuthorNameForHeading{M.~Kitagawa}

\Address{Institute of Mathematics for Industry, Kyushu University, \\ 744~Motooka, Fukuoka-shi, 819-0395, Fukuoka, Japan}
\Email{\mail{m-kitagawa@imi.kyushu-u.ac.jp}}

\ArticleDates{Received March 31, 2025, in final form October 27, 2025; Published online November 08, 2025}

\Abstract{In this paper, we study the irreducibility of $\mathcal{U}(\mathfrak{g})^{G'}$-modules on the spaces of intertwining operators in the branching problem of reductive Lie algebras, and construct a~family of finite-dimensional irreducible $\mathcal{U}(\mathfrak{g})^{G'}$-modules using the Zuckerman derived functors. We provide criteria for the irreducibility of $\mathcal{U}(\mathfrak{g})^{G'}$-modules in the cases of generalized Verma modules, cohomologically induced modules, and discrete series representations. We treat only discrete decomposable restrictions with certain dominance conditions (quasi-abelian and in the good range). To describe the $\mathcal{U}(\mathfrak{g})^{G'}$-modules, we give branching laws of cohomologically induced modules using ones of generalized Verma modules when $K'$ acts on $K/L_K$ transitively.}


\Classification{22E46; 22E47}

\section{Introduction}

The purpose of this paper is to study an algebraic structure on the spaces of intertwining operators in the branching problem of reductive Lie algebras.
We consider the irreducibility of $\univ{g}^{G'}$-modules on the spaces, and construct a family of finite-dimensional irreducible $\univ{g}^{G'}$-modules using the Zuckerman derived functors.
This work is based on a part of the author's Ph.D.\ thesis.

Let $G$ be a connected reductive algebraic group over $\CC$ and $G'$ a connected reductive subgroup of $G$.
For an irreducible $\lie{g}$-module $V$ and an irreducible $\lie{g'}$-submodule $V'$ of $V|_{\lie{g'}}$, the algebra~$\univ{g}^{G'}$ of $G'$-invariants in the universal enveloping algebra $\univ{g}$ acts on \smash{$\Hom_{\lie{g'}}\bigl(V', V|_{\lie{g'}}\bigr)$}.
It is a classical result that if \smash{$V|_{\lie{g'}}$} is completely reducible and locally finite, then the \smash{$\univ{g}^{G'}$}-module is irreducible (see, e.g., \cite[\S4.2]{GoWa09}).
This is corresponding to the case of the branching problem of real Lie groups and compact subgroups.
In this paper, we essentially treat the case of non-compact subgroups.

In general, the $\univ{g}^{G'}$-module may not be irreducible.
Some branching laws in the context of the theta lifting give typical examples (see Example~\ref{ex:NonIrreducible}).
On the other hand, we have shown in~\cite[Theorem~7.18 and its corollaries]{Ki20-2} that the module has finite length if $V$ and $V'$ are realized as holonomic $\ntDsheaf{}$-modules on the full flag varieties.
For example, $V$ and $V'$ are objects in the BGG category $\BGGcat{}{}$, or Harish-Chandra modules.
Moreover, we have treated a similar $\univ{g}^{G'}$-module $\Hom_{\lie{g'}}\bigl(V|_{\lie{g'}}, V'\bigr)$ and its cohomological version.
We concentrate on the $\univ{g}^{G'}$-module $\Hom_{\lie{g'}}\bigl(V', V|_{\lie{g'}}\bigr)$ in this paper.

We will consider generalized Verma modules, cohomologically parabolically induced modules and discrete series representations.
We give criteria for the irreducibility of the $\univ{g}^{G'}$-module for such modules.

We shall explain why we study $\univ{g}^{G'}$-modules.
Recently, explicit realizations of intertwining operators (symmetry breaking operators) in the branching problem are computed by many researchers, e.g.,~\cite{KoPe15_f_method_I,KoPe15_f_method_II,KoSp15_symmetry_breaking}.
Under some good setting, the operators are represented by special functions.
The computations use concrete realizations of representations and operators, and the operators are obtained from good differential equations and good geometric structures.

To understand the goodness of intertwining operators, we consider algebraic structures on the $\Hom$ space.
For example, if $\Hom_{\lie{g'}}\bigl(V'_1, V_1|_{\lie{g'}}\bigr)\simeq \Hom_{\lie{g'}}\bigl(V'_2, V_2|_{\lie{g'}}\bigr)$,
then we may expect that there is some relation between intertwining operators in the two different $\Hom$ spaces.
Alternatively, if $\univ{g}^{G'}$ acts on the $\Hom$ space irreducibly, then we may consider that the algebra $\univ{g}^{G'}$ controls the branching law well.
In fact, we will show in Theorem~\ref{thm:PIdegMultiplicity} that the supremum of the multiplicities is equal to a ring-theoretic invariant of $\univ{g}^{G'}$, called the PI degree.

We shall state the main results in this paper.
Let $G$ be a connected semisimple algebraic group over $\CC$ and $G'$ a connected reductive subgroup of $G$.
Let $\lie{q}$ be a parabolic subalgebra of $\lie{g}$ such that $\lie{q'} \coloneq \lie{g'}\cap \lie{q}$ is a parabolic subalgebra of $\lie{g'}$.
Suppose that $\lie{q}$ has a Levi decomposition $\lie{q} = \lie{l}\oplus \lie{u}$ such that $\lie{q'}$ has a Levi decomposition $\lie{q'} = (\lie{l}\cap \lie{g'}) \oplus (\lie{u} \cap \lie{g'})$.
Set $\lie{l'}\coloneq \lie{l}\cap \lie{g'}$ and $\lie{u}'\coloneq \lie{u} \cap \lie{g'}$.
We consider a generalized Verma module $\ind{g}{q}(F)$ induced from a finite-dimensional irreducible $\lie{l}$-module $F$, letting $\lie{u}$ act on $F$ trivially.

In general, the restriction $\ind{g}{q}(F)|_{\lie{g'}}$ may be complicated.
To study the space $\Hom_{\lie{g'}}\bigl(V', V|_{\lie{g'}}\bigr)$, we need admissibility and complete reducibility
and assume two conditions, i.e., weakly $\lie{g'}$-compatible and quasi-abelian.
See Definitions~\ref{def:WeaklyCompatible} and~\ref{def:QuasiAbelian}.
If $\lie{q}$ is weakly $\lie{g'}$-compatible, then~$\ind{g}{q}(F)|_{\lie{g'}}$ is discretely decomposable and $\lie{l'}$-admissible (or equivalently $\lie{g'}$-admissible).
If,~in~addition, $\lie{q}$ is quasi-abelian with respect to $\lie{g'}$, then $\ind{g}{q}(F)|_{\lie{g'}}$ with $F$ in the good range is completely reducible.

The weak $\lie{g'}$-compatibility is a generalization of the $\lie{g'}$-compatibility defined in~\cite{Ko12_generalized_Verma}.
If $(\lie{g}, \lie{g'})$ is a symmetric pair corresponding to an involution $\theta$ of $\lie{g}$, the notions are equivalent to the $\theta$-stability.
The notion of quasi-abelian parabolic subalgebras is defined in~\cite[p.~109]{HMSW87} for symmetric $(\lie{g}, \lie{g'})$.
Our definition is a straightforward generalization of theirs.
For example, if \smash{${\bigl[\lie{u}\cap (\lie{g'})^\perp, \lie{u'}\bigr] = 0}$}, then $\lie{q}$ is quasi-abelian.

Fix Cartan subalgebras $\lie{t'}$ of $\lie{l'}$ and $\lie{t}$ of $\lie{l}$ such that $\lie{t'}\subset \lie{t}$.
We denote by $\Delta(W, \lie{t})$ the set of all non-zero weights in a $\lie{t}$-module $W$,
and by $\rhog{\lie{u}}$ half the sum of roots in $\Delta(\lie{u}, \lie{t})$.
Suppose that $\lie{q}$ is weakly $\lie{g'}$-compatible and quasi-abelian with respect to $\lie{g'}$.

The following theorem will be proved in Theorem~\ref{thm:QuasiAbelianGoodRange} and Corollary~\ref{cor:IrreducibleVermaQuasiAbelian}.
Without the quasi-abelian property, criteria for the irreducibility of $\univ{g}^{G'}$-modules are given in corollaries of Theorem~\ref{thm:InjectivePhi}.

\begin{Theorem}\label{intro:GeneralizedVerma}
	Let $F$ be a finite-dimensional irreducible $\lie{l}$-module with infinitesimal character $\lambda \in \lie{t}^*$ in the good range, i.e.,
	\begin{align*}
		\Re(\lambda + \rhog{\lie{u}}, \alpha) < 0, \qquad \forall \alpha \in \Delta^+(\lie{u}, \lie{t}).
	\end{align*}
	Then
	\begin{enumerate}\itemsep=0pt
		\item[$1.$] $\ind{g}{q}(F)|_{\lie{g'}}$ is a direct sum of irreducible $\lie{g'}$-modules of the form $\ind{g'}{q'}(F')$ with finite-dimensional irreducible $\lie{l'}$-module $F'$ in the good range.
		\item[$2.$] The $\univ{g}^{G'}$-module \smash{$\Hom_{\lie{g'}}\bigl(\ind{g'}{q'}(F'), \ind{g}{q}(F)|_{\lie{g'}}\bigr)$} is zero or irreducible for any irreducible finite-dimensional $\lie{l'}$-module $F'$,
	\end{enumerate}
\end{Theorem}

\begin{Remark}
	The multiplicity in $\ind{g}{q}(F)|_{\lie{g'}}$ can be computed from the irreducible decomposition of a locally finite $\lie{l'}$-module.
	See Proposition~\ref{prop:FiltrationRestriction}.
\end{Remark}

The first assertion of Theorem~\ref{intro:GeneralizedVerma} was proved in~\cite[Lemma~3.1]{HMSW87} when $F|_{\lie{l'}}$ is irreducible and $(\lie{g}, \lie{g'})$ is a symmetric pair.
If $(\lie{g}, \lie{l})$ is a symmetric pair, then $\ind{g}{q}(F)$ is isomorphic to the underlying Harish-Chandra module of a holomorphic discrete series representation.
In this case, the first assertion of Theorem~\ref{intro:GeneralizedVerma} was proved in~\cite[Theorem~7.4]{Ko98} and~\cite{Ko12_generalized_Verma}.

Using the Zuckerman derived functor, we obtain a similar result to Theorem~\ref{intro:GeneralizedVerma} for cohomologically induced modules.
Let $\theta$ be an involution of $\lie{g}$ fixing $\lie{g'}$.
Let $K$ be a connected reductive algebraic group with the Lie algebra $\lie{g}^\theta$ such that $(\lie{g}, K)$ is a pair (see Definition~\ref{def:pair}).
Suppose that $\lie{q}$ and $\lie{l}$ are $\theta$-stable.
Write $K'$, $L_K$ and $L'_K$ for the analytic subgroups with the Lie algebras $\lie{k}\cap \lie{g'}$, $\lie{k}\cap \lie{l}$ and $\lie{k'}\cap \lie{l}$, respectively.
Then $L_K$ (resp.\ $L'_K$) is a Levi subgroup of $K$ (resp.\ $K'$).

For a finite-dimensional irreducible $(\lie{l}, L_K)$-module $F$, we set
\begin{align*}
	\lmod{g}{q}{S}(F) \coloneq \Dzuck{K}{L_K}{S}\bigl(\ind{g}{q}(F)\bigr),
\end{align*}
where \smash{$\Dzuck{K}{L_K}{S}$} is the $S$-th Zuckerman derived functor and $S\coloneq \dim(\lie{u}\cap \lie{k})$.
Then \smash{$\lmod{g}{q}{S}(F)$} is an irreducible $(\lie{g}, K)$-module if $F$ is in the good range.
See~\cite{KnVo95_cohomological_induction} for the cohomologically induced modules.

The following theorem is proved in Theorem~\ref{thm:DerivedFunctorModuleQuasiAbelian}.
See Theorems~\ref{thm:IrreducibleHolomorphicDisc} and~\ref{thm:DiscreteSeries} for (holomorphic) discrete series representations.

\begin{Theorem}\label{intro:CohomologicalInduction}
	Assume that $K'$ acts on $K/L_K$ transitively.
	Let $F$ be a finite-dimensional irreducible $(\lie{l}, L_K)$-module in the good range.
	\begin{enumerate}\itemsep=0pt
		\item[$1.$] $\lmod{g}{q}{S}(F)|_{\lie{g'}, K'}$ is a direct sum of irreducible $(\lie{g'}, K')$-modules of the form $\lmod{g'}{q'}{S}(F')$ with finite-dimensional irreducible $\bigl(\lie{l'}, L'_K\bigr)$-module $F'$ in the good range.
		\item[$2.$] For any finite-dimensional irreducible $\bigl(\lie{l'}, L'_K\bigr)$-module $F'$ in the good range, \smash{$\Dzuck{K'}{L'_K}{S}$} induces an isomorphism
		\begin{align*}
			\Hom_{\lie{g}', L'_K}\bigl(\ind{g'}{q'}(F'), \ind{g}{q}(F)|_{\lie{g'}, L'_K}\bigr) \xrightarrow{\simeq} \Hom_{\lie{g}', K'}\bigl(\lmod{g'}{q'}{S}(F'), \lmod{g}{q}{S}(F)|_{\lie{g'}, K'}\bigr)
		\end{align*}
		of $\univ{g}^{G'}\!$-modules, and
		each $\univ{g}^{G'}\!$-module \smash{$\Hom_{\lie{g}', K'}\bigl(\lmod{g'}{q'}{S}(F'), \lmod{g}{q}{S}(F)|_{\lie{g'}, K'}\bigr)$} is irreducible or zero.
	\end{enumerate}
\end{Theorem}

Discretely decomposable restrictions of cohomologically induced modules ($A_{\lie{q}}(\lambda)$, in particular)
are studied by many researchers.
T.~Kobayashi initiated the theory of discretely decomposable restrictions and gave examples of branching laws of cohomologically induced modules in~\mbox{\cite{Ko93_discretely_decomposable,Ko94,Ko98_discrete_decomposable_2,Ko98_discrete_decomposable_3}}.
Y.~Oshima gave in~\cite{Os13_phd,Os24} an abstract formula of the branching laws like generalized Blattner's formula using the theory of $\ntDsheaf{}$-modules, and computed all the branching laws of $A_{\lie{q}}(\lambda)$ with character $\lambda$ in the weakly fair range when $(\lie{g}, \lie{g'})$ is a symmetric pair.
See~\mbox{\cite[Corollaries 5.7 and 5.8, and Section 8]{Os24}} for the case in which $K'$ acts on $K/L_K$ transitively.
We refer the reader to~\cite{Ko24-recent} for the recent advances of the branching problem.

As an example of Theorem~\ref{intro:CohomologicalInduction}, we will deal with the case of discrete series representations and symmetric $(\lie{g}, \lie{g'})$ in Theorem~\ref{thm:DiscreteSeries}.
In the case, the branching laws were computed by several methods, e.g., transfer of $K$-types~\cite{GrWa00}, orbit method~\cite{DuVa10,Va16}, reproducing kernel~\cite{OrVa20,OrVa24}, the Radon--Penrose transform~\cite{Se13} (including $A_{\lie{q}}(\lambda)$).
Note that our method is purely algebraic.

To show Theorems~\ref{intro:GeneralizedVerma} and~\ref{intro:CohomologicalInduction}, we construct $\univ{g}^{G'}$-modules isomorphic to the $\univ{g}^{G'}$-modules \smash{$\Hom_{\lie{g}', L'_K}\bigl(\ind{g'}{q'}(F'), \ind{g}{q}(F)|_{\lie{g'}, L'_K}\bigr)$} using the Zuckerman derived functor.
Retain the above notation $\lie{g}, \lie{g'}, \lie{q}, \lie{q'}, \dots$.
Then $(\lie{g'\oplus g}, \Delta(G'))$ is a pair, where $\Delta\colon G'\rightarrow G'\times G$ is the diagonal embedding.
Set $S\coloneq \dim(\lie{u}')$.
Write $L'$ for the analytic subgroup of $G'$ with the Lie algebra $\lie{l'}$.

A $(\lie{g'\oplus g}, \Delta(G'))$-module is completely reducible as a $\Delta(G')$-module.
Hence the functor
\begin{align*}
	\Mod(\lie{g'\oplus g}, \Delta(G')) \ni V\mapsto V^{\Delta(G')} \in \Mod\bigl(\univ{g}^{G'}\bigr)
\end{align*}
is exact and sends irreducible objects to irreducible ones or zero.
See Section~\ref{subsect:GKmod}.
Therefore, what we should do is to construct irreducible $(\lie{g'\oplus g}, \Delta(G'))$-modules.
We will use classical methods to do so, e.g., the Zuckerman derived functors, the translation functors.

For a finite-dimensional $\lie{l'}$-module $F'$, we denote by \smash{$\BGGcat{g'}{q'}(F)$} the full subcategory of the relative BGG category \smash{$\BGGcat{g'}{q'}$} whose object $V$ satisfies that $V\otimes F'$ lifts to an $L'$-module.
The following theorem is proved in Theorem~\ref{thm:Exactness} and Corollary~\ref{cor:CategoryEquivalence}.

\begin{Theorem}\label{intro:Embedding}
	Let $F$ be a finite-dimensional irreducible $\lie{l}$-module with infinitesimal character $\lambda \in \lie{t}^*$ satisfying
	\begin{align*}
		\frac{2(\lambda+\rhog{\lie{u}}, \alpha)}{(\alpha, \alpha)} \not \in \set{0, 1, 2, \ldots}, \qquad \forall \alpha \in \Delta(\lie{u}, \lie{t}).
	\end{align*}
	Then the functor
	\begin{align*}
		\BGGcat{g'}{q'}(F) \ni M\mapsto \Dzuck{\Delta(G')}{\Delta(L')}{S}\bigl(M\otimes \ind{g}{q}(F)\bigr) \in \Mod(\lie{g'\oplus g}, \Delta(G'))
	\end{align*}
	is fully faithful and exact, and preserves submodules lattices.
	In particular, it sends irreducible objects to irreducible ones.
\end{Theorem}

If $(\lie{g}, \lie{k})$ is not a symmetric pair, a $(\lie{g}, \lie{k})$-module with some finiteness properties is called a~generalized Harish-Chandra module, and studied in~\cite{PeZu04_generalized_harish_chandra,PeZu04_generalized_harish_chandra2} and their sequels.
Cohomologically induced modules are treated in the papers.
In our setting, $(\lie{g}, \lie{k})$ is the pair $(\lie{g'\oplus g}, \Delta(\lie{g'}))$.
To~show Theorem~\ref{intro:Embedding}, we need this special structure.

\subsection*{Notation and convention}

In this paper, any algebra except Lie algebras is unital, associative and over $\CC$.
Any algebraic group is defined over $\CC$ and any Lie algebra is finite-dimensional.
For an algebra $\alg{A}$ (resp.\ a~pair $(\lie{g}, K)$), we denote by $\Mod(\alg{A})$ (resp.\ $\Mod(\lie{g}, K)$) the category of $\alg{A}$-modules (resp.\ \mbox{$(\lie{g}, K)$-modules}).

We express real Lie groups and their Lie algebras by Roman letters and corresponding German letters with subscript $(\cdot)_\RR$, and express complex Lie groups (or affine algebraic groups) and their Lie algebras by Roman letters and corresponding German letters, respectively.
Similarly, we express the complexification of a real Lie algebra by the same German letter as that of the real form without any subscript.
For example, the Lie algebras of real Lie groups $G_\RR$, $K_\RR$ and $H_\RR$ are denoted as $\lie{g}_\RR$, $\lie{k}_\RR$ and $\lie{h}_\RR$ with complexifications $\lie{g}$, $\lie{k}$ and $\lie{h}$, respectively.

For a $\lie{t}$-module $V$ of a commutative Lie algebra $\lie{t}$, we denote by $\Delta(V, \lie{t})$ the (multi)set of all non-zero weights in $V$.
We write $V_\alpha$ for the weight space of a weight $\alpha \in \lie{t}^*$.
We express half the sum of roots (or positive roots in some context) in $\Delta(V, \lie{t})$ by $\rhog{V}$.

For a $G$-set $X$ of a group $G$, we write $X^G$ for the subset of all $G$-invariant elements in $X$.
We use similar notation for the set of all vectors annihilated by a Lie algebra $\lie{g}$ (resp.\ an~element~$X$) as $V^{\lie{g}}$ \big(resp.\ $V^X$\big).
The coordinate ring of an affine variety $X$ is denoted by $\rring{X}$.

\section{Preliminary}

In this section, we recall several fundamental notions about $(\lie{g}, K)$-modules, generalized Verma modules and Zuckerman derived functors.

\subsection[(g,K)-module]{$\boldsymbol{(\lie{g}, K)}$-module}\label{subsect:GKmod}

We recall the notion of $(\lie{g}, K)$-modules.
Let $G$ be an affine algebraic group over $\CC$.
A linear action of $G$ on a vector space $V$ is said to be rational
if the action is locally finite and any finite-dimensional subrepresentation is a representation of an algebraic group.
In the case, $V$ is called a $G$-module.
If $G$ is reductive, any $G$-module is completely reducible.

Let us recall the definitions of pairs $(\lie{g}, K)$ and $(\lie{g}, K)$-modules.
In this paper, we treat only the case of connected reductive $K$.
We refer the reader to \cite[Chapter~I]{KnVo95_cohomological_induction} for the general case.

\begin{Definition}\label{def:pair}
	Let $\lie{g}$ be a finite-dimensional complex Lie algebra and $K$ a connected reductive algebraic group
	such that $\lie{k}$ is a subalgebra of $\lie{g}$.
	We say that $(\lie{g}, K)$ is a \emph{pair} if the adjoint action of $\lie{k}$ on $\lie{g}$ lifts to a rational action of $K$.
	A pair $(\lie{h}, L)$ is said to be a \emph{subpair} of $(\lie{g}, K)$ if $\lie{h}$ is a~subalgebra of $\lie{g}$
	and $L$ is a connected reductive subgroup of $K$.
\end{Definition}

\begin{Definition}\label{def:GKmod}
	Let $(\lie{g}, K)$ be a pair and $V$ a $\lie{g}$-module.
	We denote by $V_{\lie{k}}$ the sum of all finite-dimensional $\lie{k}$-submodules,
	and by $V_K$ the sum of all $\lie{k}$-submodules that lift to $K$-modules.
	We~say that $V$ is a \emph{$(\lie{g}, K)$-module} if $V = V_K$.
\end{Definition}

Clearly, $V_K$ is a $(\lie{g}, K)$-module.
Since we have assumed that $K$ is connected, the $K$-action on a $(\lie{g}, K)$-module $V$ is uniquely determined by the $\lie{g}$-action.

Let $(\lie{g}, K)$ be a pair and $V$ a $(\lie{g}, K)$-module.
Then $V|_K$ is completely reducible and $\univ{g}^K$ acts naturally on $V^K$.
Since $\lie{k}$ acts on $V^K$ trivially, the $\univ{g}^K$-action factors through $(\univ{g}/\univ{g}\lie{k})^K \!\simeq \univ{g}^K/(\univ{g}\lie{k})^K$.
It is well-known that if $V$ is irreducible, then the $\univ{g}^K$-module $V^K$ is irreducible or zero.
See~\cite[Theorem~4.2.1]{GoWa09}.

\begin{Proposition}\label{prop:UKaction}
	Let $(\lie{g}, K)$ be a pair.
	The functor
 \[
 \Mod(\lie{g}, K) \ni V \mapsto V^K \in \Mod\bigl((\univ{g}/\univ{g}\lie{k})^K\bigr)
 \]
 is exact and
	sends irreducible objects to irreducible ones or zero.
	In particular, the length of $(\univ{g}/\univ{g}\lie{k})^K$-module $V^K$ is less than or equal to that of $V$
	for any $V \in \Mod(\lie{g}, K)$.
\end{Proposition}

It may be hard to study \smash{$(\univ{g}/\univ{g}\lie{k})^K$}-modules directly.
Using Proposition~\ref{prop:UKaction}, we reduce many problems of \smash{$(\univ{g}/\univ{g}\lie{k})^K$}-modules
to those of $(\lie{g}, K)$-modules.
The following $(\lie{g'\oplus g}, \Delta(G'))$-module case is important for us.

Let $G$ be a connected affine algebraic group and $G'$ a connected reductive subgroup of $G$.
Then $(\lie{g' \oplus g}, \Delta(G'))$ is a pair, where $\Delta\colon G'\rightarrow G'\times G$ is the diagonal embedding.

\begin{Proposition}\label{prop:gg'G'}
	$(\univ{g'\oplus g}/\univ{g'\oplus g}\Delta(\lie{g'}))^{\Delta(G')}$
	is naturally isomorphic to $\univ{g}^{G'}$ as a $\CC$-algebra.
\end{Proposition}

\begin{proof}
	Identify $\univ{g'\oplus g}$ with $\univ{g'}\otimes \univ{g}$.
	Then the isomorphism is given by
	\begin{align*}
		\univ{g}^{G'} \ni X \mapsto 1\otimes X + \univ{g'\oplus g}\Delta(\lie{g'}) \in (\univ{g'\oplus g}/\univ{g'\oplus g}\Delta(\lie{g'}))^{\Delta(G')}.
	\end{align*}
	In face, the map is bijective by $\lie{g'}\oplus \lie{g} = \Delta(\lie{g'}) \oplus \lie{g}$ and the Poincar\'e--Birkhoff--Witt theorem,
	and the linear map is a homomorphism of $\CC$-algebras because it is given as the composition of the following natural maps:
	\begin{align*}
		\univ{g}^{G'} \rightarrow \univ{g'\oplus g}^{\Delta(G')} &\rightarrow \univ{g'\oplus g}^{\Delta(G')}/(\univ{g'\oplus g}\Delta(\lie{g'}))^{\Delta(G')} \\
		&\simeq (\univ{g'\oplus g}/\univ{g'\oplus g}\Delta(\lie{g'}))^{\Delta(G')}. \qedhere
 \tag*{\qed}
 \end{align*}
\renewcommand{\qed}{}
\end{proof}

The two algebras $\univ{g}^{G'}$ and $(\univ{g'\oplus g}/\univ{g'\oplus g}\Delta(\lie{g'}))^{\Delta(G')}$ act naturally on $\Hom_{\lie{g'}}(V', V)$ for a $\lie{g}$-module $V$ and a $\lie{g'}$-module $V'$.
The action of $\univ{g}^{G'}$ is give by
\begin{align*}
	(X\cdot \varphi)(\cdot) = X\cdot \varphi(\cdot), \qquad X \in \univ{g},\quad \varphi \in \Hom_{\lie{g'}}(V', V),
\end{align*}
and the action of $(\univ{g'\oplus g}/\univ{g'\oplus g}\Delta(\lie{g'}))^{\Delta(G')}$ is induced from that of the $\univ{g'\oplus g}^{\Delta(G')}$:
\begin{align*}
	&\biggl(\biggl(\sum_{i} A_i \otimes B_i\biggr) \cdot \varphi\biggr)(\cdot)= \sum_i B_i \varphi\big({}^tA_i \cdot\big), \\
	&\qquad \sum_i A_i \otimes B_i \in \univ{g'\oplus g}^{G'},\ \varphi \in \Hom_{\lie{g'}}(V', V)\biggr.
\end{align*}
It is easy to see that the two actions coincide under the isomorphism in Proposition~\ref{prop:gg'G'}.

\subsection{Translation functor}

Our purpose in this paper is to show the irreducibility of some $(\lie{g'\oplus g}, \Delta(G'))$-modules.
We will reduce the problem to the case of enough large infinitesimal characters.
To do so, we need the translation functor.
We refer the reader to~\cite[Chapter~VII]{KnVo95_cohomological_induction}.

Let $\lie{g}$ be a finite-dimensional complex reductive Lie algebra.
Fix a Cartan subalgebra $\lie{t}$ of $\lie{g}$ and denote by $W_G$ the Weyl group of $\lie{g}$.
We identify $\univcent{g}$ with $\univ{t}^{W_G}$ via the Harish-Chandra isomorphism.
Then characters of $\univcent{g}$ are parametrized by elements in $\lie{t}^*/W_G$.
When we say that $\chi \in \lie{t}^* / W_G$ is an infinitesimal character,
$\chi$ is identified with a homomorphism~${\univcent{g}\rightarrow \CC}$
and~$\Ker(\chi)$ is a maximal ideal of $\univcent{g}$.
For $\lambda \in \lie{t}^*$, we denote by $[\lambda]$ the corresponding infinitesimal character in $\lie{t}^*/W_G$,
and for a subcategory $\mathcal{C}$ of $\Mod(\lie{g})$, we denote by $\mathcal{C}_{[\lambda]}$ (or $\mathcal{C}_\lambda$) the full subcategory of $\mathcal{C}$ consisting of objects with the generalized infinitesimal character $[\lambda]$.

For a $\lie{g}$-module $V$ and an infinitesimal character $\chi \in \lie{t}^*/W_G$, set
\begin{align*}
	P_{\chi}(V) \coloneq \{v \in V \mid \Ker(\chi)^n v = 0, \, \exists n \in \NN\}.
\end{align*}
Note that $n$ can depend on $v$.
Then $P_{\chi}(V)$ is a $\lie{g}$-submodule of $V$, and called the \emph{primary component} with the infinitesimal character $\chi$.
When $V$ is the direct sum of all primary components, we say that the direct sum decomposition is the \emph{primary decomposition}.
The following fact is a~direct consequence of~\cite[Proposition~7.20]{KnVo95_cohomological_induction}.
For a $\lie{g}$-module $V$, we say that $V$ is locally finite $Z(\lie{g})$-module if $Z(\lie{g})v$ is finite-dimensional for any $v \in V$.

\begin{Fact}\label{fact:PrimaryDecomposition}
	Any locally $\univcent{g}$-finite $\lie{g}$-module has the primary decomposition.
\end{Fact}

We shall recall the translation functor.
Let $\mu \in \lie{t}^*$ be an algebraically integral weight of~$\lie{g}$.
Let~$F(\mu)$ be a finite-dimensional irreducible $\lie{g}$-module with the extreme weight $\mu$.
For a~$\lie{g}$-module~$V$ with an generalized infinitesimal character $[\lambda] \in \lie{t}^*/W_G$, set
\begin{align*}
	T_\lambda^{\lambda + \mu}(V) \coloneq P_{[\lambda + \mu]}(F(\mu)\otimes V).
\end{align*}
Then $P_{[\lambda + \mu]}(F(\mu)\otimes V)$ is a $\lie{g}$-module with the generalized infinitesimal character $[\lambda + \mu]$.

Note that $F(\mu)\otimes V$ is $\univcent{g}$-finite by Kostant's theorem~\cite[Theorem~7.133]{KnVo95_cohomological_induction},
and hence \smash{$T_\lambda^{\lambda + \mu}(V)$} is a direct summand of \smash{$F(\mu) \otimes V$}.
This implies that \smash{$T_\lambda^{\lambda + \mu}$} is an exact functor from \smash{$\Mod(\lie{g})_{[\lambda]}$} to \smash{$\Mod(\lie{g})_{[\lambda + \mu]}$}.
The functor is called the \emph{translation functor}.

The translation functor preserves the irreducibility under some assumption.
To state the facts, we shall recall the notion of integrally dominant weights.
Fix a Borel subalgebra $\lie{b}$ of $\lie{g}$ containing $\lie{t}$.
Let $\Delta^+$ denote the set of positive roots determined by $\lie{b}$.

\begin{Definition}\label{def:IntegrallyDominant}
	Let $\lambda\in \lie{t}$.
	We say that $\lambda$ regular if $(\lambda, \alpha) \neq 0$ for all $\alpha \in \Delta^+$.
	We say that $\lambda$ is \emph{integrally dominant} with respect to $\lie{b}$ if
	\begin{align*}
		\frac{2(\lambda, \alpha)}{(\alpha, \alpha)} \not \in \set{-1, -2, \ldots}, \qquad \forall \alpha \in \Delta^+.
	\end{align*}
	$\lambda$ is said to be \emph{integrally anti-dominant} with respect to $\lie{b}$ if
	\begin{align*}
		\frac{2(\lambda, \alpha)}{(\alpha, \alpha)} \not \in \set{1, 2, \ldots}, \qquad \forall \alpha \in \Delta^+.
	\end{align*}
\end{Definition}

It is well-known that \smash{$T_\lambda^{\lambda + \mu}$} is a good functor (e.g., isomorphism) on several good categories.
We only need results about $\lie{g}$-modules of finite length.
For the following fact, see~\cite[Corollary~7.209]{KnVo95_cohomological_induction} and its proof.
Let $\Mod_{fl}(\lie{g})$ denote the category of $\lie{g}$-modules of finite length.

\begin{Fact}\label{fact:TranslationEquvalence}
	Let $\lambda \in \lie{t}^*$ be a regular integrally dominant weight and $\mu \in \lie{t}^*$ an algebraically integral weight.
	Suppose that $\lambda + \mu$ is regular integrally dominant.
	Then \smash{$T^{\lambda+\mu}_{\lambda}$} gives an equivalence of categories from $\Mod_{fl}(\lie{g})_{\lambda}$ to $\Mod_{fl}(\lie{g})_{\lambda + \mu}$, and a quasi-inverse is given by \smash{$T^\lambda_{\lambda + \mu}$}.
\end{Fact}

\subsection{Discrete decomposability}

Our main concern in this paper is discretely decomposable restrictions of $\lie{g}$-modules.
We recall the notion of the discrete decomposability.
We refer the reader to~\cite{Ko98_discrete_decomposable_3,Ko12_generalized_Verma} for the details.

Let $\lie{g}$ be a finite-dimensional complex reductive Lie algebra.

\begin{Definition}
	Let $V$ be a $\lie{g}$-module.
	We say that $V$ is \emph{discretely decomposable} if $V$ has an exhaustive $\lie{g}$-module filtration $0 = V_0\subset V_1 \subset \cdots $
	such that any $V_i$ has finite length.
	Moreover, if all $V_i$ are in a subcategory $\mathcal{C}$ of $\Mod(\lie{g})$, we say that $V$ is \emph{discretely decomposable} in $\mathcal{C}$.
\end{Definition}

By definition, a discretely decomposable $\lie{g}$-module is locally $\univcent{g}$-finite.
Hence the following proposition is a consequence of Fact~\ref{fact:PrimaryDecomposition}.

\begin{Proposition}
	A discretely decomposable $\lie{g}$-module has the primary decomposition.
\end{Proposition}

\subsection{Generalized Verma module}\label{subsect:GeneralizedVerme}

This section contains a brief summary of generalized Verma modules.
We refer the reader to~\cite{Hu08_category_o}.
For the branching problem part, see~\cite{Ko12_generalized_Verma}.

Let $\lie{g}$ be a finite-dimensional complex reductive Lie algebra.
Fix a semisimple element $H \in \lie{g}$ such that $\ad(H)$ has only real eigenvalues.
$\lie{l}(H)$, $\lie{u}(H)$ and $\overline{\lie{u}}(H)$ denote the sums of eigenspaces of $\ad(H)$
with zero, positive and negative eigenvalues, respectively.
Then $\lie{q}(H)\coloneq \lie{l}(H)\oplus \lie{u}(H)$ is a parabolic subalgebra of $\lie{g}$.
Similarly, write $\overline{\lie{q}}(H)\coloneq \lie{l}(H) \oplus \overline{\lie{u}}(H)$ for the opposite parabolic subalgebra.
If $H$ is clear from the context, we omit `$(H)$' part as $\lie{l}$, $\lie{u}$ and $\lie{q}$.

Fix a Cartan subalgebra $\lie{t}$ of $\lie{l}$ and write $W_L$ for the Weyl group of $\lie{l}$.
Then we have $H \in \lie{t}$.
For a finite-dimensional completely reducible $\lie{t}$-module $V$, we denote by $\Delta(V, \lie{t})$ (or $\Delta(V)$) the (multi)set of all non-zero weights in $V$.
We write $\rhog{\lie{u}}$ for half the sum of all roots in $\Delta(\lie{u}, \lie{t})$.

We shall recall generalized Verma modules.
Let $F$ be an $\lie{l}$-module.
Set
\begin{align*}
	\ind{g}{q}(F)\coloneq \univ{g}\otimes_{\univ{q}}F,
\end{align*}
where we consider $F$ as a $\lie{q}$-module letting $\lie{u}$ act on $F$ trivially.
In this paper, we use this extension without mention.
It is well-known that if $F$ has an infinitesimal character $[\lambda] \in \lie{t}^* / W_L$,
then $\ind{g}{q}(F)$ has the infinitesimal character $[\lambda + \rhog{\lie{u}}]$.
If $F$ is a finite-dimensional irreducible $\lie{l}$-module,
$\ind{g}{q}(F)$ is called a \emph{generalized Verma module}.
In this case, $\ind{g}{q}(F)$ is a highest weight module with respect to some Borel subalgebra contained in $\lie{q}$.

The relative BGG category $\BGGcat{}{}$ denoted by $\BGGcat{g}{q}$ is defined as follows (see~\cite[Section~9]{Hu08_category_o}).
$\BGGcat{g}{q}$ is the full subcategory of $\Mod(\lie{g})$ whose object $V$ satisfies the following conditions:
\begin{enumerate}\itemsep=0pt
 \item $V$ is finitely generated as a $\lie{g}$-module.
 \item $V|_{\lie{l}}$ is locally finite and completely reducible.
 \item The action of $\lie{u}$ on $V$ is locally nilpotent.
\end{enumerate}
Then any generalized Verma module $\ind{g}{q}(F)$ is an object of the category $\BGGcat{g}{q}$.
It is well-known that any object in $\BGGcat{g}{q}$ has finite length and hence $\univcent{g}$-finite.

A generalized Verma module may be reducible in general.
The following result provides a~criterion for the irreducibility of a~generalized Verma module (see, e.g.,~\cite[Theorem~9.12]{Hu08_category_o}).
\begin{Fact}\label{fact:VermaIrreducible}
	Let $F$ be an irreducible finite-dimensional $\lie{l}$-module with infinitesimal character $[\lambda] \in \lie{t}^*/W_L$ satisfying
	\begin{align*}
		\frac{2(\lambda+\rhog{\lie{u}}, \alpha)}{(\alpha, \alpha)} \not \in \set{1, 2, \ldots}, \qquad \forall \alpha \in \Delta(\lie{u}, \lie{t}).
	\end{align*}
	Then the generalized Verma module $\ind{g}{q}(F)$ is irreducible.
\end{Fact}

Remark that the assumption of Fact~\ref{fact:VermaIrreducible} implies that $\lambda + \rhog{\lie{u}}$ is integrally anti-dominant for some Borel subalgebra contained in $\lie{q}$.
The image of a generalized Verma module by the translation functor is also a generalized Verma module under a dominance condition (see~\cite[Theorem~7.237]{KnVo95_cohomological_induction}).

\begin{Fact}\label{fact:TranslationVerma}
	Let $F$ be a finite-dimensional $\lie{l}$-module with infinitesimal character $[\lambda] \in \lie{t}^* / W_L$,
	and $\mu \in \lie{t}^*$ an algebraically integral weight.
	Suppose that there exists a Borel subalgebra $\lie{b}$ of $\lie{g}$ containing $\lie{t}$ such that $\lambda+\rhog{\lie{u}}$ and $\lambda+\mu+\rhog{\lie{u}}$ are regular integrally dominant for $\lie{b}$.
	Then there exists a $\lie{g}$-module isomorphism
	\begin{align*}
		T_{\lambda+\rhog{\lie{u}}}^{\lambda+\mu+\rhog{\lie{u}}}\bigl(\ind{g}{q}(F)\bigr)
		\simeq \ind{g}{q}\bigl(T_{\lambda}^{\lambda+\mu}(F)\bigr).
	\end{align*}
\end{Fact}

We recall the notion of standard filtrations in \cite[Section~9.8]{Hu08_category_o}.
\begin{Definition}\label{def:StandardFiltration}
	Let $V$ be a $\lie{g}$-module.
	We say that an exhaustive filtration $0=V_0 \subset V_1 \subset \cdots$ of~$V$ is a \emph{standard filtration} in $\BGGcat{g}{q}$ if each successive quotient $V_{i+1}/{V_{i}}$ is isomorphic to a generalized Verma module
	$\ind{g}{q}(F)$.
\end{Definition}
If the parabolic subalgebra $\lie{q}$ is clear from the context,
we simply say that ${0=V_0 \subset V_1 \subset \cdots}$ is a standard filtration.
We do not assume that a standard filtration is finite.
We give a~fundamental property of the standard filtration.

\begin{Proposition}\label{prop:CompletelyReducibleVerma}
	Let $V$ be a $\lie{g}$-module with a standard filtration $0 = V_0\subset V_1 \subset \cdots$.
	\begin{enumerate}\itemsep=0pt
		\item[$1.$] For any object $W\in \BGGcat{g}{q}$, $V \otimes W$ has an exhaustive filtration such that $\gr(V\otimes W)$ has a~standard filtration.
		\item[$2.$] If each successive quotient $V_{i+1}/V_{i}$ is irreducible, then $V$ is completely reducible and isomorphic to $\gr(V)$.
	\end{enumerate}
\end{Proposition}

\begin{proof}
	The first assertion is clear from the natural isomorphism \smash{$\ind{g}{q}(\cdot)\otimes W \simeq \ind{g}{q}(\cdot \otimes W)$}
	and that $W$ has a $\lie{q}$-module filtration whose successive quotients are irreducible and finite-dimensional.

	To prove the second assertion, we can assume that $V$ has finite length.
	In fact, if each~$V_i$~is completely reducible, then $V$ is completely reducible
	because $V$ is a sum of irreducible submod\-ules.
	Since \smash{$\Ext^1_{\BGGcat{g}{q}}(M, N)=0$} for any two irreducible generalized Verma modules~$M$ and~$N$ (see~\cite[Theorem~3.3\,(d)]{Hu08_category_o}),
	the assertion follows by induction on the length.
\end{proof}

We consider the branching problem of a generalized Verma module.
For simplicity, suppose that $\lie{g}$ is semisimple.
Let $\lie{g}'$ be a reductive subalgebra of $\lie{g}$.
We recall the notion of $\lie{g'}$-compatible parabolic subalgebras in the sense of~\cite{Ko12_generalized_Verma},
and define its generalization.

\begin{Definition}\label{def:WeaklyCompatible}
	Let $\lie{q}$ be a parabolic subalgebra of $\lie{g}$.
	We say that $\lie{q}$ is \emph{$\lie{g'}$-compatible} if there exists $H \in \lie{g'}$ such that $\lie{q} = \lie{q}(H)$.
	We say that $\lie{q}$ is \emph{weakly $\lie{g'}$-compatible} if $\lie{q}$ contains a~$\lie{g'}$-compatible parabolic subalgebra $\lie{q}(H) \subset \lie{g}$ and has a Levi decomposition $\lie{q} = \lie{l} \oplus \lie{u}$ such that
	\begin{alignat}{3}
		&\lie{q}(H) \cap \lie{g'} = \lie{q} \cap \lie{g}', \qquad&& \lie{l}(H)\cap \lie{g'} = \lie{l} \cap \lie{g}',& \nonumber \\
		&\lie{u}(H)\cap \lie{g'} = \lie{u} \cap \lie{g}', \qquad&& \lie{\overline{u}}(H)\cap \lie{g'} = \lie{\overline{u}} \cap \lie{g'},& \label{eqn:WeaklyCompatible}
	\end{alignat}
	where $\lie{\overline{u}}$ is the nilpotent radical of the opposite parabolic subalgebra $\lie{\overline{q}}$ determined by the Levi decomposition $\lie{q} = \lie{l}\oplus \lie{u}$.
\end{Definition}

Formulas \eqref{eqn:WeaklyCompatible} mean that $\lie{q}$, $\lie{l}$, $\lie{u}$ and $\lie{\overline{u}}$ have direct sum decompositions compatible with $\lie{g} = \lie{g'}\oplus (\lie{g'})^\perp$.
Remark that if $(\lie{g}, \lie{g'})$ is a symmetric pair defined by an involution $\theta$,
then $\lie{q}$ is weakly $\lie{g'}$-compatible if and only if $\lie{q}$ is $\theta$-stable,
and they are equivalent to that $\lie{q}$ is $\lie{g'}$-compatible by~\cite[Proposition~4.76]{KnVo95_cohomological_induction}.

\begin{Example}
	Let $\lie{g} = \lieSL(3, \CC)$ and set
	\begin{align*}
		\lie{g'} \coloneq \set{\begin{pmatrix} a & 0 & \hphantom{-}b \\ 0 & 0 & \hphantom{-}0 \\ c & 0 & -a \end{pmatrix} \bigg|\, a, b, c \in \CC}, \qquad H \coloneq \begin{pmatrix} 1 & 0 & \hphantom{-}0 \\ 0 & 0 & \hphantom{-}0 \\ 0 & 0 & -1 \end{pmatrix}.
	\end{align*}
	Then $\lie{q}(H)$ is a Borel subalgebra of $\lie{g}$, and every parabolic subalgebra of $\lie{g}$ containing $\lie{q}(H)$ is weakly $\lie{g'}$-compatible.
	Moreover, any $\lie{g'}$-compatible parabolic subalgebra of $\lie{g}$ is either a Borel subalgebra or $\lie{g}$ itself.
\end{Example}

Let $\lie{q} = \lie{l}\oplus \lie{u}$ be a weakly $\lie{g'}$-compatible parabolic subalgebra of $\lie{g}$
containing $\lie{q}(H)$, $H \in \lie{g'}$.
Set
\begin{align}
	\lie{l'}\coloneq \lie{l}\cap \lie{g'}, \qquad \lie{u'}\coloneq \lie{u}\cap \lie{g'}, \qquad \lie{\overline{u}'}\coloneq \lie{\overline{u}}\cap \lie{g'}, \qquad \lie{q'}\coloneq \lie{l'}\oplus \lie{u'}. \label{eqn:DefinitionQ'}
\end{align}
Then $\lie{q'}$ is a parabolic subalgebra of $\lie{g'}$.

An important fact is that any object of $\BGGcat{g}{q}$ is discretely decomposable
as a $\lie{g}'$-module in $\BGGcat{g'}{q'}$ (see~\cite[Proposition~3.8]{Ko12_generalized_Verma}).
By the same proof as in~\cite[Lemma~6.4.4]{Wa88_real_reductive_I}, we obtain the following.

\begin{Proposition}\label{prop:FiltrationRestriction}
	Let $F$ be a finite-dimensional irreducible $\lie{l}$-module.
	Then $\ind{g}{q}(F)|_{\lie{g'}}$ has a~standard filtration $0 = V_0\subset V_1 \subset \cdots$
	satisfying
	\begin{align*}
		\gr(V) \simeq \ind{g'}{q'}\bigl(F \otimes S\bigl(\lie{\overline{u}}/\lie{\overline{u}'}\bigr)\bigr),
	\end{align*}
	where $\lie{\overline{u}}/\lie{\overline{u}'}$ is regarded as a $\lie{q'}$-module
	by letting $\lie{u'}$ act on $\lie{\overline{u}}/\lie{\overline{u}'}$ trivially.
\end{Proposition}

Since $\lie{g}/\lie{q}$ is a quotient of $\lie{g}/\lie{q}(H)\simeq \lie{\overline{u}}(H)$, any eigenvalue in $\lie{g}/\lie{q}$ of $H$ is negative.
This implies that any eigenspace in $\ind{g}{q}(F)\simeq S(\lie{g}/\lie{q})\otimes F$ of $H$ is finite-dimensional.
From $H \in \lie{l'}$ and this, $\ind{g}{q}(F)$ is $\lie{l}'$-admissible, that is, the multiplicity of any irreducible $\lie{l'}$-module is finite.
The following proposition is well-known if $(\lie{g}, \lie{g'})$ is a symmetric pair.
The proof is the same as the symmetric pair case.
For the notation $(\cdot)_{\lie{l'}}$ in the proposition, see Definition~\ref{def:GKmod}.

\begin{Proposition}\label{prop:AdmissibleDual}
	Let $F$ be a finite-dimensional irreducible $\lie{l}$-module.
	Then $\ind{g}{q}(F)$ is $\lie{l}'$-admissible.
	If, in addition, $\ind{g}{q}(F)$ is irreducible, then $(\ind{g}{q}(F)^*)_{\lie{l'}}$ is isomorphic to
	$\ind{g}{\overline{q}}(F^*)$ as a $\lie{g}$-module.
\end{Proposition}

\begin{proof}
	The first assertion is shown in the above.
	For the second assertion, assume that $\ind{g}{q}(F)$ is irreducible.
	Since $\ind{g}{q}(F)$ is $\lie{l'}$-admissible, we have \smash{$\bigl(\ind{g}{q}(F)^*\bigr)_{\lie{l'}} = \bigl(\ind{g}{q}(F)^*\bigr)_{\lie{l}}$}.
	Since $\ind{g}{q}(F)$ is irreducible, \smash{$\bigl(\ind{g}{q}(F)^*\bigr)_{\lie{l}}$} is an irreducible highest weight module
	and the highest weight is the~same as that of $\ind{g}{\overline{q}}(F^*)$.
	This shows the assertion.
\end{proof}

\subsection{Zuckerman derived functor}

Here we review the Zuckerman derived functor.
We will use the Zuckerman derived functor in two ways.
One is to construct cohomologically induced module, e.g., discrete series representations.
The other is to study $\univ{g}^{G'}$-modules.

Let $(\lie{g}, K)$ be a pair with reductive $\lie{g}$,
and $M$ a connected reductive subgroup of $K$.
Recall that, for a $(\lie{g}, M)$-module $V$, $V_K$ is the sum of all $(\lie{k}, M)$-submodules that lift to $K$-modules.
Then $V_K$ is a $(\lie{g}, K)$-module,
and the functor $\Mod(\lie{g}, M) \ni V \rightarrow V_K \in \Mod(\lie{g}, K)$ is left exact.
For $i \in \NN$, we denote by $\Dzuck{K}{M}{i}$ the right derived functor of the functor $(\cdot)_K$.
The functor $\Dzuck{K}{M}{i}$ is called the $i$-th \emph{Zuckerman derived functor}.

The functor can be constructed explicitly as follows.
We refer the reader to \cite[Section~I.8]{BoWa00_continuous_cohomology}.
Let $V$ be a $(\lie{g}, M)$-module and $i \in \NN$.
Consider the cohomology group $H^i(\lie{k}, M; \rring{K}\otimes V)$.
Here $\rring{K}$ is the coordinate ring and the $(\lie{k}, M)$-cohomology is taken via the tensor product of the right regular action on $\rring{K}$ and the action on $V$.
Then the left regular action on $\rring{K}$ induces a rational $K$-action on $H^i(\lie{k}, M; \rring{K}\otimes V)$.
Moreover, it is known that $H^i(\lie{k}, M; \rring{K}\otimes V)$ admits a $(\lie{g}, K)$-module structure and it is isomorphic to \smash{$\Dzuck{K}{M}{i}(V)$}.
See~\cite[Proposition~I.8.2]{BoWa00_continuous_cohomology}.
The following commutative diagram characterizes the $\univ{g}$-action:
\begin{align*}
	\xymatrix{
		\univ{g}\otimes \Dzuck{K}{M}{i}(V) \ar[d]^{\simeq} \ar[rr]^-m && \Dzuck{K}{M}{i}(V) \ar@{=}[d] \\
		\Dzuck{K}{M}{i}(\univ{g}\otimes V) \ar[rr]^-{\Dzuck{K}{M}{i}(m)} && \Dzuck{K}{M}{i}(V),
	}
\end{align*}
where $m$'s are the multiplication maps and $\univ{g}$ is regarded as a $(\lie{g}, K)$-module via the adjoint action.
See~\cite[Lemma~6.3.1]{Wa88_real_reductive_I}.
From this and Proposition~\ref{prop:UKaction}, we have the following proposition.

\begin{Proposition}\label{prop:KInvZuckerman}
	Let $V$ be a $(\lie{g}, M)$-module.
	Then there exists a natural $\univ{g}^K$-module isomorphism $\Dzuck{K}{M}{i}(V)^K \simeq H^i(\lie{k}, M; V)$.
	In particular, the length of $H^i(\lie{k}, M; V)$ is bounded by that of $\Dzuck{K}{M}{i}(V)$ from above.
\end{Proposition}

\begin{Remark}
	The $\univ{g}^K$-module structure on $H^i(\lie{k}, M; V)$ is induced from that on $V$.
	Indeed, the multiplication $X\cdot \colon V\rightarrow V$ by $X \in \univ{g}^K$ is a $(\lie{k}, M)$-homomorphism, and $H^i(\lie{k}, M;\cdot)$ is a~functor on the category of $(\lie{k}, M)$-modules.
\end{Remark}

We summarize fundamental results about the Zuckerman derived functor.
We refer the reader to~\cite[Theorems 2.103 and 5.21, and the proof of Theorem~7.237]{KnVo95_cohomological_induction}.

\begin{Fact}\label{fact:BasicZuckerman}
	Let $V$ be a $(\lie{g}, M)$-module.
	\begin{enumerate}\itemsep=0pt
		\item[$1.$] If $V$ has an infinitesimal character $\chi$, then $\Dzuck{K}{M}{i}(V)$ has the infinitesimal character $\chi$.
		\item[$2.$] Let $W$ be a $(\lie{g}, K)$-module. Then there exists a natural isomorphism \[\Dzuck{K}{M}{i}(V)\otimes W \simeq \Dzuck{K}{M}{i}(V\otimes W).\]
		\item[$3.$] $\Dzuck{K}{M}{i}$ commutes with the translation functors $T_{\lambda}^{\lambda + \mu}$ if $\mu$ is an extreme weight of a finite-dimensional irreducible $(\lie{g}, K)$-module.
	\end{enumerate}
\end{Fact}

\begin{Lemma}\label{lem:LimitZuckerman}
	Let $V$ be a $(\lie{g}, M)$-module with an exhaustive filtration
	$0 = V_0 \subset V_1 \subset \cdots $.
	If~\smash{$\Dzuck{K}{M}{i}(V_j) = 0$} for any $j$, then \smash{$\Dzuck{K}{M}{i}(V)=0$} holds.
\end{Lemma}
	
\begin{proof}
	Recall that the cohomology $H^i(\lie{k}, M; V \otimes \rring{K})\bigl(\simeq \Dzuck{K}{M}{i}(V)\bigr)$ is computed by the standard complex
	\begin{align*}
		C^i\coloneq \Hom_{M}\bigl(\wedge^i(\lie{k}/\lie{m}), V \otimes \rring{K}\bigr)
	\end{align*}
	with $d\colon C^i \rightarrow C^{i+1}$.
	Let $\omega \in C^i$ such that $\mathrm{d}\omega = 0$.
	Then $\Im(\omega)$ is contained in $V_j \otimes \rring{K}$ for some $j$.
	Hence, by assumption, there exists $\omega' \in C^{i-1}$ such that $\mathrm{d}\omega'=\omega$.
	This shows the assertion.
\end{proof}

We shall state a part of an algebraic analogue of the Borel--Weil--Bott theorem~\cite[Corollary~4.160]{KnVo95_cohomological_induction}.
Let $(\lie{g}, G)$ be a pair.
Note that $\lie{g}$ is reductive since $G$ is reductive.
Let $\lie{q} = \lie{l} \oplus \lie{u}$ be a~parabolic subalgebra of $\lie{g}$
and fix a Cartan subalgebra $\lie{t}\subset \lie{l}$.
Write $L$ for the analytic subgroup in $G$ with the Lie algebra $\lie{l}$.
Then $L$ is a Levi subgroup of $G$.

Fix a set $\Delta^+(\lie{l}, \lie{t})$ of positive roots in $\Delta(\lie{l}, \lie{t})$
and put $\Delta^+(\lie{g}, \lie{t}) \coloneq \Delta^+(\lie{l}, \lie{t})\cup \Delta(\lie{u}, \lie{t})$.
Write $\rhog{\lie{u}}$ for half the sum of roots in $\Delta(\lie{u}, \lie{t})$.
Set $S \coloneq \dim(\lie{u})$.

\begin{Fact}\label{fact:BorelWeilBott}
	Let $F$ be an irreducible $(\lie{l}, L)$-module with infinitesimal character $[\lambda]$.
	\begin{enumerate}\itemsep=0pt
		\item[$1.$] $\Dzuck{G}{L}{S}\bigl(\ind{g}{q}(F)\bigr)$ is zero or irreducible.
		\item[$2.$] $\Dzuck{G}{L}{S}\bigl(\ind{g}{q}(F)\bigr) \neq 0$ if and only if
		\begin{align*}
			(\lambda + \rhog{\lie{u}}, \alpha) < 0, \qquad \forall \alpha \in \Delta(\lie{u}, \lie{t}).
		\end{align*}
		\item[$3.$] If an irreducible $(\lie{l}, L)$-module $F'$ satisfies
		\begin{align*}
			\Dzuck{G}{L}{S}\bigl(\ind{g}{q}(F)\bigr)\simeq \Dzuck{G}{L}{S}\bigl(\ind{g}{q}(F')\bigr) \neq 0,
		\end{align*}
		then $F\simeq F'$ holds.
	\end{enumerate}
\end{Fact}

\begin{Remark}
	The fact is the case $q = S$ of~\cite[Corollary~4.160]{KnVo95_cohomological_induction}, that is,
	\begin{align*}
		|\Delta(\lie{u}, \lie{t})| = S = \bigl|\bigl\{\alpha \in \Delta^+(\lie{u}, \lie{t}) \mid (\lambda + \rhog{\lie{u}}, \alpha) < 0\bigr\}\bigr|.
	\end{align*}
	This condition corresponds to that of 2 in the fact.
	\smash{$\Dzuck{G}{L}{S}\bigl(\ind{g}{q}(F)\bigr)$} in the fact corresponds to~\smash{$\Pi_{S}\bigl(P^{\lie{g}, L}_{\lie{q}, L}\bigl(\mathcal{F}^{\lie{q}, L}_{\lie{l}, L}(Z)\bigr)\bigr)$}.
	In our case, \smash{$P^{\lie{g}, L}_{\lie{q}, L}\bigl(\mathcal{F}^{\lie{q}, L}_{\lie{l}, L}(Z)\bigr)$} is $\ind{g}{q}(F)$,
	and the functor $\Pi_{S}$ is isomorphic to $\Dzuck{G}{L}{S}$ by the Zuckerman duality~\cite[Corollary~3.7]{KnVo95_cohomological_induction}.
\end{Remark}

\section{Cohomological induction}

In this section, we consider cohomologically induced $(\lie{g}, K)$-modules.
We deal with a vanishing theorem and irreducibility.
The results are well-known if $(\lie{g}, \lie{k})$ is a symmetric pair.
We give proofs for the results in the general setting to make this paper self-contained
even though the proofs are essentially the same as the symmetric case, e.g.,~\cite[Section~6]{Wa88_real_reductive_I}.

\subsection{Vanishing theorem}

Let $(\lie{g}, K)$ be a pair with semisimple $\lie{g}$.
Let $\lie{q} = \lie{l}\oplus \lie{u}$ be a weakly $\lie{k}$-compatible parabolic subalgebra of $\lie{g}$
containing $\lie{q}(H)$, $H \in \lie{g'}$.
See Definition~\ref{def:WeaklyCompatible}.
Set
\begin{align*}
	\lie{l}_K\coloneq \lie{l}\cap \lie{k}, \qquad \lie{u}_K\coloneq \lie{u}\cap \lie{k}, \qquad \lie{\overline{u}}_K\coloneq \lie{\overline{u}}\cap \lie{k}, \qquad \lie{q}_K\coloneq \lie{l}_K\oplus \lie{u}_K.
\end{align*}
Then $\lie{q}_K$ is a parabolic subalgebra of $\lie{k}$.
Let $\lie{\overline{q}}$ and $\lie{\overline{q}}_K$ denote the opposite parabolic subalgebras of $\lie{q}$ and $\lie{q}_K$, respectively.
Let $L_K$ denote the analytic subgroup in $K$ with the Lie algebra $\lie{l}_K$, which is a Levi subgroup of $K$.
Put $S\coloneq \dim(\lie{u}_K)$.

We shall consider cohomologically induced modules \smash{$\Dzuck{K}{L_K}{i}\bigl(\ind{g}{q}(F)\bigr)$}.
When $(\lie{g}, \lie{k})$ is a symmetric pair, such modules are classical and well studied.
See, e.g.,~\cite{KnVo95_cohomological_induction} and~\cite[Section~6]{Wa88_real_reductive_I}.
In the context of generalized Harish-Chandra modules, $\Dzuck{K}{L_K}{i}\bigl(\ind{g}{q}(F)\bigr)$ is studied in~\cite{PeZu04_generalized_harish_chandra,PeZu12_generalized_harish_chandra}.
The following fact is fundamental.
See~\cite[Lemma~6.4.2]{Wa88_real_reductive_I}.

\begin{Fact}\label{fact:StandardFiltrationVanishing}
	Let $M$ be a $(\lie{k}, L_K)$-module with a standard filtration $($see Definition $\ref{def:StandardFiltration})$.
	Then $\Dzuck{K}{L_K}{i}(M)$ vanishes for any $i < S$.
\end{Fact}

The following proposition is used in Lemma~\ref{lem:ExactnessPart}.

\begin{Proposition}\label{prop:StandardFiltrationVanishing}
	Let $M$ be a $(\lie{k}, L_K)$-module with an exhaustive filtration such that $\gr(M)$ has a standard filtration.
	Then $\Dzuck{K}{L_K}{i}(M)=0$ holds for any $i < S$.
\end{Proposition}

\begin{proof}
	The assertion follows from Lemma~\ref{lem:LimitZuckerman} and Fact~\ref{fact:StandardFiltrationVanishing}.
\end{proof}

We recall a vanishing theorem and existence of a cyclic subspace.
Our proofs are essentially the same as in~\cite[Section~6]{Wa88_real_reductive_I}.

\begin{Proposition}\label{prop:Vanishing}
	Let $F$ be an irreducible finite-dimensional $(\lie{l}, L_K)$-module.
	Then we have
	\begin{align*}
		\Dzuck{K}{L_K}{i}\bigl(\ind{g}{q}(F)\bigr) = 0, \qquad \forall i < S.
	\end{align*}
	Moreover, if $\ind{g}{q}(F)$ is irreducible, then
	\begin{align*}
		\Dzuck{K}{L_K}{i}\bigl(\ind{g}{q}(F)\bigr) = 0, \qquad \forall i \neq S.
	\end{align*}
\end{Proposition}
	
\begin{proof}
	The first assertion is a direct consequence of Proposition~\ref{prop:FiltrationRestriction} and
	Fact~\ref{fact:StandardFiltrationVanishing}.
	
	Assume that $\ind{g}{q}(F)$ is irreducible and fix $i \in \NN$.
	By Proposition~\ref{prop:AdmissibleDual}, we have
 \[
 \bigl(\ind{g}{q}(F)^*\bigr)_{L_K} \simeq \ind{g}{\overline{q}}(F^*).
 \]
	By the duality~\cite[Corollary~3.7]{KnVo95_cohomological_induction}, there exists an isomorphism
 \[
 \bigl(\Dzuck{K}{L_K}{i}\bigl(\ind{g}{q}(F)\bigr)^*\bigr)_K\simeq \Dzuck{K}{L_K}{2S-d}\bigl(\ind{g}{\overline{q}}(F^*)\bigr).
 \]
	Using the first assertion for \smash{$\Dzuck{K}{L_K}{2S-d}\bigl(\ind{g}{\overline{q}}(F^*)\bigr)$},
	we have
	\begin{align*}
		(\Dzuck{K}{L_K}{i}\bigl(\ind{g}{q}(F)\bigr)^*)_K\simeq \Dzuck{K}{L_K}{2S-i}\bigl(\ind{g}{\overline{q}}(F^*)\bigr)=0\qquad \text{if}\ 2S-i < S.
	\end{align*}
	This shows the second assertion.
\end{proof}

\subsection{Cyclic subspace}

We shall construct a cyclic subspace in the Zuckerman derived functor module
under some dominance condition.
Retain the notation in the previous subsection.

Let $\lie{t}_K$ be a Cartan subalgebra of $\lie{l}_K$
and write $\rhog{\lie{u}_K}$ for half the sum of all roots in $\Delta(\lie{u}_K, \lie{t}_K)$.
Let $W_K$ (resp.\ $W_{L_K}$) denote the Weyl group of $\lie{k}$ (resp.\ $\lie{l}_K$).
For an $(\lie{l}, L_K)$-module $F$, we have a natural homomorphism $\ind{k}{q_\mathit{K}}(F) \rightarrow \ind{g}{q}(F)$ of $(\lie{k}, L_K)$-modules, and it induces a $K$-module homomorphism
\begin{align*}
	B_F\colon\ \Dzuck{K}{L_K}{S}\bigl(\ind{k}{q_\mathit{K}}(F)\bigr) \rightarrow \Dzuck{K}{L_K}{S}\bigl(\ind{g}{q}(F)\bigr).
\end{align*}
This map is called a \emph{bottom-layer map} in~\cite[p.~365]{KnVo95_cohomological_induction} (when $(\lie{g}, \lie{k})$ is symmetric).

\begin{Lemma}\label{lem:MinimalTypeOne}
	Let $F$ be a finite-dimensional irreducible $(\lie{l}, L_K)$-module.
	Then for each irreducible $K$-submodule $V$ of $\Dzuck{K}{L_K}{S}\bigl(\ind{k}{q_\mathit{K}}(F)\bigr)$,
	the linear map
	\begin{align*}
		(B_F)_*\colon\ \Hom_{K}\bigl(V, \Dzuck{K}{L_K}{S}\bigl(\ind{k}{q_\mathit{K}}(F)\bigr)\bigr) \rightarrow \Hom_{K}\bigl(V, \Dzuck{K}{L_K}{S}\bigl(\ind{g}{q}(F)\bigr)\bigr)
	\end{align*}
	is bijective.
	In particular, $B_F$ is injective.
\end{Lemma}

\begin{proof}
	Let $V$ be an irreducible $K$-submodule of $\Dzuck{K}{L_K}{S}\bigl(\ind{k}{q_\mathit{K}}(F)\bigr)$.
	Then there exists a unique irreducible $L_K$-submodule $F_0$ of $F$ such that
	the image of the natural homomorphism
	\begin{align*}
		\Dzuck{K}{L_K}{S}\bigl(\ind{k}{q_\mathit{K}}(F_0)\bigr) \rightarrow \Dzuck{K}{L_K}{S}\bigl(\ind{k}{q_\mathit{K}}(F)\bigr)
	\end{align*}
	is equal to $V$.
	Write $[\lambda] \in \lie{t}^*_K / W_{L_K}$ for the infinitesimal character of $F_0$.
	Then $(\lambda + \rhog{\lie{u}_K}, \beta) < 0$ holds for any $\beta \in \Delta(\lie{u}_K, \lie{t}_K)$,
	and $\Dzuck{K}{L_K}{S}\bigl(\ind{k}{q_\mathit{K}}(F_0)\bigr)$ has the infinitesimal character $[\lambda+\rhog{\lie{u}_K}]$.
	See Fact~\ref{fact:BorelWeilBott}.

	Set $W\coloneq\univ{k}(1\otimes F)$ in $\ind{g}{q}(F)$, which
	is a $\lie{k}$-submodule of $\ind{g}{q}(F)$ isomorphic to $\ind{k}{q_\mathit{K}}(F)$.
	We shall show that $W$ contains the primary component \smash{$P_{[\lambda + \rhog{\lie{u}_K}]}\bigl(\ind{g}{q}(F)|_K\bigr)$}.
	If so, we have
	\begin{align*}
		P_{[\lambda + \rhog{\lie{u}_K}]}\bigl(\Dzuck{K}{L_K}{S}\bigl(\ind{g}{q}(F)\bigr)|_K\bigr) \subset \Im(B_F).
	\end{align*}
	by Fact~\ref{fact:BasicZuckerman}, and this shows the assertion.

	By Proposition~\ref{prop:FiltrationRestriction}, there exists a standard filtration
	$0 = M_0 \subset M_1 \subset \cdots$ of $M \coloneq \ind{g}{q}(F)$ such that
	\smash{$\gr(M) \simeq \ind{k}{q_\mathit{K}}(F\otimes S(\lie{\bar{u}}/\lie{\bar{u}}_K))$},
	where we let $\lie{u}_K$ act on $S(\lie{\bar{u}}/\lie{\bar{u}}_K)$ trivially.
	Since \smash{$P_{[\lambda + \rhog{\lie{u}_K}]}$} is an exact functor,
	it is enough to show
	\begin{align*}
		P_{[\lambda + \rhog{\lie{u}_K}]}\bigl(\ind{k}{q_\mathit{K}}(F\otimes S(\lie{\bar{u}}/\lie{\bar{u}}_K))\bigr)
		\subset \ind{k}{q_\mathit{K}}(F\otimes \CC).
	\end{align*}

	Let $F'$ be an irreducible $L_K$-submodule of $F\otimes S(\lie{\bar{u}}/\lie{\bar{u}}_K)$ with infinitesimal character $[\mu] \in \lie{t}_K^*/W_{L_K}$.
	Then $\mu(H) \leq \lambda(H)$, and the equality holds only if $F'\subset F\otimes \CC$ since all eigenvalues of $\ad(H)$ in $\lie{\overline{u}}$ are negative.
	Assume $[\mu + \rhog{\lie{u}_K}] = [\lambda+\rhog{\lie{u}_K}]$ in $\lie{t}_K^*/W_K$ and let us show $\mu(H) = \lambda(H)$.
	Note that $\ind{k}{q_\mathit{K}}(F')$ has the infinitesimal character $[\mu + \rhog{\lie{u}_K}] \in \lie{t}_K^*/W_K$.
	By assumption, there exists $s \in W_K$ such that $\mu + \rhog{\lie{u}_K} = s(\lambda+\rhog{\lie{u}_K})$.

	Since $\lambda+\rhog{\lie{u}_K}$ is an algebraically integral weight and dominant with respect to $-\Delta(\lie{u}_K, \lie{t}_K)$,
	there is a sum $R$ of elements in $\Delta(\lie{l}_K, \lie{t}_K)\cup \Delta(\lie{u}_K, \lie{t}_K)$ such that $s(\lambda+\rhog{\lie{u}_K})=\lambda+\rhog{\lie{u}_K} + R$.
	Hence we have
	\begin{align*}
		(\mu + \rhog{\lie{u}_K})(H) = s(\lambda+\rhog{\lie{u}_K})(H) \geq (\lambda+\rhog{\lie{u}_K})(H)
	\end{align*}
	Recall that we have $R(H)\geq 0$ by $\lie{u}_K\subset \lie{u}(H)$ (see Definition~\ref{def:WeaklyCompatible}).
	Combining this with $\mu(H) \leq \lambda(H)$, we obtain $\mu(H)=\lambda(H)$, and hence $F' \subset F\otimes \CC$.
\end{proof}

The following proposition is proved by the same way as the proof of Lemma~\ref{lem:MinimalTypeOne}.

\begin{Proposition}\label{prop:Ktype}
	Let $F$ be a finite-dimensional irreducible $(\lie{l}, L_K)$-module.
	Let $V_0$ and $V_1$ be irreducible $K$-submodules of $\Im(B_F)$ and $\Dzuck{K}{L_K}{S}\bigl(\ind{g}{q}(F)\bigr)$, respectively.
	For each $i = 0, 1$, take a representative $\lambda_i \in \lie{t}^*_K$ of the infinitesimal character of $V_i$ such that
	\begin{align*}
		(\lambda_i, \beta) < 0, \qquad \forall \beta \in \Delta(\lie{u}_K, \lie{t}_K).
	\end{align*}
	Then one has $\lambda_0(H) \geq \lambda_1(H)$ and the equality holds only if $V_1 \subset \Im(B_F)$.
\end{Proposition}

We prove that the subspace $\Im(B_F)$ in Lemma~\ref{lem:MinimalTypeOne} generates the $\lie{g}$-module $\Dzuck{K}{L_K}{S}\bigl(\ind{g}{q}(F)\bigr)$ under some dominance condition.
To do so, we need the following fact (see, e.g., \cite[Lem\-ma~6.A.1.3]{Wa88_real_reductive_I}).

\newcommand{\del}{\partial}

\begin{Fact}\label{fact:KoszulResolution}
	Let $\lie{g}$ be a complex Lie algebra and
	$\lie{k}$ a subalgebra of $\lie{g}$.
	Let $V$ be a $\lie{g}$-module.
	Set $n\coloneq\dim(\lie{g}/\lie{k})$.
	Then there exists an exact sequence
	\begin{align*}
0&\xrightarrow{\del_n} \univ{g}\otimes_{\univ{k}}(\wedge^n (\lie{g}/\lie{k})\otimes V)
		\xrightarrow{\del_{n-1}} \univ{g}\otimes_{\univ{k}}\bigl(\wedge^{n-1} (\lie{g}/\lie{k})\otimes V\bigr)\\
		& \xrightarrow{\del_{n-2}} \cdots \xrightarrow{\del_{1}} \univ{g}\otimes_{\univ{k}}(\lie{g}/\lie{k}\otimes V)
		\xrightarrow{\del_{0}} \univ{g}\otimes_{\univ{k}}V
		\xrightarrow{\epsilon} V \rightarrow 0
	\end{align*}
	of $\lie{g}$-modules.
	The last homomorphism $\epsilon$ is the multiplication map.
\end{Fact}

\begin{Lemma}\label{lem:CyclicSubspace}
	Let $F$ be a finite-dimensional irreducible $(\lie{l}, L_K)$-module.
	Suppose that
 \[
 \biggl(\lambda +\rhog{\lie{u}_K} + \sum_{\alpha \in E}\alpha, \beta\biggr) < 0
 \]
	for any $\univcent{l_\mathit{K}}$-character $[\lambda]$ in $F|_{\univcent{l_{\mathit{K}}}}$, $\beta \in \Delta(\lie{u}_K, \lie{t}_K)$ and $E \subset \Delta(\lie{q}/\lie{q}_K, \lie{t}_K)$.
	Then $\Im(B_F)$ in Lemma~$\ref{lem:MinimalTypeOne}$ generates \smash{$\Dzuck{K}{L_K}{S}\bigl(\ind{g}{q}(F)\bigr)$} as a $\lie{g}$-module.
\end{Lemma}

\begin{Remark}
	In this lemma, we regard $\Delta(\lie{q}/\lie{q}_K, \lie{t}_K)$ and $E$ as multisets.
\end{Remark}

\begin{proof}[Proof of Lemma~\ref{lem:CyclicSubspace}]
	By Fact~\ref{fact:KoszulResolution} for $\lie{q}_K\subset \lie{q}$ and $F$, we have an exact sequence
	\begin{align*}
0&\xrightarrow{\del_n} \univ{q}\otimes_{\univ{q_\mathit{K}}}(\wedge^n (\lie{q/q_\mathit{K}})\otimes F)
		\xrightarrow{\del_{n-1}} \univ{q}\otimes_{\univ{q_\mathit{K}}}(\wedge^{n-1} (\lie{q/q_\mathit{K}})\otimes F)\\
		& \xrightarrow{\del_{n-2}} \cdots \xrightarrow{\del_{1}} \univ{q}\otimes_{\univ{q_\mathit{K}}}(\lie{q/q_\mathit{K}}\otimes F)
		\xrightarrow{\del_{0}} \univ{q}\otimes_{\univ{q_\mathit{K}}}F
		\xrightarrow{\epsilon} F \rightarrow 0
	\end{align*}
	of $\lie{q}$-modules.
	Recall that $\univ{g}$ is free as a right $\univ{q}$-module.
	Applying the exact functor $\univ{g}\otimes_{\univ{q}}(\cdot)$
	to the above exact sequence, we obtain an exact sequence
	\begin{align}
0&\xrightarrow{\del_n} \univ{g}\otimes_{\univ{k}}E_n
		\xrightarrow{\del_{n-1}} \univ{g}\otimes_{\univ{k}}E_{n-1} \nonumber\\
		& \xrightarrow{\del_{n-2}} \cdots \xrightarrow{\del_{1}} \univ{g}\otimes_{\univ{k}}E_1
		\xrightarrow{\del_{0}} \univ{g}\otimes_{\univ{k}}E_0
		\xrightarrow{\epsilon} \ind{g}{q}(F) \rightarrow 0 \label{eqn:CyclicSpaceExactSequence}
	\end{align}
	of $\lie{g}$-modules, where \smash{$E_i\coloneq\ind{k}{q_\mathit{K}}\bigl(\wedge^i (\lie{q/q_\mathit{K}})\otimes F\bigr)$}.
	Set $W \coloneq \univ{k}(1\otimes F) \subset \ind{g}{q}(F)$ as in the proof of Lemma~\ref{lem:MinimalTypeOne}.
	Then we have $\epsilon(1\otimes E_0)=W$.

	We shall show that $\epsilon_S\coloneq \Dzuck{K}{L_K}{S}(\epsilon)$ is surjective.
	By the dominance assumption, Fact~\ref{fact:VermaIrreducible} and Proposition~\ref{prop:CompletelyReducibleVerma}, each $E_i$ is completely reducible.
	Since $\univ{g}\otimes_{\univ{k}}E_i$
	is isomorphic to $S(\lie{g}/\lie{k})\otimes E_i$ as a $\lie{k}$-module, Proposition~\ref{prop:Vanishing} shows
	\begin{align*}
		\Dzuck{K}{L_K}{d}\bigl(\univ{g}\otimes_{\univ{k}}E_i\bigr)
		\simeq S(\lie{g/k})\otimes \Dzuck{K}{L_K}{d}(E_i) = 0
	\end{align*}
	for any $d \neq S$.
	By \eqref{eqn:CyclicSpaceExactSequence} and the long exact sequence, we have exact sequences
	\begin{align*}
		&0 \rightarrow \Dzuck{K}{L_K}{S}(\Im(\del_0))\rightarrow \Dzuck{K}{L_K}{S}\bigl(\univ{g}\otimes_{\univ{k}}E_0\bigr)
		\xrightarrow{\epsilon_S} \Dzuck{K}{L_K}{S}\bigl(\ind{g}{q}(F)\bigr)\\
		& \hphantom{0} \rightarrow \Dzuck{K}{L_K}{S+1}(\Im(\del_0)) \rightarrow 0, \\
		&0 \rightarrow \Dzuck{K}{L_K}{d}(\Im(\del_{i})) \rightarrow \Dzuck{K}{L_K}{d+1}(\Im(\del_{i+1})) \rightarrow 0
	\end{align*}
	for any $d > S$ and $i \in \NN$.
	The second exact sequence implies $\Dzuck{K}{L_K}{d}(\Im(\del_{i}))=0$ for any $d > S$ and $i \in \NN$.
	Hence $\epsilon_S$ is surjective by the first exact sequence.

	By~\cite[Lemma~6.3.1]{Wa88_real_reductive_I}, the following diagram is commutative:
	\begin{align*}
		\xymatrix{
		\univ{g}\otimes \Dzuck{K}{L_K}{S}(E_0) \simeq \univ{g}\otimes W \ar[d]^{\simeq} \ar[drr]^{m} \\
		\Dzuck{K}{L_K}{S}(\univ{g}\otimes E_0) \ar[rr]^{\Dzuck{K}{L_K}{S}(m)} && \Dzuck{K}{L_K}{S}\bigl(\ind{g}{q}(F)\bigr),
		}
	\end{align*}
	where $m$'s are the multiplication maps and the vertical isomorphism is given by Fact~\ref{fact:BasicZuckerman}.
	Note that we regard $\univ{g}$ as a $(\lie{g}, K)$-module via the adjoint action.
	The map $\Dzuck{K}{L_K}{S}(m)$ is written as the composition
	\begin{align*}
		\Dzuck{K}{L_K}{S}(\univ{g}\otimes E_0) \rightarrow \Dzuck{K}{L_K}{S}\bigl(\univ{g}\otimes_{\univ{k}} E_0\bigr) \xrightarrow{\epsilon_S} \Dzuck{K}{L_K}{S}\bigl(\ind{g}{q}(F)\bigr).
	\end{align*}
	Since the canonical surjection $\univ{g}\otimes E_0 \rightarrow \univ{g}\otimes_{\univ{k}} E_0$ of $(\lie{k}, L_K)$-modules has a right inverse and $\epsilon_S$ is surjective, \smash{$\Dzuck{K}{L_K}{S}(m)$} is surjective.
	By the above commutative diagram, we obtain \smash{$\univ{g}W=\Dzuck{K}{L_K}{S}\bigl(\ind{g}{q}(F)\bigr)$}.
\end{proof}

\begin{Corollary}\label{cor:Irreducible}
	Retain the setting in Lemma~$\ref{lem:CyclicSubspace}$.
	Let $N$ denote the unique maximal submodule of $\ind{g}{q}(F)$ and set $M\coloneq \ind{g}{q}(F)/N$.
	Assume that $F|_{L_K}$ is irreducible and \smash{$\Dzuck{K}{L_K}{S+1}(N) = 0$} holds.
	Then \smash{$\Dzuck{K}{L_K}{S}(M)$} is irreducible.
\end{Corollary}

\begin{proof}
	Write $q\colon \ind{g}{q}(F)\rightarrow M$ for the quotient map.
	By the proof of Lemma~\ref{lem:MinimalTypeOne} and the dominance assumption, $\univ{k}(1\otimes F)$ is a direct sum of some primary components in $\ind{g}{q}(F)|_{\lie{k}}$, and hence $q(\univ{k}(1\otimes F))$ is a direct summand of $M|_{\lie{k}}$.
	Since $1\otimes F$ contains a highest weight vector in $\ind{g}{q}(F)$, $q(\univ{k}(1\otimes F))$ is non-zero.
	By Fact~\ref{fact:BorelWeilBott} and the dominance assumption in Lemma~\ref{lem:CyclicSubspace}, \smash{$\Dzuck{K}{L_K}{S}(q(\univ{k}(1\otimes F)))$} is non-zero.
	Hence we can identify \smash{$\Dzuck{K}{L_K}{S}(q(\univ{k}(1\otimes F)))$} as a non-zero subspace $V_0$ of \smash{$\Dzuck{K}{L_K}{S}(M)$}.

	The assumption \smash{$\Dzuck{K}{L_K}{S+1}(N) = 0$} implies that \smash{$\Dzuck{K}{L_K}{S}(M)$} is a quotient of \smash{$\Dzuck{K}{L_K}{S}\bigl(\ind{g}{q}(F)\bigr)$}.
As a~consequence, by Lemma~\ref{lem:CyclicSubspace}, \smash{$\Dzuck{K}{L_K}{S}(M)$} is generated by the non-zero $K$-submodule $V_0$.
	By the assumption that $F|_{L_K}$ is irreducible, $V_0$ is irreducible.
	See Fact~\ref{fact:BorelWeilBott}.

	By the duality~\cite[Corollary~3.7]{KnVo95_cohomological_induction}, there exists a $(\lie{g}, K)$-module isomorphism
 \[
 \bigl(\Dzuck{K}{L_K}{S}(M)^*\bigr)_K\simeq \Dzuck{K}{L_K}{S}\bigl((M^*)_{L_K}\bigr).
 \]
	We identify them via the isomorphism.
	The highest weight module $(M^*)_{L_K} = (M^*)_{\lie{l}}$ is isomorphic to the unique irreducible quotient of \smash{$\ind{g}{\overline{q}}(F^*)$}.
	By the same argument as for \smash{$\Dzuck{K}{L_K}{S}(M)$}, the dual \smash{$\Dzuck{K}{L_K}{S}((M^*)_{L_K})$} is generated by a unique irreducible $K$-submodule $V'_0$ isomorphic to~$V_0^*$.
	If \smash{$\Dzuck{K}{L_K}{S}(M)$} has a non-trivial quotient, then~\smash{$\Dzuck{K}{L_K}{S}(M^*_{L_K})$} has a non-trivial \mbox{submodule} \mbox{containing} $V'_0$, and this is contradiction.
	Therefore, \smash{$\Dzuck{K}{L_K}{S}(M)$} is irreducible.
\end{proof}

\section{Embedding of categories}\label{sect:Embedding}

We have seen in Propositions~\ref{prop:UKaction} and~\ref{prop:gg'G'} that irreducible $\univ{g}^{G'}$-modules are constructed from irreducible $(\lie{g'\oplus g}, \Delta(G'))$-modules.
In this section, we give an embedding from the relative BGG category $\BGGcat{g}{q}$ to $\Mod(\lie{g'\oplus g}, \Delta(G'))$.
As a consequence, we obtain a family of finite-dimensional irreducible $\univ{g}^{G'}$-modules.

\subsection{Setting and main theorem}\label{subsect:EmbeddingSetting}

The main concern in this section is to construct a category embedding
from the BGG category to the category of generalized Harish-Chandra modules.
We shall give a setting and state the main result.

Let $(\lie{g}, G')$ be a pair with semisimple $\lie{g}$.
Then $(\lie{g'\oplus g}, \Delta(G'))$ is a pair, where $\Delta\colon G'\rightarrow G'\times G'$ is the diagonal embedding.
Let $\lie{q} = \lie{l}\oplus \lie{u}$ be a weakly $\lie{g'}$-compatible parabolic subalgebra of $\lie{g}$
containing $\lie{q}(H)$, $H \in \lie{g'}$ (see Definition~\ref{def:WeaklyCompatible}).
Define $\lie{q'} = \lie{l'} \oplus \lie{u'} = \lie{q}\cap \lie{g'}$ as \eqref{eqn:DefinitionQ'}.
Then their Levi decompositions determine opposite parabolic subalgebras $\lie{\overline{q}} = \lie{l}\oplus \lie{\overline{u}}$ and $\lie{\overline{q}}' = \lie{l}\oplus \lie{\overline{u}}'$ for $\lie{q}$ and $\lie{q'}$, respectively.
Note that $\lie{q'\oplus q}$ and \smash{$\lie{\overline{q'} \oplus \overline{q}}$} are weakly $\Delta(\lie{g'})$-compatible.
Write $L'$ for the centralizer of $H$ in $G'$.

Fix a Cartan subalgebra $\lie{t}'$ of $\lie{l}'$ and write $W_{G'}$ and $W_{L'}$ for the Weyl groups.
Let $\rhog{\lie{u}'}$ denote half the sum of roots in $\Delta(\lie{u}', \lie{t}')$.
Similarly, we define $\lie{t}$, $W_{G}$, $W_{L}$ and $\rhog{\lie{u}}$ for $\lie{g}$ to satisfy $\lie{t'}\subset \lie{t}$.

We have the three Lie algebras $\lie{g'}\oplus 0$, $0\oplus \lie{g'}$ and $\Delta(\lie{g'})$ in $\lie{g'}\oplus \lie{g}$ isomorphic to $\lie{g'}$.
If there is no need to distinguish them, we omit $\Delta$.
We use similar notation for subalgebras of $\lie{g'}$.

For a finite-dimensional $\lie{l}'$-module $F$,
we denote by $\BGGcat{g'}{q'}(F)$ the full subcategory of \smash{$\BGGcat{g'}{q'}$} whose object $M$ satisfies that $M\otimes F$ lifts to an $L'$-module.
Then \smash{$\BGGcat{g'}{q'}(F)$} is closed under taking finite direct sums and subquotients.

Our main purpose in this section is to show the following theorem.
We denote by $\lattice(M)$ the \emph{lattice of submodules} of a module $M$.
Set $S\coloneq \dim(\lie{u'})$.
For a finite-dimensional irreducible $\lie{l}$-module $F$, we denote by $T_F$ (or $T$ for short) the functor \smash{$\Dzuck{\Delta(G')}{\Delta(L')}{S}\bigl(\cdot \otimes \ind{g}{q}(F)\bigr)$} from \smash{$\BGGcat{g'}{q'}(F)$} to $\Mod(\lie{g'\oplus g}, \Delta(G'))$.

\begin{Theorem}\label{thm:Exactness}
	Let $F$ be a finite-dimensional irreducible $\lie{l}$-module with infinitesimal character $[\lambda] \in \lie{t}^*/W_L$
	satisfying
	\begin{align*}
		\frac{2(\lambda+\rhog{\lie{u}}, \alpha)}{(\alpha, \alpha)} \not \in \set{0, 1, 2, \ldots}, \qquad \forall \alpha \in \Delta(\lie{u}, \lie{t}).
	\end{align*}
	Then the functor $T$ is exact and preserves submodule lattices.
\end{Theorem}

Note that `preserves submodule lattices' means that $T$ induces a lattice isomorphism from $\lattice(M)$ to $\lattice(T(M))$ for any $M \in \BGGcat{g'}{q'}(F)$.
The proof of Theorem~\ref{thm:Exactness} is postponed to Section~\ref{subsect:ProofEmbedding}.
We shall show a corollary assuming Theorem~\ref{thm:Exactness}.

\begin{Corollary}\label{cor:CategoryEquivalence}
	Under the assumption in Theorem~$\ref{thm:Exactness}$,
	the functor $T$ is fully faithful, and maps irreducible objects to irreducible ones.
\end{Corollary}

\begin{proof}
	Since $T$ preserves submodule lattices, the second assertion is clear.

	Let \smash{$M_1, M_2 \in \BGGcat{g'}{q'}(F)$}.
	First, we shall prove the faithfulness of $T$.
	Let $f \in \Hom_{\lie{g}'}(M_1, M_2)$ such that $T(f)=0$.
	By the exactness of $T$, we have two exact sequences
	\begin{align*}
		M_1 \xrightarrow{f} \Im(f) \rightarrow 0, \qquad T(M_1) \xrightarrow{T(f)} T(\Im(f)) \rightarrow 0.
	\end{align*}
	$T(f)=0$ implies $T(\Im(f))=0$.
	Since $T$ preserves submodule lattices, we obtain $\Im(f)=0$, and hence $f=0$.

	Next, we shall show that $T$ is full.
	Let $p_i$, $i = 1, 2$, denote the projection from $M_1 \oplus M_2$ to~$M_i$.
	Let $f$ be a morphism from $T(M_1)$ to $T(M_2)$.
	Consider
	\begin{align*}
		\id_{T(M_1)}\oplus f\colon\ T(M_1) \rightarrow T(M_1) \oplus T(M_2).
	\end{align*}
	$\Im\bigl(\id_{T(M_1)} \oplus f\bigr)$ is the graph of $f$.
	Since $T$ preserves submodule lattices, there exists a unique submodule $N \subset M_1 \oplus M_2$
	such that $T(N)=\Im\bigl(\id_{T(M_1)}\oplus f\bigr)$.
	Then $N$ is a graph.
	In fact, since~$T$ is faithful and $T(p_1|_N)$ is bijective,
	$p_1|_N$ is also bijective.

	We set $f'\coloneq p_2\circ(p_1|_N)^{-1}$.
	Then we have
 \[
 T(f')=T(p_2)\circ T\bigl((p_1|_N)^{-1}\bigr)=T(p_2)\circ \bigl(\id_{T(M_1)}\oplus f\bigr)=f.
 \]
	Hence $T$ is full.
	We have proved the corollary.
\end{proof}

The exactness in Theorem~\ref{thm:Exactness} is easy from the general result in the previous section.

\begin{Lemma}\label{lem:ExactnessPart}
	Let $F$ be a finite-dimensional irreducible $\lie{l}$-module.
	Assume that $\ind{g}{q}(F)$ is irreducible.
	Then \smash{$\Dzuck{\Delta(G')}{\Delta(L')}{d}\bigl(M\otimes \ind{g}{q}(F)\bigr) = 0$} holds
	for any \smash{$M \in \BGGcat{g'}{q'}(F)$} and $d \neq S$.
\end{Lemma}
	
\begin{proof}
	Let \smash{$M \in \BGGcat{g'}{q'}(F)$}.
	The assertion for $d < S$ follows from Propositions~\ref{prop:CompletelyReducibleVerma},~\ref{prop:FiltrationRestriction} and~\ref{prop:StandardFiltrationVanishing}.
	Let $d > S$.
	Since $M\otimes \ind{g}{q}(F)$ is $\lie{l'}$-admissible by Proposition~\ref{prop:AdmissibleDual}, we have
	\begin{align*}
		\bigl(\bigl(M\otimes \ind{g}{q}(F)\bigr)^*\bigr)_{\lie{l'}} \simeq (M^*)_{\lie{l'}} \otimes \ind{g}{\overline{q}}(F^*).
	\end{align*}
	As in the proof of Proposition~\ref{prop:Vanishing}, we obtain
	\begin{align*}
		&\Bigl(\Dzuck{\Delta(G')}{\Delta(L')}{d}\bigl(M\otimes \ind{g}{q}(F)\bigr)^*\Bigr)_{\Delta(G')} \\
		&\qquad{}\simeq \Dzuck{\Delta(G')}{\Delta(L')}{2S - d}\bigl((M^*)_{\lie{l'}} \otimes \ind{g}{\overline{q}}(F^*)\bigr) = 0.
		\qedhere
 \tag*{\qed}
 \end{align*}
 \renewcommand{\qed}{}
\end{proof}

\begin{proof}[Proof of the exactness in Theorem~\ref{thm:Exactness}]
	Note that $\ind{g}{q}(F)$ is irreducible under the assumption in Theorem~\ref{thm:Exactness} by Fact~\ref{fact:VermaIrreducible}.
	Hence the exactness follows from Lemma~\ref{lem:ExactnessPart} and the long exact sequence.
\end{proof}

\begin{Corollary}\label{cor:LatticeHom}
	Retain the setting in Lemma~$\ref{lem:ExactnessPart}$.
	Let~\smash{$M \in \BGGcat{g'}{q'}(F)$}.
	Then $T$ induces a lattice homomorphism $T^* \colon \lattice(M)\rightarrow \lattice(T(M))$.
\end{Corollary}

\begin{proof}
	The map $T^*$ is given by $T^*(N) = T(N)$, $N \in \lattice(M)$.
	Here we identify $T(N)$ with a submodule of $T(M)$ using the exactness of $T$.
	It is obvious that $T^*$ is order-preserving.
	Let $a\colon M\oplus M\rightarrow M$ denote the addition.
	Then $T(a)\colon T(M)\oplus T(M)\rightarrow T(M)$ is also the addition.
	For any two submodules $M_1, M_2 \subset M$, we have $T(M_1 \cap M_2) = T(M_1) \cap T(M_2)$ and $T(M_1 + M_2) = T(M_1)+T(M_2)$ from the following commutative diagram:
	\begin{align*}
		\begin{array}{ccccccccc}
			0 &\rightarrow& T(M_1 \cap M_2) &\rightarrow &T(M_1) \oplus T(M_2)& \rightarrow& T(M_1 + M_2) &\rightarrow& 0 \\
			 & & \rotatebox{90}{$\supset$} & & \rotatebox{90}{$\supset$} & & \rotatebox{90}{$\supset$} & & \\
			0 &\rightarrow& T(M) &\rightarrow &T(M) \oplus T(M)& \xrightarrow{T(a)}& T(M) &\rightarrow& 0.
		\end{array}
	\end{align*}
	This implies that $T^*$ is a lattice homomorphism.
\end{proof}

\subsection{Proof of Theorem~\ref{thm:Exactness}}\label{subsect:ProofEmbedding}

In this subsection, we prove the remaining part of Theorem~\ref{thm:Exactness}.
Retain the notation in the previous subsection.
We shall prove that $T^*\colon \lattice(M) \rightarrow \lattice(T(M))$ is bijective for any \smash{$M \in \BGGcat{g'}{q'}(F)$}.

Fix a simply-connected connected semisimple algebraic group $G$ with the Lie algebra $\lie{g}$.
Replacing $G'$ with its finite covering, we may assume that the homomorphism $G'\rightarrow \Ad(G)$
factors through a homomorphism $G'\rightarrow G$.
Then any $G$-module is a $(\lie{g}, G')$-module.
Remark that this operation does not affect Theorem~\ref{thm:Exactness}.

Let $[\mu] \in (\lie{t}')^*/W_{G'}$ be an infinitesimal character and $F$ a finite-dimensional $\lie{l'}$-module.
We~denote by \smash{$\BGGcat{g'}{q'}(F, [\mu])$} the full subcategory of \smash{$\BGGcat{g'}{q'}(F)$} consisting objects
with the generalized infinitesimal character $[\mu]$.
Since any object in \smash{$\BGGcat{g'}{q'}(F)$} has the primary decomposition, we have
\begin{align*}
	\BGGcat{g'}{q'}(F) = \bigoplus_{[\mu]} \BGGcat{g'}{q'}(F, [\mu]).
\end{align*}
Hence it is enough to show Theorem~\ref{thm:Exactness} replacing \smash{$\BGGcat{g'}{q'}(F)$} with \smash{$\BGGcat{g'}{q'}(F, [\mu])$}.

Note that \smash{$\BGGcat{g'}{q'}(F, [\mu])$} contains finitely many generalized Verma modules \smash{$\ind{g'}{q'}(F')$} up to isomorphism.
We denote by $\Lambda(F, [\mu])$ the set of infinitesimal characters of all finite-dimensional irreducible $\lie{l'}$-modules $F'$ such that \smash{$\ind{g'}{q'}(F') \in \BGGcat{g'}{q'}(F, [\mu])$}.
Then $[\mu' + \rhog{\lie{u'}}] = [\mu]$ holds for any~${[\mu'] \in \Lambda(F, [\mu])}$.

To show Theorem~\ref{thm:Exactness}, we reduce the problem to the case of enough large characters $F = \CC_\nu$ using the translation functor.
The reduction is standard and postponed to Theorem~\ref{thm:ReductionToScalar}.
First, we shall consider the case of enough large characters.

Let $[\mu] \in (\lie{t}')^*/W_{G'}$ be an infinitesimal character and $\nu$ a character of $\lie{l}$.
Then the infinitesimal character of $\CC_{\nu}$ is $[\nu + \rhog{\lie{l}}]$, where $\rhog{\lie{l}}$ is half the sum of positive roots of $\lie{l}$.
The choice of the set of positive roots is not important here.
Assume that
\begin{align}
	&\frac{2(\nu + \rhog{\lie{l}} + \rhog{\lie{u}}, \beta)}{(\beta, \beta)} \not\in \set{0, 1, 2, \ldots}, \label{eqn:EmbeddingDominant1} \\
	&\biggl(\mu' + \nu|_{\lie{t'}} + \rhog{\lie{u'}} + \sum_{\alpha \in E} \alpha, \beta'\biggr) < 0 \label{eqn:EmbeddingDominant2}
\end{align}
for any $[\mu'] \in \Lambda(\CC_{\nu}, [\mu])$, $\beta \in \Delta(\lie{u}, \lie{t})$, $\beta' \in \Delta(\lie{u'}, \lie{t'})$ and $E\subset \Delta(\lie{u}, \lie{t'})$.
We shall show Theorem~\ref{thm:Exactness} for $F = \CC_\nu$.
We use the notation $T$ in Theorem~\ref{thm:Exactness}.

We recall the bottom-layer map.
Let $F'$ be a finite-dimensional $\lie{l}'$-module such that $\smash{\ind{g'}{q'}(F')}\! \in \smash{\BGGcat{g'}{q'}(\CC_\nu, [\mu])}$.
Then we have the bottom-layer map
\begin{align*}
	B_{F'\otimes \CC_{\nu}}\colon\ \Dzuck{G'}{L'}{S}\bigl(\ind{g'}{q'}(F'\otimes \CC_{\nu})\bigr) \rightarrow \Dzuck{\Delta(G')}{\Delta(L')}{S}\bigl(\ind{g'}{q'}(F')\otimes \ind{g}{q}(\CC_{\nu})\bigr).
\end{align*}
For short, we write $B_{F'}$ for $B_{F'\otimes \CC_{\nu}}$.
Note that \smash{$\Dzuck{G'}{L'}{S}\bigl(\ind{g'}{q'}(F'\otimes \CC_{\nu})\bigr)$} is irreducible by Fact~\ref{fact:BorelWeilBott} and \eqref{eqn:EmbeddingDominant2}.

\begin{Lemma}\label{lem:PreserveIrreducible}
	Let \smash{$\ind{g'}{q'}(F') \in \BGGcat{g'}{q'}(\CC_\nu, [\mu])$} be a generalized Verma module.
	Then \smash{$T\bigl(\ind{g'}{q'}(F')\bigr)$}~is generated by $\Im(B_{F'})$.
	If, moreover, $M$ is the unique irreducible quotient of \smash{$\ind{g'}{q'}(F')$},
	then~$T(M)$ is the unique irreducible quotient of \smash{$T\bigl(\ind{g'}{q'}(F')\bigr)$}.
\end{Lemma}

\begin{proof}
	The assumption of Lemma~\ref{lem:CyclicSubspace} is fulfilled by \eqref{eqn:EmbeddingDominant2}.
	Hence the assertion is a special case of Lemma~\ref{lem:CyclicSubspace} and Corollary~\ref{cor:Irreducible}.
	Remark that $\ind{g}{q}(\CC_\nu)$ is irreducible by \eqref{eqn:EmbeddingDominant1} and Fact~\ref{fact:VermaIrreducible}.
\end{proof}

\begin{Lemma}\label{lem:Injectivity}
	Let \smash{$M_1, M_2 \in \BGGcat{g'}{q'}(\CC_\nu, [\mu])$} be irreducible objects.
	If $T(M_1) \simeq T(M_2)$, then $M_1 \simeq M_2$ holds.
\end{Lemma}

\begin{proof}
	It is well-known that $M_i$, $i=1,2$, is isomorphic to the unique irreducible quotient of a generalized Verma module \smash{$\ind{g'}{q'}\bigl(F'_i\bigr)$}, and hence $T(M_i)$ is isomorphic to the unique irreducible quotient of \smash{$T\bigl(\ind{g'}{q'}\bigl(F'_i\bigr)\bigr)$} by Lemma~\ref{lem:PreserveIrreducible}.
	Suppose $T(M_1) \simeq T(M_2)$.
	By Proposition~\ref{prop:Ktype} and Lemma~\ref{lem:PreserveIrreducible}, we have $\Im\bigl(B_{F'_0}\bigr) \simeq \Im\bigl(B_{F'_1}\bigr)$.
	Recall that \smash{$\Im\bigl(B_{F'_i}\bigr) \simeq \Dzuck{G'}{L'}{S}\bigl(\ind{g'}{q'}\bigl(F'_i\otimes \CC_{\nu}\bigr)\bigr)$}.
	By~Fact~\ref{fact:BorelWeilBott}, we obtain $F'_0 \simeq F'_1$ and hence $M_1 \simeq M_2$.
\end{proof}

\begin{Lemma}
	Let $M \in \BGGcat{g'}{q'}(\CC_\nu, [\mu])$.
	Then $T^*\colon \lattice(M) \rightarrow \lattice(T(M))$ defined in Corollary~$\ref{cor:LatticeHom}$ is bijective.
\end{Lemma}

\begin{proof}
	We have shown that $T^*$ is a lattice homomorphism in Corollary~\ref{cor:LatticeHom}.
	Since $T$ is exact and sends irreducible objects to irreducible ones by Lemma~\ref{lem:PreserveIrreducible}, $T(M)$ has the same length as~$M$.

	We shall show that $T^*$ is surjective by induction on the length of $M$.
	If $M$ is completely reducible (or, in particular, irreducible), then the surjectivity follows from Lemma~\ref{lem:Injectivity}.
	Assume that~$M$ has length $2$ and is not completely reducible, and let us show that $T(M)$ is indecomposable.
	Then $M$ or $(M^*)_{\lie{l}'}$ is a highest weight module.
	If $M$ is a highest weight module,~then~$T(M)$ is indecomposable by Lemma~\ref{lem:PreserveIrreducible}.
	Similarly, if $(M^*)_{\lie{l}'}$ is a highest weight~module, \smash{$\Dzuck{\Delta(G')}{\Delta(L')}{S}\bigl((M^*)_{\lie{l}'} \otimes \ind{g}{\overline{q}}(\CC_{-\nu})\bigr)$} is indecomposable, and hence its $\Delta(G')$-finite dual $T(M)$ is also indecomposable.

	Assume that $M$ has length strictly greater than $2$.
	Let $N\subset M$ and $L\subset T(M)$ be irreducible submodules.
	By the induction hypothesis, we have $(L+T(N))/T(N) \in T^*(\lattice(M/N))$.
	This implies that there exists a submodule $L'\subset M$ such that $T(L') = L + T(N)$.
	Since the length of~$L'$ is less than or equal to $2$, we have $L \in \lattice(T(L')) = T^*(\lattice(L'))\subset T^*(\lattice(M))$ by the induction hypothesis.

	Let $L\subset T(M)$ be a non-irreducible submodule.
	By the irreducible case, there exists~an irreducible submodule $N \subset M$ such that $T(N)\subset L$.
	By the induction hypothesis, we have $L/T(N) \in T^*(\lattice(M/N))$, and hence $L \in T^*(\lattice(M))$.
\end{proof}

We have proved Theorem~\ref{thm:Exactness} for $F = \CC_{\nu}$ with the assumptions \eqref{eqn:EmbeddingDominant1} and \eqref{eqn:EmbeddingDominant2}.
To reduce the general cases, we need the following lemma, which assures the existence of such a character~$\nu$.

\begin{Lemma}\label{lem:PositivityRho}
	One has
	\begin{alignat*}{3}
		&(\rhog{\lie{u}}, \alpha)> 0, \qquad&& \forall \alpha \in \Delta(\lie{u}, \lie{t}),&\\
		&(\rhog{\lie{u}}|_{\lie{t'}}, \alpha)> 0, \qquad&& \forall \alpha \in \Delta(\lie{u}', \lie{t}').&
	\end{alignat*}
\end{Lemma}

Lemma~\ref{lem:PositivityRho} is a special case of the following proposition.
Only here, let $\lie{g}$ be a complex reductive Lie algebra and $\lie{t}$ a Cartan subalgebra of $\lie{g}$.

\begin{Proposition}\label{prop:PositivityRho}
	Let $V$ be a finite-dimensional $\lie{g}$-module and $\alpha \in \Delta(\lie{g}, \lie{t})$.
	Let $W$ be a $\lie{t}$-submodule of $V$ such that $\lie{g}_{\alpha} \cdot W \subset W$.
	Then \smash{$\sum_{\beta \in \Delta(W, \lie{t})}(\beta, \alpha) \geq 0$} holds.
	If, moreover, $\lie{g}\subset V$, $\lie{g}_{\alpha} \subset W$ and $\lie{g}_{-\alpha} \not\subset W$,
	then \smash{$\sum_{\beta \in \Delta(W, \lie{t})}(\beta, \alpha) > 0$} holds.
\end{Proposition}

\begin{Remark}
	In Proposition~\ref{prop:PositivityRho}, we consider $\Delta(W, \lie{t})$ as a multiset.
\end{Remark}
	
\begin{proof}[Proof of Proposition~\ref{prop:PositivityRho}]
	Fix an $\lie{sl}_2$-triple $\set{H_\alpha, X_\alpha, Y_\alpha}$ for the root $\alpha$ such that $X_\alpha \in \lie{g}_{\alpha}$ and $Y_\alpha \in \lie{g}_{-\alpha}$.
	Set $\lie{s}\coloneq \spn{H_\alpha, X_\alpha, Y_\alpha}$.

	Since $V|_{\lie{s}}$ is completely reducible, $W$ has the eigenspace decomposition
	$W = \bigoplus_{i = -N}^N W_i$ for~$H_\alpha$, where $W_i$ is zero or the eigenspace of eigenvalue $i \in \ZZ$.
	Since $W$ is $X_\alpha$-stable, the map $X_\alpha^i\cdot \colon W_{-i} \rightarrow W_{i}$
	is well-defined and injective for each $i \in \NN$ by the representation theory of~$\lie{sl}(2,\CC)$.
	Hence we obtain the desired inequality from
	\begin{align*}
		\frac{2}{(\alpha, \alpha)} \sum_{\beta \in \Delta(W, \lie{t})} (\beta, \alpha) = \sum_{i=-N}^N i\dim(W_i)
	 = \sum_{i=1}^N i(\dim(W_i) - \dim(W_{-i})) \geq 0.
	\end{align*}
	
	The last assertion follows from that $X_\alpha^2\cdot \colon W_{-2} \rightarrow W_{2}$ is not surjective
	if $\lie{g}_{\alpha} \subset W$ and $\lie{g}_{-\alpha} \not\subset W$.
\end{proof}

\begin{Theorem}\label{thm:ReductionToScalar}
	Let $F$ be a finite-dimensional irreducible $\lie{l}$-module with infinitesimal character~$[\lambda]$ satisfying the assumption in Theorem~$\ref{thm:Exactness}$,
	and \smash{$M \in \BGGcat{g'}{q'}(F)$}.
	Then the lattice homomorphism $T^*\colon \lattice(M) \rightarrow \lattice(T(M))$ defined in Corollary~$\ref{cor:LatticeHom}$ is bijective.
\end{Theorem}

\begin{proof}
	Since $M$ has the primary decomposition, we may assume that $M$ has a generalized infinitesimal character.
	Fix a set $\Delta^+(\lie{l}, \lie{t})$ of positive roots of $\Delta(\lie{l}, \lie{t})$.
	Then $\Delta^+(\lie{g}, \lie{t})\coloneq\Delta^+(\lie{l}, \lie{t}) \cup \Delta(\lie{u}, \lie{t})$ is a set of positive roots
	of $\Delta(\lie{g}, \lie{t})$.
	Replacing $\lambda$ with some representative of $[\lambda]$, we may assume that $\lambda$ is regular anti-dominant with respect to $\Delta^+(\lie{l}, \lie{t})$.
	By assumption, $\lambda + \rhog{\lie{u}}$ is regular and integrally anti-dominant (see Definition~\ref{def:IntegrallyDominant})
	with respect to $\Delta^+(\lie{g}, \lie{t})$.

	Set $\lie{l}_{ss}\coloneq[\lie{l}, \lie{l}]$
	and $\lie{t}_{ss}\coloneq \lie{t}\cap \lie{l}_{ss}$.
	Since $-(\lambda+\rhog{\lie{l}})|_{\lie{t}_{ss}}$ is an algebraically integral weight of $\lie{l}_{ss}$
	and $\lie{l}$ is a Levi subalgebra of $\lie{g}$,
	there exists an algebraically integral weight $\mu \in \lie{t}^*$ of $\lie{g}$
	such that $\mu|_{\lie{t}_{ss}}=-(\lambda+\rhog{\lie{l}})|_{\lie{t}_{ss}}$.
	Set $\nu\coloneq\lambda+\mu+\rhog{\lie{l}}$.
	Then we have \smash{$T_{\lambda}^{\lambda+\mu}(F) \simeq \CC_{\nu}$}.
	It is easy to see \smash{$\BGGcat{g'}{q'}(F)=\BGGcat{g'}{q'}(\CC_{\nu})$} because the $\lie{l}'$-module $F\otimes \CC_{-\nu}$ lifts to an $L'$-module.
	Remark that we have assumed that $G$ is simply-connected.

	By Lemma~\ref{lem:PositivityRho}, if necessary, replacing $\mu$ by $\mu-2m\rhog{\lie{u}}$ for enough large integer $m$,
	we may assume that $\nu=\lambda+\mu+\rhog{\lie{l}}$ satisfies \eqref{eqn:EmbeddingDominant1} and \eqref{eqn:EmbeddingDominant2}.
	Using the translation functors and by Facts~\ref{fact:TranslationVerma} and~\ref{fact:BasicZuckerman},
	we have
	\begin{align*}
		&T_{\lambda+\rhog{\lie{u}}}^{\lambda+\mu+\rhog{\lie{u}}}\bigl(\ind{g}{q}(F)\bigr)
 \simeq \ind{g}{q}(\CC_{\nu}),\\
		&T_{\lambda+\rhog{\lie{u}}}^{\lambda+\mu+\rhog{\lie{u}}}\Bigl(\Dzuck{\Delta(G')}{\Delta(L')}{S}\bigl(M\otimes \ind{g}{q}(F)\bigr)\Bigr)
 \simeq \Dzuck{\Delta(G')}{\Delta(L')}{S}\bigl(M\otimes \ind{g}{q}(\CC_{\nu})\bigr), \\
		&T_{\lambda+\mu+\rhog{\lie{u}}}^{\lambda+\rhog{\lie{u}}}\Bigl(\Dzuck{\Delta(G')}{\Delta(L')}{S}\bigl(M\otimes \ind{g}{q}(\CC_{\nu})\bigr)\Bigr)
 \simeq \Dzuck{\Delta(G')}{\Delta(L')}{S}\bigl(M\otimes \ind{g}{q}(F)\bigr).
	\end{align*}
	Note that the isomorphisms are natural in $M$.

	We have shown that \smash{$\Dzuck{G'}{L'}{S}\bigl(M\otimes \ind{g}{q}(\CC_{\nu})\bigr)$} has finite length in Lemma~\ref{lem:Injectivity}.
According to Fact~\ref{fact:TranslationEquvalence}, the~translation functors \smash{$T^{\lambda+\rhog{\lie{u}}}_{\lambda+\mu+\rhog{\lie{u}}}$} and \smash{$T^{\lambda+\mu+\rhog{\lie{u}}}_{\lambda+\rhog{\lie{u}}}$} give equivalences of categories on the images of the functors $T_F$ and $T_{\CC_{\nu}}$.
	Therefore, the assertion follows from Lemma~\ref{lem:Injectivity}.
\end{proof}

\subsection{Equivalence of categories in special case}

We shall consider the case in which $\lie{q}$ is a Borel subalgebra and $\lie{g'} = \lie{g}$.
We show that the functor~$T_F$ defined in Theorem~\ref{thm:Exactness} gives an equivalence of categories in this setting.

Let $G$ be a simply-connected connected semisimple algebraic group and $T$ a maximal torus of $G$.
Fix a Borel subalgebra $\lie{b}$ of $\lie{g}$ with Levi decomposition $\lie{b}=\lie{t}\oplus \lie{u}$, where $\lie{u}$ is the nilpotent radical of $\lie{b}$.
Then $(\lie{g\oplus g}, \Delta(G))$ is a pair.
Set $S\coloneq \dim(\lie{u})$.

\newcommand{\functF}{\mathcal{F}}
\newcommand{\functG}{\mathcal{G}}

We need the following equivalence of categories proved in~\cite[Theorem~5.9]{BeGe80_projective_functor}.
Let $\lambda \in \lie{t}^*$.
Two covariant functors $\functF$ and $\functG$ are defined by
\begin{align*}
	&\functF(M)\coloneq \Hom_{\CC}\bigl(\ind{g}{b}(\CC_{\lambda}), M\bigr)_{\Delta(G)}, \qquad M \in \BGGcat{g}{b}(\CC_{-\lambda}), \\
	&\functG(N)\coloneq N \otimes_{\univ{g}} \ind{g}{b}(\CC_{\lambda}), \qquad N \in \Mod(\lie{g\oplus g}, \Delta(G))_{-\lambda-\rhog{\lie{u}}}.
\end{align*}
Here, to define $\functG$, we consider the $(\lie{g}\oplus \lie{g}, \Delta(G))$-module $N$
as a $(\univ{g}, \univ{g})$-bimodule via
\begin{align*}
	A\cdot n \cdot B \coloneq \bigl(A \otimes {}^t\!B\bigr)n, \qquad A,B \in \univ{g} \text{ and }n \in N,
\end{align*}
and $\Mod(\lie{g\oplus g}, \Delta(G))_{-\lambda-\rhog{\lie{u}}}$ is the full subcategory of $\Mod(\lie{g\oplus g}, \Delta(G))$
consisting of objects with the infinitesimal character $-\lambda-\rhog{\lie{u}}$ with respect to the subalgebra $0\oplus \lie{g}\subset \lie{g\oplus g}$.
For the functor $(\cdot)_{\Delta(G)}$, see Definition~\ref{def:GKmod}.
The following fact is proved in~\cite[Theorem~5.9]{BeGe80_projective_functor}.

\begin{Fact}\label{fact:BGequivalence}
	If $\lambda+\rhog{\lie{u}}$ is regular and integrally dominant $($see Definition $\ref{def:IntegrallyDominant})$,
	then $\functF$ gives an equivalence of categories and $\functG$ is a quasi-inverse of $\functF$.
\end{Fact}

\begin{Theorem}\label{thm:CategoryEquivalence}
	Let $\lambda \in \lie{t}^*$ and assume that $\lambda + \rhog{\lie{u}}$ is regular and integrally anti-dominant.
	Then the functor $\mathcal{T}\colon \BGGcat{g}{b}(\CC_{\lambda}) \rightarrow \Mod(\lie{g\oplus g}, \Delta(G))_{\lambda+\rhog{\lie{u}}}$
	defined by $\mathcal{T}(M) = \Dzuck{\Delta(G)}{\Delta(T)}{S}\bigl(M \otimes \ind{g}{b}(\CC_{\lambda})\bigr)$ gives an equivalence of categories.
\end{Theorem}

\begin{Remark}
	This theorem was conjectured and partially proved by T.J.~Enright in~\cite{En81_lecture}.
\end{Remark}

\begin{proof}[Proof of Theorem~\ref{thm:CategoryEquivalence}]
	By Theorem~\ref{thm:Exactness} and Corollary~\ref{cor:CategoryEquivalence}, the functor $\mathcal{T}$ is exact and fully faithful.
	Hence what we need to show is that $\mathcal{T}$ is dense (i.e., for any object $N$ in the codomain,
	there is an object~$M$ in the domain with $\mathcal{T}(M) \simeq N$).
	In fact, a functor $F$ gives an equivalence of categories if and only if $F$ is fully faithful and dense.

	We set $\lambda'=-\lambda-2\rhog{\lie{u}}$.
	Then we have \smash{$\BGGcat{g}{b}(\CC_{-\lambda'})=\BGGcat{g}{b}(\CC_{\lambda})$},
	and $\lambda' + \rhog{\lie{u}}$ is regular and integrally dominant.
	We consider $\functG$ for the weight \smash{$\lambda'$}.
	By Theorem~\ref{thm:Exactness} and Fact~\ref{fact:BGequivalence},
	we have an exact fully faithful endofunctor $\functG\circ \mathcal{T}$ on \smash{$\BGGcat{g}{b}(\CC_\lambda)$},
	and $\functG\circ \mathcal{T}$ preserves the irreducibility.

	We shall prove that $\functG\circ \mathcal{T}$ is a dense functor.
	We denote by $\Irr$ the set of all isomorphism classes of irreducible objects in $\BGGcat{g}{b}(\CC_\lambda)$.
	Then $\functG \circ \mathcal{T}$ induces a permutation on $\Irr$.
	Since $\functG\circ \mathcal{T}$ preserves infinitesimal characters,
	the cardinality of each orbit on $\Irr$ is bounded by $|W_G|$.
	Hence there exists a positive integer $k$ such that $(\functG\circ T)^k$ acts on $\Irr$ trivially.
	Obviously, $(\functG\circ \mathcal{T})^k$ is dense if and only if $\functG\circ \mathcal{T}$ is dense.

	Set $\mathcal{E}\coloneq(\functG\circ \mathcal{T})^k$ and we shall show that $\mathcal{E}$ is dense.
	Let $P$ be an indecomposable projective object in $\BGGcat{g}{b}(\CC_\lambda)$.
	Then $P$ is a projective cover of an irreducible object $L$ (see~\cite[Section 3.9]{Hu08_category_o}).
	Let $\pi \colon P\rightarrow L$ denote the surjective morphism.
	Remark that $L$ is a unique irreducible quotient of $P$.
	Since $\mathcal{E}$ is exact, $\mathcal{E}(\pi)\colon \mathcal{E}(P) \rightarrow \mathcal{E}(L)$ is surjective.
	Since $\mathcal{E}$ acts on $\Irr$ trivially, we can identify $\mathcal{E}(L)$ with $L$.
	Since $P$ is projective, there exists a homomorphism $\tau\colon P \rightarrow \mathcal{E}(P)$ such
	that the following diagram commutes:
	\begin{align*}
		\xymatrix{
		& P\ar[dl]_{\tau}\ar[d]^\pi& \\
		\mathcal{E}(P)\ar[r]_{E(\pi)}& \mathcal{E}(L) = L\ar[r]& 0 \text{ (exact).}
		}
	\end{align*}
	Recall that the functor $\mathcal{E}$ induces a lattice isomorphism between $\lattice(P)$ and $\lattice(\mathcal{E}(P))$.
	This implies that $\mathcal{E}(L)$ is a unique irreducible quotient of $\mathcal{E}(P)$, and hence $\tau$ is surjective.
	Since $P$ and $\mathcal{E}(P)$ have the same length, we have $P \simeq \mathcal{E}(P)$.

	Since $\BGGcat{g}{b}(\CC_\lambda)$ has enough projectives~\cite[Section 3.9]{Hu08_category_o},
	any object of $\BGGcat{g}{b}(\CC_\lambda)$
	can be written as a quotient of a projective object, which is a direct sum of indecomposable projectives.
	Therefore~$\mathcal{E}$ is dense.
	This completes the proof.
\end{proof}

\subsection{Embedding 2}

We use the notation in Section~\ref{subsect:EmbeddingSetting}.
We shall show a similar result to Theorem~\ref{thm:Exactness} swapping the roles of $\lie{g}$ and $\lie{g'}$.

Let $F'$ be a finite-dimensional irreducible $\lie{l'}$-module.
We denote by $\BGGcat{g}{q}(F')$ the full subcategory of $\BGGcat{g}{q}$ whose object $M$ satisfies that $M\otimes F'$ lifts to an $L'$-module.
Let $T'$ denote the~func\-tor~\smash{$\Dzuck{\Delta(G')}{\Delta(L')}{S}\bigl(\ind{g'}{q'}(F') \otimes \cdot\bigr)$} from $\BGGcat{g}{q}(F')$ to $\Mod(\lie{g'\oplus g}, \Delta(G'))$.

\begin{Theorem}\label{thm:Embedding2}
	Let $F'$ be a finite-dimensional irreducible $\lie{l'}$-module with infinitesimal character $[\lambda'] \in (\lie{t}')^*/W_{L'}$ satisfying
	\begin{align*}
		\frac{2(\lambda'+\rhog{\lie{u'}}, \alpha)}{(\alpha, \alpha)} \not \in \set{0, 1, 2, \ldots}\qquad \text{for any }\alpha \in \Delta(\lie{u'}, \lie{t'}).
	\end{align*}
	Then the functor $T'$ is exact and preserves submodule lattices.
\end{Theorem}

\begin{proof}
	Let $M \in \BGGcat{g}{q}(F')$.
	We have assumed that $\lie{q}$ is weakly $\lie{g'}$-compatible.
	By Propositions~\ref{prop:FiltrationRestriction} and~\ref{prop:AdmissibleDual}, $M|_{\lie{g'}}$ is discretely decomposable in \smash{$\BGGcat{g'}{q'}(F')$}, and each primary component~\smash{$P_\chi(M|_{\lie{g'}})$}~has finite length.
	This implies that $M|_{\lie{g'}}$ is a direct sum of objects in \smash{$\BGGcat{g'}{q'}(F')$}.
	By Theorem~\ref{thm:Exactness} for~${\lie{g} = \lie{g'}}$, $T'$ is exact and faithful.

	To show that $T'$ preserves submodule lattices, let $L$ be a $(\lie{g' \oplus g}, \Delta(G'))$-submodule of $T'(M)$.
	By Theorem~\ref{thm:Exactness} and the above discussion, there exists a unique $\lie{g'}$-submodule $L' \subset M$ such that \smash{$\Dzuck{\Delta(G')}{\Delta(L')}{S}\bigl(\ind{g'}{q'}(F') \otimes L'\bigr) = L$}.
	We shall show that $L'$ is $\lie{g}$-stable.

	For a direct sum $N$ of objects in \smash{$\BGGcat{g'}{q'}(F')$}, set \smash{$T'(N) = \Dzuck{\Delta(G')}{\Delta(L')}{S}\bigl(\ind{g'}{q'}(F') \otimes N\bigr)$} by abuse of notation.
	By~\cite[Lemma~6.3.1]{Wa88_real_reductive_I}, we have a commutative diagram
	\begin{align*}
		\xymatrix{
			W\otimes T'(L') = W\otimes L \ar[r]^-m \ar[d]^{\simeq} & T'(M) \ar@{=}[d]\\
			T'(W\otimes L') \ar[r]^-{T'(m)} & T'(M)
		}
	\end{align*}
	for any finite-dimensional $G$-submodule $W\subset \univ{g}$, where $m$'s are the multiplication maps.
	This implies that $T'(\univ{g}L') = \univ{g}T'(L') = L$
	and hence $\univ{g}L' = L'$ by the uniqueness of $L'$.
	Therefore, $L'$ is $\lie{g}$-stable and this shows the assertion.
\end{proof}

\section[Irreducibility of U(g)\^{}{G'}-modules]{Irreducibility of $\boldsymbol{\univ{g}^{G'}}$-modules}

In this section, we consider the irreducibility of $\univ{g}^{G'}$-modules in the branching problem.
We treat generalized Verma modules, cohomologically induced modules and discrete series representations.

\subsection{Generalized Verma module}

In this subsection, we consider the branching problem of generalized Verma modules.
We will give a criterion for the irreducibility of $\univ{g}^{G'}$-modules on $\Hom$ spaces.

Let $(\lie{g}, G')$ be a pair with semisimple $\lie{g}$.
Let $\lie{q} = \lie{l}\oplus \lie{u}$ be a weakly $\lie{g'}$-compatible parabolic subalgebra of $\lie{g}$
containing $\lie{q}(H)$ ($H \in \lie{g'}$).
See Definition~\ref{def:WeaklyCompatible}.
Set
\begin{align*}
	\lie{l'}\coloneq \lie{l}\cap \lie{g'}, \qquad \lie{u'}\coloneq \lie{u}\cap \lie{g'}, \qquad \lie{\overline{u}'}\coloneq \lie{\overline{u}}\cap \lie{g'}, \qquad \lie{q'}\coloneq \lie{l'}\oplus \lie{u'}.
\end{align*}
Then $\lie{q'}$ is a parabolic subalgebra of $\lie{g'}$.
Write $L'$ for the centralizer of $H$ in $G'$.

Fix a Cartan subalgebra $\lie{t}$ of $\lie{l}$ and write $W_{G}$ and $W_L$ for the Weyl groups.
Let $\rhog{\lie{u}}$ denote half the sum of roots in $\Delta(\lie{u}, \lie{t})$.
Similarly, we define $\lie{t'}$, $W_{G'}$, $W_{L'}$ and $\rhog{\lie{u'}}$ for $\lie{g'}$.
Fix a~set $\Delta^+(\lie{l'}, \lie{t'})$ of positive roots in $\Delta(\lie{l'}, \lie{t'})$
and $\rhog{\lie{l'}}$ denotes half the sum of the positive roots.
Set~${\rhog{\lie{g'}} \coloneq \rhog{\lie{l'}} + \rhog{\lie{u'}}}$.

In this subsection, using the results in Section~\ref{sect:Embedding}, we consider the irreducibility of the $\univ{g}^{G'}$-module
$\Hom_{\lie{g'}}\bigl(\ind{g'}{q'}(F'), \ind{g}{q}(F)|_{\lie{g'}}\bigr)$.
For a finite-dimensional irreducible $\lie{l}'$-module $F$, we set
\begin{align*}
	\fdual{F} \coloneq F^* \otimes \CC_{-2\rhog{\lie{u'}}}.
\end{align*}
Then the multiplicity of $\CC_{-2\rhog{\lie{u'}}}$ in $\fdual{F}\otimes F$ is $1$, and
\begin{align*}
	(-2\rhog{\lie{u'}} - \rhog{\lie{l'}} + \rhog{\lie{u'}}, \beta) = (-\rhog{\lie{l'}} - \rhog{\lie{u'}}, \beta) < 0, \qquad \forall \beta \in \Delta(\lie{u'}, \lie{t'}).
\end{align*}
Note that $-2\rhog{\lie{u'}} - \rhog{\lie{l'}}$ is the infinitesimal character of $\CC_{-2\rhog{\lie{u'}}}$
and $-\rhog{\lie{l'}} - \rhog{\lie{u'}} = -\rhog{\lie{g'}}$ is that of the $\lie{g'}$-module $\CC$.

\begin{Lemma}\label{lem:UniqueTrivialCharacter}
	Let $F$ be a finite-dimensional irreducible $\lie{l}'$-module and $F'$ an irreducible $\lie{l}'$-submodule of $\fdual{F}\otimes F$.
	If \smash{$\ind{g'}{q'}(F')$} has the infinitesimal character $[\rhog{\lie{g'}}]$, then $F' \simeq \CC_{-2\rhog{\lie{u'}}}$ holds.
\end{Lemma}

\begin{proof}
	Write $\lambda\in (\lie{t'})^*$ for the representative of the infinitesimal character of $F'$ such that~$\lambda$ is anti-dominant with respect to $\Delta^+(\lie{l'}, \lie{t'})$.
	Then $[\lambda + \rhog{\lie{u'}}]$ is the infinitesimal character of~\smash{$\ind{g'}{q'}(F')$}.
	By assumption, there exists $w \in W_{G'}$ such that $w(\lambda + \rhog{\lie{u'}}) = -\rhog{\lie{g'}}$.
	Since~$F'$ is a~submodule of $\fdual{F}\otimes F$, we have $\lambda(H) = (-2\rhog{\lie{u'}} - \rhog{\lie{l'}})(H)$, and hence
	\begin{align}
		(\lambda + \rhog{\lie{u'}})(H) = -\rhog{\lie{g'}}(H). \label{eqn:UniqueTrivialCharacter}
	\end{align}

	Let $\lie{b'}$ denote the Borel subalgebra of $\lie{g'}$ such that $-\rhog{\lie{g'}}$ is dominant with respect to $\lie{b'}$.
	Then $\lambda + \rhog{\lie{u'}}$ is dominant with respect to $w^{-1}\lie{b'}$.
	This implies that $\lambda + \rhog{\lie{u'}}$ is half the sum of roots in $\Delta\bigl(w^{-1}\lie{b'}, \lie{t'}\bigr)$.
	By \eqref{eqn:UniqueTrivialCharacter}, we have $\lie{\overline{u}'} \subset w^{-1}\lie{b'}$.
	Since $\lambda$ is dominant with respect to $-\Delta^+(\lie{l'}, \lie{t'})$, this implies $w^{-1}\lie{b'} = \lie{b'}$
	and hence $w = e$.
	Therefore, we obtain $\lambda = -2\rhog{\lie{u'}} - \rhog{\lie{l'}}$ and $F'\simeq \CC_{-2\rhog{\lie{u'}}}$.
\end{proof}

By Lemma~\ref{lem:UniqueTrivialCharacter} and the proof of Lemma~\ref{lem:MinimalTypeOne}, we obtain the following.

\begin{Lemma}\label{lem:ExistenceInvariantVector}
	Let $F$ be a finite-dimensional irreducible $\lie{l}'$-module and $M$ a non-zero quotient of~${\ind{g'}{q'}\bigl(\fdual{F}\bigr)\otimes \ind{g'}{q'}(F)}$.
	\begin{enumerate}\itemsep=0pt
		\item[$1.$] \smash{$P_{-\rhog{\lie{l'}} - \rhog{\lie{u'}}}\bigl(M|_{\Delta(\lie{g'})}\bigr) \simeq \ind{g'}{q'}\bigl(\CC_{-2\rhog{\lie{u'}}}\bigr)$}.
		\item[$2.$] \smash{$\Dzuck{\Delta(G')}{\Delta(L')}{S}(M)^{\Delta(G')} \simeq \CC$}.
		\item[$3.$] \smash{$\Dzuck{\Delta(G')}{\Delta(L')}{i}(M)^{\Delta(G')} = 0$} for any $i \neq S$.
	\end{enumerate}
\end{Lemma}

By Lemma~\ref{lem:ExistenceInvariantVector}, the functor \smash{$\Dzuck{\Delta(G')}{\Delta(L')}{S}\bigl(\ind{g'}{q'}\bigl(\fdual{F}\bigr)\otimes (\cdot)\bigr)^{\Delta(G')}$} is regarded as an analogue of the functor \smash{$\Hom_{\lie{g'}}\bigl(\ind{g'}{q'}(F), \cdot\bigr)$}.
We shall compare the two functors.

Let $F'$ a finite-dimensional irreducible $\lie{l'}$-module and $M \in \BGGcat{g}{q}$.
Suppose that the $\lie{l}'$-module $M \otimes (F')^*$ lifts to an $L'$-module.
This implies that the $\lie{l}'$-module \smash{$\ind{g'}{q'}\bigl(\fdual{F'}\bigr)\otimes M$} lifts to an~\mbox{$L'$-module}.
Remark that this condition follows from \smash{$\Hom_{\lie{g'}}\bigl(\ind{g'}{q'}(F'), M\bigr) \neq 0$} if $M$ is indecomposable.

There is a canonical $\lie{g'}$- and $\univ{g}^{G'}$-homomorphism
\begin{align*}
\Phi \colon\ \ind{g'}{q'}(F') \otimes \Hom_{\lie{g'}}\bigl(\ind{g'}{q'}(F'), M\bigr)
\rightarrow M
\end{align*}
defined by $\Phi(v \otimes \varphi) = \varphi(v)$.
Let $M'$ be a non-zero quotient of \smash{$\ind{g'}{q'}\bigl(\fdual{F'}\bigr)$}.
Applying Lem\-ma~\ref{lem:ExistenceInvariantVector}\,(2) to \smash{$M'\otimes \ind{g'}{q'}(F')$}, we have isomorphisms
\begin{align}
	&\Hom_{\lie{g'}}\bigl(\ind{g'}{q'}(F'), M\bigr) \simeq \CC\otimes \Hom_{\lie{g'}}\bigl(\ind{g'}{q'}(F'), M\bigr) \nonumber \\
	&\qquad{}\simeq \Dzuck{\Delta(G')}{\Delta(L')}{S}\bigl(M'\otimes \ind{g'}{q'}(F')\bigr)^{\Delta(G')}\otimes \Hom_{\lie{g'}}\bigl(\ind{g'}{q'}(F'), M\bigr) \nonumber \\
	&\qquad{}\simeq \Dzuck{\Delta(G')}{\Delta(L')}{S}\bigl(M'\otimes \ind{g'}{q'}(F') \otimes \Hom_{\lie{g'}}\bigl(\ind{g'}{q'}(F'), M\bigr)\bigr)^{\Delta(G')} \label{eqn:IsomHomZuck}
\end{align}
of $\univ{g}^{G'}$-modules.
Compositing this isomorphism with \smash{$\Dzuck{\Delta(G')}{\Delta(L')}{S}(\id\otimes \Phi)^{\Delta(G')}$} gives a~$\univ{g}^{G'}$-homomorphism
\begin{align*}
	\Phi_{M'} \colon\ \Hom_{\lie{g'}}\bigl(\ind{g'}{q'}(F'), M\bigr)\rightarrow \Dzuck{\Delta(G')}{\Delta(L')}{S}(M'\otimes M)^{\Delta(G')}.
\end{align*}
If $\Phi_{M'}$ is injective, we can study the space $\Hom_{\lie{g'}}\bigl(\ind{g'}{q'}(F'), M\bigr)$ using the Zuckerman derived functor module.
We shall consider the injectivity of $\Phi_{M'}$.

\begin{Lemma}\label{lem:VanishingInvariantKernel}
	\smash{$\Dzuck{\Delta(G')}{\Delta(L')}{i}(M'\otimes \Ker(\Phi))^{\Delta(G')} = 0$} holds for any $i$.
\end{Lemma}

\begin{proof}
	Let $N$ denote the unique maximal submodule of $\ind{g'}{q'}(F')$.
	By Lemma~\ref{lem:ExistenceInvariantVector}, it is enough to show
	\begin{align}
	\Ker(\Phi) \subset N \otimes \Hom_{\lie{g'}}\bigl(\ind{g'}{q'}(F'), M\bigr). \label{eqn:VanishingInvariantKernel}
	\end{align}

	Let $\mu$ denote the highest weight of $\ind{g'}{q'}(F')$.
	Assume that \eqref{eqn:VanishingInvariantKernel} does not hold.
	Set $X \coloneq \Hom_{\lie{g'}}\bigl(\ind{g'}{q'}(F'), M\bigr)$.
	Then we have
	\begin{align*}
		0\neq \Ker(\Phi)/(\Ker(\Phi) \cap (N \otimes X)) \subset \bigl(\ind{g'}{q'}(F') / N\bigr) \otimes X.
	\end{align*}
	Since \smash{$\bigl(\ind{g'}{q'}(F')/N\bigr) \otimes X$} is isomorphic to a direct sum of some copies of the irreducible $\lie{g'}$-mod\-ule $\ind{g'}{q'}(F')/N$ as a $\lie{g'}$-module, there is a~weight vector $v \in \Ker(\Phi)$ with the weight $\mu$.
	Then $v$ is a~$\lie{b'}$-eigenvector.

	Write $v= \sum_{i=0}^r v_i\otimes \varphi_i$, where \smash{$\set{v_i}_{i=0}^r \subset \ind{g'}{q'}(F')$}
	and $\set{\varphi_i}_{i=0}^r \subset X$
	such that~$\set{v_i}_{i=0}^r$ and~$\set{\varphi_i}_{i=0}^r$ are linearly independent.
	Then each $v_i$ is a highest weight vector of \smash{$\ind{g'}{q'}(F')$}.
	This implies $r=0$.
	By $v \in \Ker(\Phi)$, we have $\varphi_0(v_0)=0$.
	Since $v_0$ is a cyclic vector of \smash{$\ind{g'}{q'}(F')$}, we have $\varphi_0 = 0$ and hence $v = 0$.
	This contradicts to $v \neq 0$.
	Therefore, we have shown \eqref{eqn:VanishingInvariantKernel}.
\end{proof}

Consider the following two exact sequences:
\begin{align*}
	&0 \rightarrow \Ker(\Phi)\rightarrow \ind{g'}{q'}(F') \otimes \Hom_{\lie{g'}}\bigl(\ind{g'}{q'}(F'), M\bigr) \xrightarrow{\Phi} \Im(\Phi) \rightarrow 0, \\
	&0 \rightarrow \Im(\Phi)\rightarrow M \rightarrow \Coker(\Phi) \rightarrow 0.
\end{align*}
Applying the functor \smash{$\Dzuck{\Delta(G')}{\Delta(L')}{i}(M'\otimes (\cdot))^{\Delta(G')}$}, and by Lemma~\ref{lem:VanishingInvariantKernel} and \eqref{eqn:IsomHomZuck}, we obtain
\begin{align*}
	&\Dzuck{\Delta(G')}{\Delta(L')}{S}(M'\otimes \Im(\Phi))^{\Delta(G')}\simeq \Hom_{\lie{g'}}\bigl(\ind{g'}{q'}(F'), M\bigr), \\
	&\Dzuck{\Delta(G')}{\Delta(L')}{S - 1}(M'\otimes \Coker(\Phi))^{\Delta(G')}
 \rightarrow \Dzuck{\Delta(G')}{\Delta(L')}{S}(M'\otimes \Im(\Phi))^{\Delta(G')} \\
	&\hphantom{\Dzuck{\Delta(G')}{\Delta(L')}{S - 1}(M'\otimes \Coker(\Phi))^{\Delta(G')}}{}
 \xrightarrow{\Phi_{M'}} \Dzuck{\Delta(G')}{\Delta(L')}{S}(M'\otimes M)^{\Delta(G')} \qquad (\text{exact})
\end{align*}
as $\univ{g}^{G'}$-modules.
Therefore, we have shown the following criterion.

\begin{Theorem}\label{thm:InjectivePhi}
	$\Phi_{M'}$ is injective if and only if \smash{$\Dzuck{\Delta(G')}{\Delta(L')}{S - 1}(M'\otimes \Coker(\Phi))^{\Delta(G')} = 0$}.
\end{Theorem}

\begin{Remark}\label{rmk:InjectivePhiStandard}
	If $(M'\otimes \Coker(\Phi))|_{\Delta(\lie{g'})}$ has an exhaustive filtration whose associated graded module has a standard filtration, we have \smash{$\Dzuck{\Delta(G')}{\Delta(L')}{S - 1}(M'\otimes \Coker(\Phi))^{\Delta(G')} = 0$} by Proposition~\ref{prop:StandardFiltrationVanishing}.
\end{Remark}

By Proposition~\ref{prop:CompletelyReducibleVerma}, \smash{$\bigl(\ind{g'}{q'}\bigl(\fdual{F'}\bigr)\otimes \Coker(\Phi)\bigr)|_{\Delta(\lie{g'})}$} has an exhaustive filtration whose associated graded module has a standard filtration.
Hence if \smash{$M' = \ind{g'}{q'}\bigl(\fdual{F'}\bigr)$}, then $\Phi_{M'}$ is injective.
In~general, the $(\lie{g'\oplus g}, \Delta(G'))$-module structure on \smash{$\Dzuck{\Delta(G')}{\Delta(L')}{S}(M'\otimes M)$} is not easy even if \smash{$M' = \ind{g'}{q'}\bigl(\fdual{F'}\bigr)$}.
Using Theorem~\ref{thm:Embedding2}, we obtain a partial result.

\begin{Corollary}\label{cor:BoundLength2}
	Let $F'$ be a finite-dimensional irreducible $\lie{l'}$-module with infinitesimal character~$[\lambda']$ satisfying
	\begin{align*}
		\frac{2(-\lambda'-\rhog{\lie{u}}, \alpha)}{(\alpha, \alpha)} \not \in \set{0, 1, 2, \ldots}\qquad \text{for any }\alpha \in \Delta(\lie{u'}, \lie{t'}).
	\end{align*}
	Then the length of $\Hom_{\lie{g'}}\bigl(\ind{g'}{q'}(F'), M\bigr)$ is less than or equal to that of $M$.
	In particular, if $M$ is irreducible, then $\Hom_{\lie{g'}}\bigl(\ind{g'}{q'}(F'), M\bigr)$ is irreducible or zero.
\end{Corollary}

\begin{proof}
	By Theorem~\ref{thm:Embedding2}, the length of \smash{$\Dzuck{\Delta(G')}{\Delta(L')}{S}\bigl(\ind{g'}{q'}(F'^d\bigr)\otimes M)$} is less than or equal to that of $M$.
	Hence the assertion follows from Theorem~\ref{thm:InjectivePhi}.
\end{proof}

\begin{Remark}
	Under the assumption of Corollary~\ref{cor:BoundLength2}, \smash{$\ind{g'}{q'}(F')$} is projective in \smash{$\BGGcat{g'}{q'}$}.
\end{Remark}

We shall estimate the length of the $\univ{g}^{G'}$-module \smash{$\Hom_{\lie{g'}}\bigl(\ind{g'}{q'}(F'), M\bigr)$} for generalized Verma modules $M$ with a dominance condition.
Let $F$ be a finite-dimensional irreducible $\lie{l}$-module and assume that $F$ satisfies the assumption of Theorem~\ref{thm:Exactness}, that is,
\begin{align*}
	\frac{2(\lambda+\rhog{\lie{u}}, \alpha)}{(\alpha, \alpha)} \not \in \set{0, 1, 2, \ldots}\qquad \text{for any }\alpha \in \Delta(\lie{u}, \lie{t}),
\end{align*}
where $[\lambda] \in \lie{t}^*/W_L$ is the infinitesimal character of $F$.
Then the submodule lattice $\smash{\Dzuck{\Delta(G')}{\Delta(L')}{S}}(M' \otimes \ind{g}{q}(F))$ is isomorphic to that of $M'$
by Theorem~\ref{thm:Exactness}.

\begin{Corollary}\label{cor:BoundLength}
	The length of the $\univ{g}^{G'}$-module \smash{$\Hom_{\lie{g'}}\bigl(\ind{g'}{q'}(F'), \ind{g}{q}(F)\bigr)$} is less then or equal to that of \smash{$\ind{g'}{q'}\bigl(\fdual{F'}\bigr)$}.
	In particular, if \smash{$\ind{g'}{q'}\bigl(\fdual{F'}\bigr)$} is irreducible, then \smash{$\Hom_{\lie{g'}}\bigl(\ind{g'}{q'}(F'), \ind{g}{q}(F)\bigr)$} is irreducible or zero.
\end{Corollary}

\begin{proof}
	The assertion follows from Theorems~\ref{thm:Exactness} and~\ref{thm:InjectivePhi} by putting
\[
M' = \ind{g'}{q'}(F')\qquad \text{and}\qquad {M = \ind{g}{q}(F)}.
\tag*{\qed}
\]
\renewcommand{\qed}{}
\end{proof}

\begin{Corollary}\label{cor:GeneralizedVermaIrreducible}
	Assume that $\Im(\Phi)$ for $M = \ind{g}{q}(F)$ is a direct summand of $\ind{g}{q}(F)|_{\lie{g'}}$.
	Then the $\univ{g}^{G'}$-module $\Hom_{\lie{g'}}\bigl(\ind{g'}{q'}(F'), \ind{g}{q}(F)\bigr)$ is irreducible or zero.
\end{Corollary}

\begin{Remark}
	For example, if $\ind{g}{q}(F)|_{\lie{g'}}$ is completely reducible, the assumption of the corollary holds.
\end{Remark}

\begin{proof}[Proof of Corollary~\ref{cor:GeneralizedVermaIrreducible}]
	Let $[\lambda']$ denote the infinitesimal character of $\ind{g'}{q'}(F')$.
	By Proposition~\ref{prop:AdmissibleDual}, the primary component $P_{[\lambda']}(M|_{\lie{g'}})$ has finite length.
	Since $M|_{\lie{g'}}$ has a standard filtration, so does $P_{[\lambda']}(M|_{\lie{g'}})$.
	By~\cite[Proposition~3.7 and Theorem~9.8\,(3)]{Hu08_category_o}, any direct summand of $P_{[\lambda']}(M|_{\lie{g'}})$ has a standard filtration.
	In particular, $\Coker(\Phi)$ has a standard filtration by assumption.
	Therefore, the assertion follows from Theorems~\ref{thm:Exactness} and~\ref{thm:InjectivePhi} by letting $M'$ be an~irreducible quotient of \smash{$\ind{g'}{q'}\bigl(\fdual{F'}\bigr)$}.
	See Proposition~\ref{prop:CompletelyReducibleVerma} and Remark~\ref{rmk:InjectivePhiStandard} for the vanishing of~\smash{$\Dzuck{\Delta(G')}{\Delta(L')}{S - 1}(M'\otimes \Coker(\Phi))^{\Delta(G')}$}.
\end{proof}

\subsection{Quasi-abelian parabolic subalgebra}

In this subsection, we give a sufficient condition for completely reducibility of $\ind{g}{q}(F)|_{\lie{g'}}$
to apply Corollary~\ref{cor:GeneralizedVermaIrreducible}.
Retain the notation in the previous subsection.

Set $\lie{u}''\coloneq \lie{u} \cap (\lie{g}')^\perp$ and $\lie{\overline{u}}''\coloneq \lie{\overline{u}} \cap (\lie{g}')^\perp$.
Then $\lie{u}''$ is $\lie{u}'$-stable and $\lie{u} = \lie{u}' \oplus \lie{u}''$ holds by the definition of weakly $\lie{g'}$-compatible parabolic subalgebras (Definition~\ref{def:WeaklyCompatible}).
The notion of quasi-abelian parabolic subalgebras is defined in~\cite[p.~109]{EnPaWaWo85} for symmetric $(\lie{g}, \lie{g'})$.
We consider a~straightforward generalization of the definition.

\begin{Definition}\label{def:QuasiAbelian}
	$\lie{q}$ is said to be \emph{quasi-abelian} with respect to $\lie{g'}$ if
	$(\alpha, \beta)\geq 0$ holds for any $\alpha \in \Delta(\lie{u'}, \lie{t'})$ and $\beta \in \Delta(\lie{u''}, \lie{t'})$.
\end{Definition}

It is clear that if $\lie{p}$ is a weakly $\lie{g'}$-compatible parabolic subalgebra containing $\lie{q}$ and $\lie{q}$ is quasi-abelian, then $\lie{p}$ is also quasi-abelian.
In fact, the nilpotent radical of $\lie{p}$ is smaller than that of~$\lie{q}$.

If the nilpotent radical of $\lie{q}$ is abelian, $\lie{q}$ is quasi-abelian.
More generally, we have the following criterion.

\begin{Proposition}\label{prop:QuasiAbelianCondition}
	If $[\lie{u}', \lie{u}'']=0$ holds, then $\lie{q}$ is quasi-abelian.
	In particular, a parabolic subalgebra with abelian nilpotent radical is quasi-abelian.
\end{Proposition}
	
\begin{proof}
	Let $\alpha \in \Delta(\lie{u'}, \lie{t'})$ and $\beta \in \Delta(\lie{u''}, \lie{t'})$.
	Assume $(\alpha, \beta) < 0$.
	This implies $\alpha + \beta \in \Delta(\lie{u''}, \lie{t'})$
	and \smash{$\bigl[\lie{u}'_{\alpha}, \lie{u}''_{\beta}\bigr]\neq 0$}.
	By assumption, \smash{$\bigl[\lie{u}'_{\alpha}, \lie{u}''_{\beta}\bigr]\subset [\lie{u}', \lie{u}'']=0$} holds and this is contradiction.
	Therefore, $\lie{q}$ is quasi-abelian.
\end{proof}

\begin{Example}
	If $\lie{g'} \simeq \lie{sl}(2,\CC)$, then $\lie{q}$ is quasi-abelian as follows.
	$\lie{q}$ contains $\lie{q}(H)$ defined by a semisimple element $H\in \lie{g'}$.
	Since any eigenvalue in $\lie{u}(H)$ of $\ad(H)$ is positive, $\lie{q}(H)$ is quasi-abelian.
	As noted immediately after Definition~\ref{def:QuasiAbelian}, $\lie{q}$ is also quasi-abelian.

	Take an $\lie{sl}_2$-triple $\set{H, X, Y}\subset \lie{g'}$ with $[X, Y] = H$.
	If $\ad(H)$ has an eigenvalue greater than or equal to $3$ in $\lie{g}$, then we have $[\lie{u'}, \lie{u}''] \neq 0$
	by the representation theory of $\lie{sl}(2,\CC)$.
	Take a~finite-dimensional irreducible $\lie{sl}(2,\CC)$-module $F$ of dimension greater than or equal to $3$.
	Then the embedding $\lie{g'} = \lie{sl}(2,\CC) \hookrightarrow \lie{sl}(F) = \lie{g}$ gives an example
	of quasi-abelian parabolic subalgebra with $[\lie{u'}, \lie{u''}] \neq 0$.
	This means that the converse of Proposition~\ref{prop:QuasiAbelianCondition} does not hold in general.
\end{Example}

A criterion for the completely reducibility of $\ind{g}{q}(F)|_{\lie{g'}}$ is given in~\cite[Lemma~3.1]{EnPaWaWo85} if $(\lie{g}, \lie{g'})$ is a symmetric pair.
We extend the criterion to our setting.

\begin{Lemma}\label{lem:QuasiAbelianCompletelyReducible}
	Let $F$ be a finite-dimensional irreducible $\lie{l}$-module.
	Assume that for any irreducible submodule of $F|_{\lie{l}'}$, its infinitesimal character $[\lambda'] \in (\lie{t'})^*/W_{L'}$ satisfies
	\begin{align*}
		\frac{2(\lambda'+\rho(\lie{u}'), \alpha)}{(\alpha, \alpha)} \not \in \set{1, 2, \ldots}, \qquad \forall \alpha \in \Delta(\lie{u'}, \lie{t'}).
	\end{align*}
	Suppose that $\lie{q}$ is quasi-abelian.
	Then \smash{$\ind{g}{q}(F)|_{\lie{g}'}$} is completely reducible and each irreducible component is isomorphic to an irreducible generalized Verma module.
\end{Lemma}

\begin{proof}
	By Proposition~\ref{prop:FiltrationRestriction}, \smash{$\ind{g}{q}(F)|_{\lie{g}'}$} has a standard filtration with
	\begin{align*}
		\gr\bigl(\ind{g}{q}(F)|_{\lie{g}'}\bigr) \simeq \ind{g'}{q'}\bigl(F\otimes S\bigl(\lie{\overline{u}}''\bigr)\bigr),
	\end{align*}
	where the $\lie{u}'$-action on $S(\lie{\overline{u}}'')$ is trivial.

	Let $V$ be an irreducible $\lie{l'}$-submodule of $F\otimes S(\lie{\overline{u}}'')$.
	By the Weyl character formula, the infinitesimal character of $V$ is of the form $[\lambda' - R]$,
	where $R$ is a sum of elements in $\Delta(\lie{u''}, \lie{t'})$ and~$[\lambda']$ is the infinitesimal character of some irreducible submodule of $F|_{\lie{l'}}$.
	Since $\lie{q}$ is quasi-abelian, we have
	\begin{align*}
		\frac{2(R, \alpha)}{(\alpha, \alpha)} \in \set{0, 1, 2, \ldots}, \qquad \forall \alpha \in \Delta(\lie{u'}, \lie{t'}).
	\end{align*}
	This and the assumption imply
	\begin{align*}
		\frac{2(\lambda'-R+\rho(\lie{u}'), \alpha)}{(\alpha, \alpha)} \not \in \set{1, 2, \ldots}\qquad \text{for any }\alpha \in \Delta(\lie{u'}, \lie{t'}).
	\end{align*}
	Hence \smash{$\ind{g'}{q'}(V)$} is irreducible by Fact~\ref{fact:VermaIrreducible}.
	Therefore, Proposition~\ref{prop:CompletelyReducibleVerma} shows the assertion.
\end{proof}

The condition in Lemma~\ref{lem:QuasiAbelianCompletelyReducible} is given in terms of $\lie{t'}$-weights.
We shall give a criterion for the completely reducibility in terms of $\lie{t}$-weights.
To do so, we prepare several convexity results.
Note that J.A.~Vargas gave a similar estimate in~\cite[\S(1.5)]{Va16} in the context of the branching problem of discrete series representations.
Our result (Theorem~\ref{thm:QuasiAbelianGoodRange}) contains his estimate.
See also Theorems~\ref{thm:DerivedFunctorModuleQuasiAbelian} and~\ref{thm:DiscreteSeries}.

Suppose $\lie{t'}\subset \lie{t}$.
Let $\lie{t}''$ be the orthogonal complement of $\lie{t}'$ in $\lie{t}$.
Consider $(\lie{t'})^*$ as a subspace of $\lie{t}^*$
using the direct sum decomposition $\lie{t}=\lie{t'} \oplus \lie{t''}$.
For $\alpha \in \Delta(\lie{u'}, \lie{t'})$, define
\begin{align*}
	\Delta(\alpha)\coloneq \set{\beta \in \Delta(\lie{u}, \lie{t})\mid \beta|_{\lie{t}'}=\alpha}.
\end{align*}
For a subset $S$ of a real vector space, we denote by $\Co(S)$ the convex hull of $S$.

\begin{Lemma}\label{lem:ConvexityRestriction}
	Let $\alpha \in \Delta(\lie{u'}, \lie{t'})$.
	Then $\alpha \in \Co(\Delta(\alpha))$ holds.
\end{Lemma}

\begin{proof}
	Fix a Cartan involution $\theta$ of $\lie{g}$ such that
	$\lie{g}', \lie{t}$ and $\lie{t}'$ are $\theta$-stable.
	In fact, since $\lie{g}'$ is reductive in $\lie{g}$, such an involution exists.

	Take a root vector $X \in \lie{u}'_{\alpha}$.
	By $X \in \lie{u}$, we can write \smash{$X=\sum_{\beta \in \Delta(\alpha)}X_{\beta}$} with $X_{\beta} \in \lie{u}_{\beta}$.
	Consider the inner product $\langle \cdot, \cdot \rangle \coloneq -(\cdot, \theta(\cdot))$ on $\lie{g}$.
	Note that $[\lie{t''}, \lie{g'}] \subset (\lie{g}')^\perp$.
	Then the root spaces are mutually orthogonal, and hence we have
	\begin{align*}
		&\|X\|^2= \sum_{\beta \in \Delta(\alpha)} \|X_\beta\|^2, \\
		&\sum_{\beta \in \Delta(\alpha)} \langle \beta(Z')X_\beta, X_\beta\rangle= \langle \alpha(Z')X, X\rangle, \\
		&\sum_{\beta \in \Delta(\alpha)} \langle \beta(Z'')X_\beta, X_\beta\rangle= \langle [Z'', X], X\rangle = 0 = \langle \alpha(Z'')X, X\rangle
	\end{align*}
	for any $Z' \in \lie{t'}$ and $Z'' \in \lie{t''}$.
	This shows
	\begin{align*}
		\alpha = \sum_{\beta \in \Delta(\alpha)} \frac{\|X_\beta\|^2}{\|X\|^2} \beta, \qquad \sum_{\beta \in \Delta(\alpha)} \frac{\|X_\beta\|^2}{\|X\|^2} = 1.
	\end{align*}
	Therefore, we obtain $\alpha \in \Co(\Delta(\alpha))$.
\end{proof}

\begin{Lemma}\label{lem:SumConvex}
	Fix a Borel subalgebra $\lie{b}_L = \lie{t}\oplus \lie{u}_L$ of $\lie{l}$.
	Let $\lambda_1$ and $\lambda_2$ be dominant integral weights of $\lie{l}$.
	Then we have $\Co(W_{L}\lambda_1)+\Co(W_{L}\lambda_2) = \Co(W_{L}(\lambda_1+\lambda_2))$.
\end{Lemma}
	
\begin{proof}
	$\Co(W_{L}\lambda_1)+\Co(W_{L}\lambda_2) \supset \Co(W_{L}(\lambda_1+\lambda_2))$ is obvious.
	We shall show the converse inclusion.
	Let $s \in W_L$ and it is enough to prove $\lambda_1 +s(\lambda_2) \in \Co(W_{L}(\lambda_1+\lambda_2))$.

	For a dominant integral weight $\lambda$ of $\lie{l}$, we denote by $F(\lambda)$ the irreducible $\lie{l}$-module with the highest weight $\lambda$.
	Then there exists a unique non-zero homomorphism $p\colon F(\lambda_1) \otimes F(\lambda_2) \rightarrow F(\lambda_1 + \lambda_2)$ up to scalar.
	For each $i = 1, 2$, fix a highest weight vector $v_i$ of $F(\lambda_i)$.
	Then we have $p(v_1\otimes v_2) \neq 0$.
	
	Take a weight vector $v'_2 \in F(\lambda_2)$ with the extreme weight $s(\lambda_2)$.
	Then there exists $X \in \univ{u_{\mathit{L}}}$ such that $Xv'_2=v_2$.
	Since $v_1$ is $\lie{u}_L$-invariant, we have $X(v_1\otimes v'_2)=v_1 \otimes v_2$.
	This implies
	\begin{align*}
		Xp\bigl(v_1\otimes v'_2\bigr) = p\bigl(X\bigl(v_1\otimes v'_2\bigr)\bigr) = p(v_1\otimes v_2)\neq 0,
	\end{align*}
	and hence $p\bigl(v_1\otimes v'_2\bigr) \neq 0$.
	
	Since any weight of $F(\lambda_1+\lambda_2)$ belongs to $\Co(W_{L}(\lambda_1+\lambda_2))$,
	this shows $\lambda_1+s(\lambda_2) \in \Co(W_{L}(\lambda_1+\lambda_2))$.
	We have proved the lemma.
\end{proof}

Fix a basis $\set{v_1, \ldots, v_{r'}}$ of $\lie{t'}$ and extend it to a basis $\set{v_1, \ldots, v_r}$ of $\lie{t}$.
The bases determine lexicographical orders on $(\lie{t'})^*$ and $\lie{t}^*$.
Let $\Delta^+(\lie{l'}, \lie{t'})$ and $\Delta^+(\lie{l}, \lie{t})$ denote the sets of positive roots given by the orders.
Write $\rhog{\lie{l'}}$ (resp.\ $\rhog{\lie{l}}$) for half the sum of roots in $\Delta^+(\lie{l'}, \lie{t'})$ (resp.\ $\Delta^+(\lie{l}, \lie{t})$).
Then $\rhog{\lie{l}}|_{\lie{t'}}$ is a dominant integral weight of $\lie{l}'$.

\begin{Theorem}\label{thm:QuasiAbelianGoodRange}
	Let $F$ be a finite-dimensional irreducible $\lie{l}$-module with infinitesimal character $[\lambda] \in \lie{t}^*/W_L$
	in the good range, namely,
	\begin{align*}
		\Re(\lambda+\rho(\lie{u}), \alpha) < 0, \qquad \forall \alpha \in \Delta(\lie{u}, \lie{t}).
	\end{align*}
	Suppose that $\lie{q}$ is quasi-abelian.
	Then \smash{$\ind{g}{q}(F)|_{\lie{g}'}$} is completely reducible,
	and each irreducible direct summand is of the form \smash{$\ind{g'}{q'}(F')$} such that $F'$ is a finite-dimensional irreducible $\lie{l'}$-module in the good range.
\end{Theorem}

\begin{proof}
	Let $F'$ be an irreducible submodule of $F|_{\lie{l}'}$ and $[\lambda']$ the infinitesimal character of $F'$.
	We shall prove
	$\Re(\lambda'+\rhog{\lie{u}'}, \alpha) < 0$, $\forall \alpha \in \Delta(\lie{u'}, \lie{t'})$.
	If we prove this, the assertion follows by the same way as Lemma~\ref{lem:QuasiAbelianCompletelyReducible}.

	Replacing the representatives $\lambda$ and $\lambda'$, we may assume that $\lambda$ and $\lambda'$ are dominant with respect to $\Delta^+(\lie{l}, \lie{t})$ and $\Delta^+(\lie{l'}, \lie{t'})$, respectively.
	By Lemma~\ref{lem:SumConvex}, we have
	\begin{align*}
		\lambda' &\in \Co(W_L(\lambda - \rhog{\lie{l}}))|_{\lie{t'}} + \rhog{\lie{l'}}
		 \subset \Co(W_{L}(\lambda-\rhog{\lie{l}}))|_{\lie{t'}}+\Co(W_{L}\rhog{\lie{l}})|_{\lie{t'}} - \rhog{\lie{l}}|_{\lie{t'}} + \rhog{\lie{l'}} \\
		&= \Co(W_{L}(\lambda))|_{\lie{t'}} - \rhog{\lie{l}}|_{\lie{t'}} + \rhog{\lie{l'}}.
	\end{align*}
	From this, we write
	\begin{align*}
		\lambda' = \sum_{s \in W_L} a_s s(\lambda)|_{\lie{t'}} - \rhog{\lie{l}}|_{\lie{t'}} + \rhog{\lie{l'}}
	\end{align*}
	with $\sum_{s\in W_L}a_s = 1$ and $a_s \geq 0$.

	Let $\alpha \in \Delta(\lie{u'}, \lie{t'})$.
	By Lemma~\ref{lem:ConvexityRestriction}, $\alpha$ is written as
	\begin{align*}
		\alpha = \sum_{\beta \in \Delta(\alpha)} c_{\beta} \beta
	\end{align*}
	with $\sum_{\beta \in \Delta(\alpha)}c_{\beta} = 1$ and $c_{\beta} \geq 0$.
	Then we have
	\begin{align}
		&(\lambda'+\rhog{\lie{u}'}, \alpha) + (\rhog{\lie{u}}|_{\lie{t'}} + \rhog{\lie{l}}|_{\lie{t'}} - \rhog{\lie{l'}} - \rhog{\lie{u'}}, \alpha) \nonumber \\
		&\qquad{}=\biggl(\sum_{s \in W_L} a_s s(\lambda) + \rhog{\lie{u}}, \alpha\biggr) \nonumber \\
		&\qquad{}=\sum_{s \in W_L, \beta \in \Delta(\alpha)}a_s c_{\beta}\bigl(\lambda + \rhog{\lie{u}}, s^{-1}(\beta)\bigr). \label{eqn:PositivityTo}
	\end{align}
	By $\Delta(\alpha)\subset \Delta(\lie{u}, \lie{t})$ and the assumption, the real part of \eqref{eqn:PositivityTo} is negative.
	By Proposition~\ref{prop:PositivityRho}, we have
	\begin{align*}
		(\rhog{\lie{u}}|_{\lie{t'}} + \rhog{\lie{l}}|_{\lie{t'}} - \rhog{\lie{l'}} - \rhog{\lie{u'}}, \alpha) \geq 0.
	\end{align*}
	Therefore, we obtain $(\lambda'+\rhog{\lie{u}'}, \alpha) < 0$.
\end{proof}

Combining Corollary~\ref{cor:GeneralizedVermaIrreducible} with Theorem~\ref{thm:QuasiAbelianGoodRange}, we obtain the following.

\begin{Corollary}\label{cor:IrreducibleVermaQuasiAbelian}
	Let $F$ be a finite-dimensional irreducible $\lie{l}$-module with infinitesimal character $[\lambda] \in \lie{t}^*/W_L$
	in the good range, and $F'$ a finite-dimensional irreducible $\lie{l}'$-module.
	Suppose that~$\lie{q}$ is quasi-abelian.
	Then the $\univ{g}^{G'}$-module \smash{$\Hom_{\lie{g'}}\bigl(\ind{g'}{q'}(F'), \ind{g}{q}(F)\bigr)$} is irreducible or zero.
\end{Corollary}

The estimate in the proof of Theorem~\ref{thm:QuasiAbelianGoodRange} is independent of the assumption that $\lie{q}$ is quasi-abelian.

\begin{Corollary}
	Let $F$ be a finite-dimensional irreducible $\lie{l}$-module in the good range.
	Let $F'$ be an irreducible submodule of $F|_{\lie{l}'}$ and $[\lambda']$ the infinitesimal character of $F'$.
	Then $\lambda'$ is in the good range, i.e.,
	\begin{align*}
		\Re(\lambda'+\rhog{\lie{u}'}, \alpha) < 0, \qquad \forall \alpha \in \Delta(\lie{u'}, \lie{t'}).
	\end{align*}
\end{Corollary}

\subsection{Holomorphic discrete series representation}

A typical example of Corollary~\ref{cor:IrreducibleVermaQuasiAbelian} is the case of holomorphic discrete series representations.
In the case, the branching laws are studied in~\cite{JaVe79} (at the unitary representation level),
and in~\mbox{\cite{Ko98,Ko08}}.
Theorem~\ref{thm:QuasiAbelianGoodRange} has been proved in~\cite[Theorem~7.4]{Ko98}.
Corollary~\ref{cor:IrreducibleVermaQuasiAbelian} for holomorphic discrete series representations of some classical groups is reduced to the compact group case by the Howe duality and see-saw pair (see, e.g.,~\cite[Lemma~2.6]{Ma16_dual_pair}).

Let $(\lie{g}, G')$ be a pair with semisimple $\lie{g}$,
and $\theta$ an involution of $\lie{g}$ fixing $\lie{g'}$.
Set $\lie{k}\coloneq \lie{g}^\theta$ and~$\lie{k'}\coloneq (\lie{g'})^\theta$.
Suppose that there exists a semisimple element $H \in \lie{g}$ such that
$H \in \lie{k'}$, $\lie{g}^H = \lie{k}$ and~$\ad(H)$ has eigenvalues $-1$, $0$, $1$ in $\lie{g}$.
Define parabolic subalgebras $\lie{q} = \lie{q}(H) = \lie{k} \oplus \lie{u}(H)$
and $\lie{q'} \coloneq \lie{k}'\oplus (\lie{u}(H)\cap \lie{g'})$ (see Section~\ref{subsect:GeneralizedVerme} for the notation).
Then $\lie{q}$ is quasi-abelian with respect to~$\lie{g'}$ by Proposition~\ref{prop:QuasiAbelianCondition}.

Fix a real form $\lie{g}_\RR$ of $\lie{g}$ such that $\theta|_{\lie{g}_\RR}$ is a Cartan involution.
Set $\lie{k}_\RR\coloneq \lie{g}_\RR \cap \lie{k}$, $\lie{g}'_\RR\coloneq \lie{g}'\cap \lie{g}_\RR$ and~${\lie{k}'_\RR\coloneq \lie{g}'_\RR \cap \lie{k}}$.
Let $G_\RR$ be a simply-connected connected semisimple Lie group with the Lie algebra~$\lie{g}_\RR$.
Let $K_\RR$, $G'_\RR$ and $K'_\RR$ denote the analytic subgroups of $G_\RR$ with the Lie algebras~$\lie{k}_\RR$,~$\lie{g}'_\RR$ and~$\lie{k}'_\RR$, respectively.

Let $F$ be an irreducible unitary representation of $K_\RR$.
Then $F$ is finite-dimensional.
By Harish-Chandra's classification of holomorphic discrete series representations, $F$ is in the good range if and only if
$\ind{g}{q}(F)$ is unitarizable and isomorphic to the underlying Harish-Chandra module of a holomorphic discrete series representation of $G_\RR$ (with respect to $\lie{q}$).
For simplicity, irreducible finite-dimensional unitary representations are regarded as holomorphic discrete series representations.
We refer the reader to~\cite[Chapter~VI]{Kn86} for holomorphic discrete series representations.

In this setting, Theorem~\ref{thm:QuasiAbelianGoodRange} and Corollary~\ref{cor:IrreducibleVermaQuasiAbelian} are rephrased as follows.
Although the following theorem has an overlap with Theorem~\ref{thm:DiscreteSeries} (the discrete series case),
we state it explicitly because in the present setting $K_\RR$ is non-compact and $(\lie{g}, \lie{g'})$ need not be a symmetric pair.

\begin{Theorem}\label{thm:IrreducibleHolomorphicDisc}
	Suppose $V$ $($resp.\ $W)$ is a holomorphic discrete series representation of $G_\RR$ $($resp.~$G'_\RR)$.
	Then $V|_{G'_\RR}$ is discretely decomposable and a direct sum of holomorphic discrete series representations of $G'_\RR$.
	Moreover, \smash{$\univ{g}^{G'}$} acts on \smash{$\Hom_{G'_\RR}\bigl(W, V|_{G'_\RR}\bigr)$} irreducibly if the space is non-zero.
\end{Theorem}

\begin{Remark}
	The discrete decomposability (or, more strongly, $G'_\RR$-admissibility) is well-known in~\cite[Theorem~7.4]{Ko98}.
\end{Remark}

\begin{proof}[Proof of Theorem~\ref{thm:IrreducibleHolomorphicDisc}]
	We denote by $V_{K_\RR}$ the subspace of all $K_\RR$-finite vectors in $V$.
	Note that the center~$Z$ of~$G_\RR$ acts on $V$ via a character and $K_\RR/Z$ is compact.
	This implies that $V_{K_\RR}$ behaves as a~underlying Harish-Chandra module for a~semisimple Lie group with finite center.

	Hence $V_{K_\RR}$ is a $\lie{g}$-submodule in $V^\infty$ the space of smooth vectors.
	Then $V_{K_\RR}$ is isomorphic to a generalized Verma module \smash{$\ind{g}{q}(F)$} with finite-dimensional irreducible $\lie{k}$-module $F$ in the good range.
	Similarly, $W_{K'_\RR}$ is defined and isomorphic to a generalized Verma module \smash{$\ind{g'}{q'}(F')$} with finite-dimensional irreducible $\lie{k}$-module $F'$ in the good range.

	Remark that \smash{$\ind{g}{q}(F)$} is $\lie{k'}$-admissible by Proposition~\ref{prop:AdmissibleDual}.
	This implies \smash{$V_{K_\RR} = V_{K'_\RR}$}.
	Hence the irreducible decomposition of \smash{$V_{K_\RR}|_{\lie{g'}}$} gives the irreducible decomposition of $V|_{G'_\RR}$
	by taking completion.
	In particular, the restriction map gives an isomorphism
	\begin{align*}
		\Hom_{G'_\RR}\bigl(W, V|_{G'_\RR}\bigr) \xrightarrow{\simeq} \Hom_{\lie{g'}}\bigl(W_{K'_\RR}, V_{K_\RR}|_{\lie{g'}}\bigr)
	\end{align*}
	of vector spaces.
	This shows the assertion.
	See~\cite[Proposition~1.6]{Ko98_discrete_decomposable_3} for this isomorphism.
\end{proof}

In general (including non-holomorphic cases), $\univ{g}^{G'}$ may not act on the Hilbert space \smash{$\Hom_{G'_\RR}\bigl(W, V|_{G'_\RR}\bigr)$} as an algebra.
In fact, an element in $\univ{g}^{G'}$ acts on the space as a densely defined closable operator, and the action is defined without using \smash{$V_{K_\RR}$} and \smash{$W_{K'_\RR}$}.
See, e.g.,~\cite[Theorem~10.24]{Ki20}.
In the current setting, \smash{$\Hom_{G'_\RR}\bigl(W, V|_{G'_\RR}\bigr)$} is finite-dimensional, and hence $\univ{g}^{G'}$ actually acts on \smash{$\Hom_{G'_\RR}\bigl(W, V|_{G'_\RR}\bigr)$}.

Remark that Theorem~\ref{thm:IrreducibleHolomorphicDisc} does not hold for unitary highest weight modules.
We shall give an example in which the $\univ{g}^{G'}$-module is not irreducible.

\begin{Example}\label{ex:NonIrreducible}
	Suppose that $\lie{g} = \lieSp(2n,\CC)$ and $\lie{g'} = \lieSp(n, \CC)$.
	Here $\lie{g'}$ is the diagonal embedding in $\lieSp(n,\CC)\oplus \lieSp(n, \CC) \subset \lieSp(2n,\CC)$.
	Take a real form $\lie{g}_\RR = \lieSp(2n,\RR)$ of $\lie{g}$ such that~${\lie{g}'_\RR \coloneq \lie{g'}\cap \lie{g}_\RR}$ is $\lieSp(n,\RR)$.
	Let $(\omega, V)$ be the even part of the Segal--Shale--Weil representation of~$\lieSp(2n,\RR)$.
	Then $V|_{\lie{g'}}$ is completely reducible.
	Take an irreducible $\lie{g'}$-submodule $V'$ of~$V|_{\lie{g'}}$.

	By the classical invariant theory, $\omega(\univ{g})^{G'}$ is commutative.
	Hence, if $\Hom_{\lie{g'}}(V', V)$ is irreducible for any $V'$, $V|_{\lie{g'}}$ should be multiplicity-free.
	It is well-known that $V|_{\lie{g'}}$ is not multiplicity-free.
	In fact, $\mathrm{O}(2)$ that is the commutator of $\mathrm{Sp}(n,\RR)$ in $\mathrm{Sp}(2n,\RR)$ acts on $\Hom_{\lie{g'}}(V', V)$ irreducibly.
	This implies that $\Hom_{\lie{g'}}(V', V)$ is not irreducible excluding special $V'$.
	
	Using the theta lifting, we can obtain many concrete examples of $\univ{g}^{G'}$-modules.
	See~\cite{KaVe78} and~\cite{Ho89} for the branching law and the theta lifting.
\end{Example}

\subsection{Zuckerman derived functor module}

As an application of Corollary~\ref{cor:IrreducibleVermaQuasiAbelian},
we consider the branching problem of Zuckerman derived functor modules induced from
quasi-abelian parabolic subalgebras.

Let $(\lie{g}, K)$ be a pair with semisimple $\lie{g}$, and $(\lie{g'}, K')$ its subpair with reductive $\lie{g'}$ in $\lie{g}$.
Suppose that there exists an involution $\theta$ of $\lie{g}$ such that $\lie{k} = \lie{g}^\theta$ and $\lie{k'} = (\lie{g'})^\theta$,
and there exists a connected reductive algebraic group $G$ with the Lie algebra $\lie{g'}$ such that $(\lie{g}, G')$ is a pair.
Fix $H \in \lie{k}$ semisimple in $\lie{g}$ with real eigenvalues
and define a parabolic subalgebra $\lie{q}\coloneq \lie{q}(H)$.
Then $\lie{q}$ is $\theta$-stable.
Suppose that $\lie{q}$ is weakly $\lie{g'}$-compatible, and set $\lie{q'}\coloneq \lie{g'}\cap \lie{q}$.
We use the notation $\lie{u}, \lie{l}, \lie{u'}, \lie{l'}, \ldots$ as in Section~\ref{subsect:GeneralizedVerme}.

Write $L_K$ for the centralizer of $H$ in $K$ and set $L'_K\coloneq K'\cap L_K$.
Fix $\theta$-stable Cartan subalgebras $\lie{t'}$ and $\lie{t}$ of $\lie{l'}$ and $\lie{l}$, respectively, such that $\lie{t'}\subset \lie{t}$.
Let $\rhog{\lie{u}}$ denote half the sum of roots in $\Delta(\lie{u}, \lie{t})$.

In this subsection, we consider Zuckerman derived functor modules defined by
\begin{align*}
	\lmod{g}{q}{i}(V)\coloneq \Dzuck{K}{L_K}{i}\bigl(\ind{g}{q}(V)\bigr)
\end{align*}
for an $(\lie{l}, L_K)$-module $V$.
Then \smash{$\lmod{g}{q}{i}(V)$} is a $(\lie{g}, K)$-module.
In~\cite{KnVo95_cohomological_induction}, this module is defined by the Bernstein functor $\Pi$,
and our parametrization is so-called unnormalized version, which is written as $\!^{\mathrm{u}}\mathcal{L}$.

Set $S\coloneq \dim(\lie{u}\cap \lie{k})$.
For an $(\lie{l}, L_K)$-module $F$ with infinitesimal character $\lambda$, we say that $\lambda$ (or $V$) is in the good range if $\lambda$ satisfies $\Re(\lambda + \rhog{\lie{u}}, \alpha) < 0$ for any $\alpha \in \Delta(\lie{u}, \lie{t})$.
The following fact is a fundamental result about $\lmod{g}{q}{i}(V)$ (see~\cite[Theorem~0.50]{KnVo95_cohomological_induction}).

\begin{Fact}\label{fact:ZuckermanDerivedFunctorModule}
	Let $V$ be an irreducible $(\lie{l}, L_K)$-module with infinitesimal character $[\lambda]$.
	\begin{enumerate}\itemsep=0pt
		\item[$1.$] \smash{$\lmod{g}{q}{i}(V)$} has the infinitesimal character $[\lambda+\rhog{\lie{u}}]$.
		\item[$2.$] If $\lambda$ is in the good range, \smash{$\lmod{g}{q}{i}(V)$} is zero for $i\neq S$
		and non-zero irreducible for $i=S$.
	\end{enumerate}
\end{Fact}

\begin{Proposition}\label{prop:UniqueZuckermanModule}
	Let $F_i$, $i=1,2$, be a finite-dimensional irreducible $(\lie{l}, L_K)$-module with infinitesimal character $[\lambda_i]$ in the good range.
	If \smash{$\lmod{g}{q}{S}(F_1) \simeq \lmod{g}{q}{S}(F_2)$}, then $F_1 \simeq F_2$ holds.
\end{Proposition}

\begin{proof}
	Set $\lie{t}^*_\RR\coloneq \mathrm{span}_{\RR}\Delta(\lie{g}, \lie{t})$.
	Then the symmetric form $(\cdot, \cdot)$ is an inner product on $\lie{t}^*_\RR$
	and $\lie{t}^*_\RR$ is stable under the action of the Weyl group $W_G$ of $\lie{g}$.
	For $\mu \in \lie{t}^*$, we denote by $\Re(\mu)$ the real part of $\mu$ with respect to the real form $\lie{t}^*_\RR$.

	Fix a set $\Delta^+(\lie{l}, \lie{t})$ of positive roots of $\Delta(\lie{l}, \lie{t})$.
	We may assume that $\lambda_1$ and $\lambda_2$ are anti-dominant with respect to $\Delta^+(\lie{l}, \lie{t})$.
	By Fact~\ref{fact:ZuckermanDerivedFunctorModule}, $\lmod{g}{q}{S}(F_i)$, $i=1,2$, has the infinitesimal character~$[\lambda_i + \rhog{\lie{u}}]$.
	By assumption, $\Re(\lambda_i + \rhog{\lie{u}})$, $i = 1,2$, is regular and anti-dominant with respect to~${\Delta^+(\lie{l}, \lie{t}) \cup \Delta(\lie{u}, \lie{t})}$.

	Since $\lmod{g}{q}{S}(F_1) \simeq \lmod{g}{q}{S}(F_2)$, there exists $s \in W_G$ such that $s(\lambda_1 + \rhog{\lie{u}}) = \lambda_2 + \rhog{\lie{u}}$.
	By the anti-dominance, $s$ is identity and hence we obtain $\lambda_1 = \lambda_2$.
\end{proof}

In a special case, the functor \smash{$\Dzuck{M}{M \cap L_K}{i}$} can be computed by \smash{$\Dzuck{K'}{K' \cap L_K}{i}$}.
It is a generalization of~\cite[Lemma~7]{GrWa00}, which is for discrete series representations.

\begin{Lemma}\label{lem:ZuckermanKtoH}
	Let $V$ be an $(\lie{l}, L_K)$-module and $M$ a connected reductive subgroup of $K$ acting on $K/L_K$ transitively.
	Fix $i \in \NN$.
	Then there exists an isomorphism
	\begin{align*}
		\lmod{g}{q}{i}(V) \simeq \Dzuck{M}{M \cap L_K}{i}\bigl(\ind{g}{q}(V)\bigr)
	\end{align*}
	of $(\lie{g}, M)$-modules.
\end{Lemma}

\begin{proof}
	By the construction of \smash{$\Dzuck{K}{L_K}{i}$} (see before Proposition~\ref{prop:KInvZuckerman}), we have
	\begin{align*}
		\lmod{g}{q}{i}(V) \simeq H^i\bigl(\lie{k}, L_K; \ind{g}{q}(V)\otimes \rring{K}\bigr).
	\end{align*}
	Recall that $H^i\bigl(\lie{k}, L_K; \ind{g}{q}(V)\otimes \rring{K}\bigr)$ is the cohomology of the complex given by
	\begin{align*}
		C^j\bigl(\lie{k}, L_K; \ind{g}{q}(V)\otimes \rring{K}\bigr) = \Hom_{L_K}\bigl(\wedge^j (\lie{k}/\lie{l}_K),\ind{g}{q}(V)\otimes \rring{K}\bigr).
	\end{align*}
	See~\cite[Section~I.8]{BoWa00_continuous_cohomology}.
	Since $M$ acts on $K/L_K$ transitively, the restriction map induces a bijection
	\begin{align*}
		(\rring{K}\otimes W)^{L_K} \xrightarrow{\simeq} (\rring{M} \otimes W)^{M\cap L_K}
	\end{align*}
	for any $L_K$-module $W$.
	Note that $M\cap L_K$ is reductive since $M/(M\cap L_K) \simeq K/L_K$ is affine.
	Hence we have an isomorphism
	\begin{align*}
		&\Hom_{L_K}(\wedge^j \bigl(\lie{k}/\lie{l}_K),\ind{g}{q}(V)\otimes \rring{K}\bigr) \\
		&\qquad{}\simeq \Hom_{M\cap L_K}\bigl(\wedge^j \bigl(\lie{m}/\lie{m}\cap \lie{l}_K\bigr), \ind{g}{q}(V)\otimes \rring{M}\bigr)
	\end{align*}
	of vector spaces.
	It is easy to see that this isomorphism induces an isomorphism of complexes.
	Taking the cohomology, we obtain the lemma.
\end{proof}

\begin{Example}
	If $K$ is simple, there are few tuples $(K, M, L_K)$ such that $M$ acts on $K/L_K$ transitively,
	\begin{align*}
		(K, M, L_K)={}& (\lieSL(2n,\CC), \lieSp(n,\CC), \lieGL(2n-1, \CC)), \\
		&(\lieSO(2n,\CC), \lieSO(2n-1,\CC), \lieGL(n,\CC)), \\
		&\bigl(\lieSO(7,\CC), \lieG, \lieSO(5,\CC)\times \CC^\times\bigr)
	\end{align*}
	are examples.
	See~\cite[Section 3 and Example 3.7]{Ko11}.
\end{Example}

We shall consider the branching problem of \smash{$\lmod{g}{q}{S}(F)$}.
To apply Corollary~\ref{cor:IrreducibleVermaQuasiAbelian} to \smash{$\lmod{g}{q}{S}(F)$}, we need an additional assumption.
We assume that $K'$ acts on $K/L_K$ transitively and $\lie{q}$ is quasi-abelian with respect to $\lie{g}'$.

\begin{Theorem}\label{thm:DerivedFunctorModuleQuasiAbelian}
	Let $F$ be a finite-dimensional irreducible $(\lie{l}, L_K)$-module in the good range.
	\begin{enumerate}\itemsep=0pt
		\item[$1.$] \smash{$\lmod{g}{q}{S}(F)|_{(\lie{g}', K')}$} is decomposed into a direct sum of irreducible modules of the form \smash{$\lmod{g'}{q'}{S}(F')$} with finite-dimensional irreducible \smash{$\bigl(\lie{l}', L'_K\bigr)$}-module $F'$ in the good range.
		\item[$2.$] \smash{$\ind{g}{q}(F)|_{\lie{g'}, L'_K}$} is decomposed into a direct sum of irreducible modules of the form \smash{$\ind{g'}{q'}(F')$} with finite-dimensional irreducible \smash{$\bigl(\lie{l}', L'_K\bigr)$}-module $F'$ in the good range.
		\item[$3.$] For any finite-dimensional irreducible \smash{$\bigl(\lie{l'}, L'_K\bigr)$}-module $F'$ in the good range, \smash{$\Dzuck{K'}{L'_K}{S}$} induces an isomorphism
		\begin{align*}
			\Hom_{\lie{g}', L'_K}\bigl(\ind{g'}{q'}(F'), \ind{g}{q}(F)|_{\lie{g'}, L'_K}\bigr) \xrightarrow{\simeq} \Hom_{\lie{g}', K'}\bigl(\lmod{g'}{q'}{S}(F'), \lmod{g}{q}{S}(F)|_{\lie{g'}, K'}\bigr)
		\end{align*}
		of \smash{$\univ{g}^{G'}$}-modules, and
		each \smash{$\univ{g}^{G'}$}-module \smash{$\Hom_{\lie{g}', K'}\bigl(\lmod{g'}{q'}{S}(F'), \lmod{g}{q}{S}(F)|_{\lie{g'}, K'}\bigr)$} is ir\-re\-du\-cible or zero.
	\end{enumerate}
\end{Theorem}

\begin{Remark}
	The abstract branching law of \smash{$\lmod{g}{q}{S}(F)|_{(\lie{g}', K')}$} is known in~\cite[Corollaries 5.7 and 5.8]{Os24} in the Grothendieck group level.
	His result for $A_{\lie{q}}(\lambda)$ (i.e., $F$ is a character) is done for~$\lambda$ in the weakly fair range.
	His results are proved under weaker assumptions than ours.
\end{Remark}

\begin{proof}[Proof of Theorem~\ref{thm:DerivedFunctorModuleQuasiAbelian}]
	The assertion~2 has been proved in Theorem~\ref{thm:QuasiAbelianGoodRange}.
	By Lemma~\ref{lem:ZuckermanKtoH}, there exists an isomorphism \smash{$\lmod{g}{q}{S}(F) \simeq \Dzuck{K'}{L'_K}{S}\bigl(\ind{g}{q}(F)\bigr)$} of $(\lie{g}, K')$-modules.
	Then we have
	\begin{align*}
		\ind{g}{q}(F)|_{(\lie{g}', L'_K), \univ{g}^{G'}} \simeq \bigoplus_{F'} \ind{g'}{q'}(F') \otimes \Hom_{\lie{g}', L'_K}\bigl(\ind{g'}{q'}(F'), \ind{g}{q}(F)\bigr),
	\end{align*}
	and hence
	\begin{align*}
		\lmod{g}{q}{S}(F)|_{(\lie{g}', K'), \univ{g}^{G'}} \simeq \bigoplus_{F'} \lmod{g'}{q'}{S}(F') \otimes \Hom_{\lie{g}', L'_K}\bigl(\ind{g'}{q'}(F'), \ind{g}{q}(F)\bigr).
	\end{align*}
	The sum is taken over all finite-dimensional irreducible \smash{$\big(\lie{l}', L'_K\big)$}-modules in the good range.
	
	By Fact~\ref{fact:ZuckermanDerivedFunctorModule}, each \smash{$\lmod{g'}{q'}{S}(F')$} is non-zero and irreducible.
	By Proposition~\ref{prop:UniqueZuckermanModule}, \smash{$\lmod{g'}{q'}{S}(F'')$} is not~isomorphic to \smash{$\lmod{g'}{q'}{S}(F')$} for $F''\not \simeq F'$ in the good range.
	This implies
	\begin{align*}
		\Hom_{\lie{g}', L'_K}\bigl(\ind{g'}{q'}(F'), \ind{g}{q}(F)|_{\lie{g'}, L'_K}\bigr) \simeq \Hom_{\lie{g}', K'}\bigl(\lmod{g'}{q'}{S}(F'), \lmod{g}{q}{S}(F)|_{\lie{g'}, K'}\bigr).
	\end{align*}
	This isomorphism is induced by the functor $\Dzuck{K'}{L'_K}{S}$.
	The irreducibility of the $\univ{g}^{G'}$-module follows from Corollary~\ref{cor:IrreducibleVermaQuasiAbelian}.
	We have proved the theorem.
\end{proof}

We shall state an important case in which the assumption of Theorem~\ref{thm:DerivedFunctorModuleQuasiAbelian} is fulfilled.
We assume that there exist two non-trivial ideals $\lie{k}_1$, $\lie{k}_2$ of $\lie{k}$ such that
\begin{align*}
	\lie{k} = \lie{k}_1 \oplus \lie{k}_2, \qquad \lie{k'} = \lie{k}_1 \oplus (\lie{k'}\cap \lie{k}_2), \qquad H \in \lie{k}_1.
\end{align*}
Set $\lie{k}'_2\coloneq \lie{k'}\cap \lie{k}_2$.
Then we have $\lie{k}_2 \subset \lie{l}_K$ and $\lie{u}\cap \lie{k}\subset \lie{k}_1$.
Let $K_i$, $i = 1,2$, denote the analytic subgroup of $K$ with the Lie algebra $\lie{k}_i$.
It is clear that $K'$ acts on $K/L_K$ transitively by $K_2 \subset L_K$ and $K_1\subset K'$.

For example, $(\lie{g}, \lie{g'}) = (\lieU(p, q), \lieU(p, q'))$, $q' < q$, satisfies the assumption.
In the case, ${\lie{k}_1 \simeq \lieU(p)}$, $\lie{k}_2 \simeq \lieU(q)$ and $\lie{k}'_2 \simeq \lieU(q')$.

We shall give a criterion for the condition that $\lie{q}$ is quasi-abelian with respect to $\lie{g'}$.

\begin{Proposition}\label{prop:QuasiAbelianZuckerman}
	If $\lie{q}$ is quasi-abelian with respect to $\lie{k}_1$, then $\lie{q}$ is quasi-abelian with respect to~$\lie{g'}$.
\end{Proposition}

\begin{proof}
	Replacing $\lie{t'}$, we may assume that $\lie{t'}$ contains a Cartan subalgebra of $\lie{k}_1$.
	Set $\lie{t}'_{c}\coloneq\lie{t}'\cap \lie{k}$.

	Recall that we have assumed $H \in \lie{k}_1 \subset \lie{g'}$.
	Set $\lie{u}''\coloneq \lie{u}\cap (\lie{g'})^{\perp}$.
	Then we have $\lie{u} = \lie{u}'' \oplus \lie{u}'$.
	Since $\lie{\overline{u}}\cap \lie{k}\subset \lie{k}_1\subset \lie{g'}$ and \smash{$\lie{u} = (\lie{u}\cap \lie{k}) \oplus \bigl(\lie{u} \cap \lie{k}^\perp\bigr)$}, we have $\lie{u}'' \subset \lie{k}^{\perp}$.

	Assume that $\lie{q}$ is not quasi-abelian with respect to $\lie{g}'$.
	Then there exist $\alpha \in \Delta(\lie{u'}, \lie{t'})$ and $\beta \in \Delta(\lie{u''}, \lie{t'})$
	such that $(\alpha, \beta) < 0$.
	Since $(\lie{g'})^{\perp}$ is a $\lie{g}'$-module, this implies $\bigl[\lie{u}'_{\alpha}, \lie{u}''_{\beta}\bigr] \neq 0$.
	Since $\lie{q}'$ is $\theta$-stable, there are two possibilities:
	\begin{enumerate}\itemsep=0pt
		\item[(1)] $\lie{u}'_{\alpha} \subset \lie{k}^{\perp}$,
		\item[(2)] $\bigl(\lie{u}'_{\alpha}\oplus \lie{u}'_{\theta(\alpha)}\bigr) \cap \lie{k}'\neq 0$.
	\end{enumerate}

	Assume (1).
	Using $[\lie{u}', \lie{u}'']\subset \lie{u}''$, $\lie{u}'' \subset \lie{k}^{\perp}$ and $\bigl[\lie{k}^{\perp}, \lie{k}^{\perp}\bigr] \subset \lie{k}$,
	we have
	\begin{align*}
	\bigl[\lie{u}'_{\alpha}, \lie{u}''_{\beta}\bigr] &\subset \lie{u}'' \cap \lie{k}
	\subset \lie{k}^{\perp} \cap \lie{k} =0.
	\end{align*}
	This contradicts $\bigl[\lie{u}'_{\alpha}, \lie{u}''_{\beta}\bigr] \neq 0$.

	Assume (2).
	By $\lie{u}''_{\beta} \subset \lie{k}^{\perp}$, $\beta$ is $\theta$-invariant.
	Hence we have
	\begin{align*}
	\bigl(\alpha|_{\lie{t}'_{c}}, \beta|_{\lie{t}'_{c}}\bigr) = (\alpha, (\beta + \theta(\beta)) / 2) = (\alpha, \beta) < 0. 
	\end{align*}
	By the assumption (1), $\alpha|_{\lie{t}'_{c}}$ belongs to $\Delta(\lie{u}'\cap \lie{k'}, \lie{t}'_{c}) = \Delta(\lie{u}\cap \lie{k}_1, \lie{t}'_c)$.
	Therefore, this contradicts the assumption that $\lie{q}$ is quasi-abelian with respect to $\lie{k}_1$.
	We have proved the assertion.
\end{proof}

\subsection{Discrete series representation}

We assumed in Theorem~\ref{thm:DerivedFunctorModuleQuasiAbelian} that $\lie{q}$ is quasi-abelian with respect to $\lie{g}'$ and $K'$ acts on $K/L_K$ transitively.
Note that there exist many tuples $(\lie{g}, K, \lie{g'}, K', \lie{q})$ such that $\lmod{g}{q}{S}(F)|_{(\lie{g}', K')}$ is discretely decomposable and $K'$ does not act on $K/L_K$ transitively.
See the classification in~\cite{KoOs12}.
On the other hand, many discretely decomposable restrictions of discrete series representations satisfy the assumption.

Let $G_\RR$ be a connected non-compact simple Lie group with finite center and a Cartan involution $\theta$,
and $G'_{\RR}$ be a $\theta$-stable connected non-compact reductive subgroup of $G_{\RR}$.
Set $K_{\RR}\coloneq G^{\theta}_{\RR}$ and \smash{$K'_{\RR}\coloneq \bigl(G'_{\RR}\bigr)^{\theta}$}.
We assume that \smash{$\bigl(G_{\RR}, G'_{\RR}\bigr)$} is a symmetric pair and $\rank(\lie{g}) = \rank(\lie{k})$.
Then~$G_\RR$ has discrete series representations.

Fix a maximal torus $T_{\RR}$ of $K_{\RR}$ satisfying that $T_{\RR} \cap K'_{\RR}$ is a maximal torus of $K'_{\RR}$.
Then $T_{\RR}$ is a fundamental Cartan subgroup of $G_{\RR}$.
Take a $\theta$-stable Borel subalgebra $\lie{b}$ of $\lie{g}$ containing $\lie{t}$.
Let $\lie{n}$ denote the nilpotent radical of $\lie{b}$.
Set $S'\coloneq \dim(\lie{n})$.

For a unitary character $\CC_{\lambda}$ of $T_{\RR}$ in the good range with respect to $\lie{b}$,
$\lmod{g}{b}{S'}(\CC_{\lambda})$ is unitarizable and its completion is a discrete series representation of $G_{\RR}$.
By the classification (\cite{KoOs12}) of discretely decomposable $A_{\lie{q}}(\lambda)$'s,
we can see the following fact.
See also~\cite{DuVa10} and~\cite[Section 8]{Os24} for the details.

\begin{Fact}\label{fact:DiscDecompDiscSeries}
	Let $\CC_{\lambda}$ be a unitary character of $T_{\RR}$ in the good range with respect to $\lie{b}$.
	Suppose that $\lmod{g}{b}{S'}(\CC_{\lambda})|_{(\lie{g}', K')}$ is discretely decomposable.
	Then there exists an element $H \in \lie{g}' \cap \sqrt{-1}\lie{t}_\RR$ such that
	the corresponding $\theta$-stable parabolic subalgebra $\lie{q}\coloneq \lie{q}(H) = \lie{l}\oplus \lie{u}$
	satisfies the following conditions:
	\begin{enumerate}\itemsep=0pt
		\item[$1.$] $\lie{q} \supset \lie{b}$.
		\item[$2.$] $\lie{l}\subset \lie{k}$.
		\item[$3.$] $\lie{q}$ is quasi-abelian with respect to $\lie{k}$.
		\item[$4.$] There exists an ideal $\lie{k}_1$ of $\lie{k}$ such that
		$\lie{u} \cap \lie{k} \subset \lie{k}_1\subset \lie{g}'$ and $H \in \lie{k}_1$.
		\item[$5.$]
There exists a finite-dimensional irreducible $(\lie{l}, L_K)$-module $F$ in the good range such that
		\begin{align*}
			\lmod{g}{b}{S'}(\CC_{\lambda}) \simeq \lmod{g}{q}{S}(F),
		\end{align*}
		where $S=\dim(\lie{u})$ and $L_K$ is the centralizer of $H$ in $K$.
	\end{enumerate}
\end{Fact}

More precisely, $\lie{u}$ in Fact~\ref{fact:DiscDecompDiscSeries} is at most 2-step nilpotent and satisfies $[\lie{u}, \lie{u}]\subset \lie{k}$.
In this case, $\lmod{g}{q}{S}(F)$ is called a small discrete series representation in~\cite{GrWa00}.
The ideal $\lie{k}_1$ in Fact~\ref{fact:DiscDecompDiscSeries} is generated by $\lie{k}\cap \lie{u}$.
In~\cite[Corollary~1]{DuVa10}, for a symmetric subgroup $G'_\RR\subset G_\RR$, the restriction $\lmod{g}{b}{S'}(\CC_{\lambda})|_{(\lie{g}', K')}$ is discretely decomposable (or equivalently $\lie{g'}$-admissible) if and only if $\lie{k}_1\subset \lie{k}'$.
Note that we have proved the if part in Theorem~\ref{thm:DerivedFunctorModuleQuasiAbelian}.

Assume that the assumption of Fact~\ref{fact:DiscDecompDiscSeries}, and take $\lie{q}$ in the fact.
Replacing $\lie{b}$ with its conjugate by an inner automorphism of $\lie{l}$, we may assume that $\lie{b'}\coloneq \lie{b}\cap \lie{g'}$ is a $\theta$-stable Borel subalgebra of $\lie{g'}$.
Set $S''\coloneq \dim(\lie{n} \cap \lie{g'})$.

\begin{Theorem}\label{thm:DiscreteSeries}
	Let $\CC_{\lambda}$ be a unitary character of $T_{\RR}$ in the good range with respect to $\lie{b}$.
	Then $\lmod{g}{b}{S'}(\CC_{\lambda})|_{\lie{g'}, K'}$ is decomposed into a direct sum of irreducible modules of the form \smash{$\lmod{g'}{b'}{S''}(\CC_{\lambda'})$} with~$\CC_{\lambda'}$ in the good range.
	Moreover, \smash{$\univ{g}^{G'}$} acts on \smash{$\Hom_{\lie{g'},K'}[$.$]\bigl(V', \lmod{g}{b}{S'}(\CC_{\lambda})|_{\lie{g'}, K'}[$.$]\bigr)$} irreducibly for any irreducible submodule $V'$ in \smash{$\lmod{g}{b}{S'}(\CC_{\lambda})|_{\lie{g'}, K'}$}.
\end{Theorem}

\begin{Remark}
	The branching law of $\lmod{g}{b}{S'}(\CC_{\lambda})|_{\lie{g'}, K'}$ is described by the branching law of a~generalized Verma module $\ind{g}{q}(F)|_{\lie{g'}, L'_K}$ as we have seen in Theorem~\ref{thm:DerivedFunctorModuleQuasiAbelian}.
\end{Remark}

\begin{Remark}
	The branching law of $\lmod{g}{b}{S'}(\CC_{\lambda})|_{\lie{g'}, K'}$ is computed by several researchers.
	In~\cite{GrWa00} and~\cite{DuVa10}, the multiplicities are given by an alternating sum like the Blattner's
	formula.
	In~\cite{Va16} and~\cite{OrVa24}, the branching law is reduced to the $K$-type formula of another discrete series representation of $G''_\RR$, where \smash{$\bigl(G_\RR, G''_\RR\bigr)$} is the associated symmetric pair of \smash{$\bigl(G_\RR, G'_\RR\bigr)$}.
	In~\cite[Section~8]{Os24}, the branching law is computed explicitly.
\end{Remark}

\begin{proof}[Proof Theorem~\ref{thm:DiscreteSeries}]
	By Proposition~\ref{prop:QuasiAbelianZuckerman}, $\lie{q}$ is quasi-abelian with respect to $\lie{g'}$.
	The assertion is reduced to Theorem~\ref{thm:DerivedFunctorModuleQuasiAbelian} by Fact~\ref{fact:DiscDecompDiscSeries}.
	Note that for any finite-dimensional irreducible $(\lie{l}, L_K)$-module~$F'$ in the good range, \smash{$\lmod{g'}{q'}{S}(F')$} is isomorphic to some \smash{$\lmod{g'}{b'}{S''}(\CC_{\lambda'})$} with a unitary character $\CC_{\lambda'}$ of $T'_\RR$ in the good range.
	In fact, by induction in stages (see~\cite[Corollary~11.86]{KnVo95_cohomological_induction}) and the vanishing theorem, we have
\[
\lmod{g'}{q'}{S}\bigl(\lmod{l'}{b'\cap l'}{S''-S}(\CC_{\lambda'})\bigr) \simeq \lmod{g'}{b'}{S''}(\CC_{\lambda'}),
\]
	taking $\CC_{\lambda'}$ to satisfy \smash{$\lmod{l'}{b'\cap l'}{S''-S}(\CC_{\lambda'}) \simeq F'$}.
\end{proof}

\subsection{Polynomial identity}

One of the motivations of our study of $\univ{g}^{G'}$-modules is to relate the multiplicities to algebraic properties of $\univ{g}^{G'}$.
We shall consider polynomial identities as one of the properties.

We recall the notion of PI degree, which estimates non-commutativity of algebras.
We refer the reader to~\cite[Section~13]{McRo01_noncommutative}.

\begin{Definition}
	Let $\alg{A}$ be a (unital associative) $\CC$-algebra.
	For a $\ZZ$-coefficient non-commutative polynomial $f$ with $n$ indeterminates, we say that $f$ is a \emph{polynomial identity} of $\alg{A}$ if
	\begin{align*}
		f(X_1, X_2, \ldots, X_n) = 0, \qquad \forall X_i \in \alg{A}.
	\end{align*}
	We denote by $\PI(\alg{A})$ the set of all polynomial identities of $\alg{A}$, and set
	\begin{align*}
		\PIdeg(\alg{A}) \coloneq \sup\set{n \in \NN \mid \PI(\alg{A})\subset \PI(M_n(\CC))},
	\end{align*}
	which is called the \emph{PI degree} of $\alg{A}$.
\end{Definition}

It is clear that if there exists a surjective homomorphism $\alg{A}\rightarrow M_n(\CC)$, then we have $\PIdeg(\alg{A}) \geq n$.
Roughly speaking, $\PIdeg(\alg{A})$ is the maximum dimension of irreducible $\alg{A}$-modules.
This is true under a mild assumption (see Proposition~\ref{prop:PIdegMaximumDim}).

Let $s_n$ denote the $\ZZ$-coefficient non-commutative polynomial with $n$ indeterminates defined~by
\begin{align*}
	s_n(X_1, X_2, \ldots, X_n) \coloneq \sum_{w \in \Sn_n} \sgn(w)X_{w(1)}X_{w(2)}\cdots X_{w(n)},
\end{align*}
where $\Sn_n$ is the symmetric group of degree $n$ and $\sgn$ is the signature character of $\Sn_n$.

\begin{Fact}[{\cite[Proposition~3.2 and Theorem~3.3]{McRo01_noncommutative}}]
	Fix an integer $n > 0$.
	\begin{enumerate}
		\item[$1.$] $($Amitsur--Levitzki$)$ $s_m \in \PI(M_n(\CC))$ for any $2n\leq m$.
		\item[$2.$] Conversely, $s_m \not\in \PI(M_n(\CC))$ for any $m < 2n$.
	\end{enumerate}
	In particular, we have $\PIdeg(M_n(\CC)) = n$.
\end{Fact}

Let $\alg{A}$ be a $\CC$-algebra.
Assume that $\alg{A}$ has at most countable dimension.

\begin{Proposition}\label{prop:PIdegMaximumDim}\qquad
	\begin{enumerate}\itemsep=0pt
		\item[$1.$] Let $\set{M_i}_{i \in I}$ be a family of irreducible $\alg{A}$-modules.
		Suppose that $\bigcap_i \Ann_{\alg{A}}(M_i)$ is zero.
		Then one has
		\begin{align*}
			\PIdeg(\alg{A}) = \sup\set{\dim(M_i)\mid i \in I}.
		\end{align*}
		\item[$2.$] If $\alg{A}$ is semiprimitive, then $\PIdeg(\alg{A})$ is the maximum dimension of irreducible $\alg{A}$-modules.
	\end{enumerate}
\end{Proposition}

\begin{proof}
	The second assertion easily follows from the first assertion.
	We shall show the first assertion.

	Assume that there exists $i \in I$ such that $\dim(M_i) = \infty$.
	Here $\infty$ means the cardinality $|\NN|$.
	By the Jacobson density theorem, $\alg{A}/\Ann_{\alg{A}}(M_i)$ contains a subalgebra isomorphic to $M_n(\CC)$ for any~${n \in \NN}$.
	This implies $\PIdeg(\alg{A}) = \infty = \sup\set{\dim(M_i)\mid i \in I}$.
	Hence if ${\sup\set{\dim(M_i)\mid i \in I} = \infty}$, then the assertion holds.

	Assume that $\sup\set{\dim(M_i)\mid i \in I} < \infty$.
	Then $\alg{A}/\Ann_{\alg{A}}(M_i)$ is isomorphic to $M_{\dim(M_i)}(\CC)$ for each $i \in I$.
	This implies $\PIdeg(\alg{A}) \geq \sup\set{\dim(M_i)\mid i \in I}$.
	By assumption, there is an injective homomorphism $\alg{A} \hookrightarrow \prod_{i \in I} \End_{\CC}(M_i)$.
	Hence we obtain the converse inequality $\PIdeg(\alg{A}) \leq \sup\set{\dim(M_i)\mid i \in I}$.
	We have shown the proposition.
\end{proof}

Let $(\lie{g}, G')$ be a pair.
The following proposition is an easy consequence of Proposition~\ref{prop:PIdegMaximumDim}.

\begin{Proposition}\label{prop:PIdegMultiplicity}
	Let $V$ be a $\lie{g}$-module such that $V|_{\lie{g'}}$ is completely reducible.
	Write $\Irr\bigl(V|_{\lie{g'}}\bigr)$ for the set of all isomorphism classes of irreducible $\lie{g'}$-submodules in $V|_{\lie{g'}}$.
	Assume that the $\univ{g}^{G'}$-module $\Hom_{\lie{g'}}\bigl(V', V|_{\lie{g'}}\bigr)$ is irreducible for any $V' \in \Irr\bigl(V|_{\lie{g'}}\bigr)$.
	Then one has
	\begin{align*}
		\sup_{V' \in \Irr(V|_{\lie{g'}})} \dim \Hom_{\lie{g'}}\bigl(V', V|_{\lie{g'}}\bigr) = \PIdeg\bigl(\bigl(\univ{g}/\Ann_{\univ{g}}(V)\bigr)^{G'}\bigr).
	\end{align*}
\end{Proposition}

\begin{proof}
	The assertion follows from Proposition~\ref{prop:PIdegMaximumDim} for
\[
\alg{A} = \bigl(\univ{g}/\Ann_{\univ{g}}(V)\bigr)^{G'}
	\qquad\text{and}\qquad
\set{M_i}_{i \in I} = \set{\Hom_{\lie{g'}}\bigl(V', V|_{\lie{g'}}\bigr)}_{V' \in \Irr(V|_{\lie{g'}})}.
\tag*{\qed}
\]
\renewcommand{\qed}{}
\end{proof}

Combining Corollary~\ref{cor:IrreducibleVermaQuasiAbelian}, Theorem~\ref{thm:DerivedFunctorModuleQuasiAbelian} and Proposition~\ref{prop:PIdegMultiplicity}, we obtain the following.

\begin{Theorem}\label{thm:PIdegMultiplicity}
	Let $V$ be a generalized Verma module $\ind{g}{q}(F)$ in Theorem~$\ref{thm:QuasiAbelianGoodRange}$ or a cohomologically induced module $\lmod{g}{q}{S}(F)$ in Theorem~$\ref{thm:DerivedFunctorModuleQuasiAbelian}$.
	Then one has
	\begin{align*}
		\sup_{V' \in \Irr(V|_{\lie{g'}})} \dim \Hom_{\lie{g'}}\bigl(V', V|_{\lie{g'}}\bigr) = \PIdeg\bigl(\bigl(\univ{g}/\Ann_{\univ{g}}(V)\bigr)^{G'}\bigr).
	\end{align*}
	In particular, $V|_{\lie{g'}}$ is multiplicity-free if and only if \smash{$\bigl(\univ{g}/\Ann_{\univ{g}}(V)\bigr)^{G'}$} is commutative.
\end{Theorem}

Under the assumption of Theorem~\ref{thm:PIdegMultiplicity}, the algebra \smash{$\alg{A}\coloneq \bigl(\univ{g}/\Ann_{\univ{g}}(V)\bigr)^{G'}$} is semiprimitive.
Then $\PIdeg(\alg{A}) = 1$ if and only if any irreducible $\alg{A}$-module is one-dimensional, that is, $\alg{A}$ is commutative (see Proposition~\ref{prop:PIdegMaximumDim}).
Note that $X_1X_2 - X_2X_1$ is a polynomial identity of~$M_1(\CC)$.

If $V$ is an irreducible $\lie{g}$-module and $V|_{\lie{g'}}$ is locally finite and completely reducible,
then the assumption of Proposition~\ref{prop:PIdegMultiplicity} is fulfilled.
In this case, an analogue of Theorem~\ref{thm:PIdegMultiplicity} appeared in the proof of~\cite[Theorem~4.3]{PeSe12}.

In general, $\Hom_{\lie{g'}}\bigl(V', V|_{\lie{g'}}\bigr)$ \big(or $\Hom_{\lie{g'}}\bigl(V|_{\lie{g'}}, V'\bigr)$\big) may not be irreducible as we have seen in Example~\ref{ex:NonIrreducible}.
We gave in~\cite{Ki20} a similar uniform estimate of multiplicities by $\PIdeg$,
and a~relation between the finiteness of $\PIdeg$ and coisotropic actions on nilpotent orbits.

\subsection*{Acknowledgement}

This work is based on a part of my Ph.D.\ Thesis under the supervision of Professor Toshiyuki Kobayashi.
I am deeply grateful to him for his continuous support, invaluable advice and constant encouragement.
I thanks the anonymous referees for valuable comments that improved the exposition and clarified the results.
This work was supported by JSPS KAKENHI Grant Number JP23K12963.

\pdfbookmark[1]{References}{ref}
\LastPageEnding

\end{document}